# A large-deviations principle for all the components in a sparse inhomogeneous random graph

Luisa Andreis[1] · Wolfgang König[2] · Heide Langhammer[3] ·
Robert I. A. Patterson[3]

## Abstract

We study an inhomogeneous sparse random graph, $\mathcal{G}_N$, on $[N] = \{1, \ldots, N\}$ as introduced in a seminal paper by Bollobás et al. (Random Struct Algorithms 31(1):3–122, 2007): vertices have a type (here in a compact metric space $\mathcal{S}$), and edges between different vertices occur randomly and independently over all vertex pairs, with a probability depending on the two vertex types. In the limit $N \to \infty$, we consider the sparse regime, where the average degree is $O(1)$. We prove a large-deviations principle with explicit rate function for the statistics of the collection of all the connected components, registered according to their vertex type sets, and distinguished according to being microscopic (of finite size) or macroscopic (of size $\asymp N$). In doing so, we derive explicit logarithmic asymptotics for the probability that $\mathcal{G}_N$ is connected. We present a full analysis of the rate function including its minimizers. From this analysis we deduce a number of limit laws, conditional and unconditional, which provide comprehensive information about all the microscopic and macroscopic components of $\mathcal{G}_N$. In particular, we recover the criterion for the existence of the phase transition given in Bollobás et al. (2007).

✉ Heide Langhammer
  langhammer@wias-berlin.de

  Luisa Andreis
  luisa.andreis@polimi.it

  Wolfgang König
  koenig@wias-berlin.de

  Robert I. A. Patterson
  patterson@wias-berlin.de

[1] Dipartimento di Matematica, Politecnico di Milano, Piazza Leonardo da Vinci, 32, 20133 Milano, Italy

[2] TU Berlin and WIAS, Mohrenstraße 39, 10117 Berlin, Germany

[3] WIAS, Mohrenstraße 39, 10117 Berlin, Germany

## 1 Introduction

In this paper, we study the *inhomogeneous random graph* model as introduced in the seminal paper [7], that is, an Erdős–Rényi graph whose vertices have types. We consider the limit of a large number of vertices and concentrate on the *sparse setting*, where each vertex has a number of edges that is of order one. This setting is famous for the emergence of a *giant cluster*. This phase transition was detected and characterized in [7] with the help of a branching process. We consider the case in which the type set is a compact metric space, but our analysis builds on the proof for the type set being any finite set.

In the present paper, we analyze the model from the view point of *large deviations* in a detailed way. We go beyond existing results by (1) considering the joint statistics of all the clusters, both microscopic and macroscopic, (2) registering the types within the clusters (not only their sizes), and (3) giving a joint large-deviations principle (LDP) for all this information. In particular, we recover the limiting quantities and the resulting phase transition in great detail, giving a lot of additional information. Our main results are Theorems 1.1 and 3.1 (the LDPs for the type set being a compact metric space and a finite type set, respectively) and Theorems 2.3 and 2.1, where we deduce consequences for the phase transition. A building block for our study is Theorem 3.6, which gives explicit logarithmic asymptotics for the probability of a macroscopic subgraph being connected and which is of independent interest.

The remainder of this section is organized as follows. In Sect. 1.1 we introduce the inhomogeneous random graph, in Sect. 1.2 we present our first main result, the large-deviations principle, and in Sect. 1.3 we give an interpretation of the result. Asymptotic results on the connectivity of graphs are highlighted in Sect. 1.4. In Sect. 1.5 we compare our results with the existing literature.

The structure of the rest of the paper is as follows. Our second main result concerns consequences of the large-deviations principle for the limiting behaviour of the model, i.e., conditional and unconditional laws of large numbers. This relies on explicit variational analysis of the rate function of Theorem 1.1 and it is explained in Sect. 2, together with the giant-cluster phase transition and a comparison to the results of [7]. Additionally, we explain the deep connection with an important inhomogeneous coagulation process and derive a solution to a spatial version of the Flory equation. In Sect. 3 we present the proof for our main result, the LDP, in the case of a finite type set. One key ingredient that we use are asymptotics for connection probabilities for different cluster sizes. They require more extensive proofs and are therefore contained separately in Sect. 4. In Sect. 3 we introduce the notation and the results in the framework of a finite type set, where many objects of interest have a simpler representation.

Reading only this section is a gentle, yet complete, introduction to our results and it is certainly suitable to readers with limited time or only interested in graphs with finite type sets. In Sect. 5, we derive the LDP in the general setting, i.e., for compact metric type spaces via a discrete approximation in the spirit of the Dawson–Gärtner theorem. In Sect. 6 we derive a full characterization of the minimizers of the microscopic part of the rate function, and in Sect. 7 we analyze all the other parts of the rate function.

### 1.1 Inhomogeneous random graphs

We are going to define the random graph model that we study in this paper. It is called an *inhomogeneous random graph* in [7], while at full length it is sometimes named *inhomogeneous random Erdős–Rényi graph*.

Let $N \in \mathbb{N}$; we consider a random graph on the vertex set $[N] = \{1, \ldots, N\}$ and fix a non-void set $\mathcal{S}$, the *type set*. We take a type vector $\mathbf{x} = \mathbf{x}^{(N)} = (x_1, \ldots, x_N) \in \mathcal{S}^N$ and interpret $x_i$ as the type of vertex $i$ for any $i \in [N]$. The edges of this graph are undirected and randomly drawn; self and multiple edges are excluded. The $\binom{N}{2}$ possible edges are sampled independently. The probability to draw an edge between two vertices with types $r$ and $s$, is called $p_{r,s}$; this defines a map $p \colon \mathcal{S} \times \mathcal{S} \to [0, 1]$. The resulting random graph is denoted $\mathcal{G}(N, \mathbf{x}, p)$ and is called the inhomogeneous random graph on $[N]$ with type vector $\mathbf{x}$ and function of probabilities $p$.

In this paper, we are interested in the limit as $N \to \infty$ in the *sparse case*, i.e., in the case where the number of edges per vertex is of finite order. This is the case if the probabilities $p_{r,s}$ are of order $1/N$. Actually, we now impose that they are given by

$$p_{r,s} = 1 \wedge \frac{1}{N} \kappa_N(r, s), \qquad r, s \in \mathcal{S},$$

where $\kappa_N \colon \mathcal{S} \times \mathcal{S} \to [0, \infty)$ is a symmetric non-negative bounded function, called a *kernel*. Without loss of generality, we are from now on assuming that $\frac{1}{N} \kappa_N \leq 1$, and hence $\kappa_N(r, s)/N$ is a probability for any $N \in \mathbb{N}$ and any $r, s \in \mathcal{S}$. We will study the graph $\mathcal{G}_N = \mathcal{G}(N, \mathbf{x}, \frac{1}{N} \kappa_N)$ in the limit $N \to \infty$.

We are interested in the structure of all the components[1] of $\mathcal{G}_N$, depending on the types, but not the indices of the vertices. Hence it is sufficient to consider the *empirical measure* $\mu_N = \frac{1}{N} \sum_{i=1}^N \delta_{x_i}$ of the type vector $\mathbf{x}$ and $N \mu_N(R)$ is the number of vertices with type in $R \subset \mathcal{S}$. We will assume that, as $N \to \infty$, $\mu_N$ converges weakly to a given probability measure $\mu$ on $\mathcal{S}$. One can conceive $\mathcal{G}_N$ as a graph on $\mathcal{S}$ where in each point $x \in \mathcal{S}$ there sit precisely $N \mu_N(\{x\}) \in \mathbb{N}$ vertices, which are all distinguished and labelled. Further, we will work under the assumption that $\kappa_N$ converges to a limiting kernel $\kappa \colon \mathcal{S} \times \mathcal{S} \to [0, \infty)$ that is continuous.

We denote by $\{\mathcal{C}_j\}_j$ the collection of all the vertex sets of the connected components of $\mathcal{G}_N$. This collection is a random decomposition of $[N]$. Since we are only interested in the statistics of these components, the labeling of the vertices in a component is irrelevant. Actually, we are even only interested in the statistics of the types of all the vertices, counted with multiplicity. To this end, we introduce the *type-registering*

---

[1] We use "cluster" and "component" synonymously.

*empirical-measure function*

$$\eta_{\mathbf{x}} \colon \mathcal{P}([N]) \to \mathcal{M}_{\mathbb{N}_0}(\mathcal{S}), \qquad \eta_{\mathbf{x}}(A) = \sum_{i \in A} \delta_{x_i}, \tag{1.1}$$

where $\mathcal{P}([N])$ is the set of subsets of $[N]$, and $\mathcal{M}_{\mathbb{N}_0}(\mathcal{S})$ is the space of finite measures on $\mathcal{S}$ with values in $\mathbb{N}_0$. We call every element of $\mathcal{M}_{\mathbb{N}_0}(\mathcal{S})$ a *type-configuration* and will denote it by $k$. In words, $\eta_{\mathbf{x}}(A)(R)$ is, for any set $R \subset \mathcal{S}$, the number of vertices in $A \subset [N]$ with type in $R$. In particular, $\eta_{\mathbf{x}}([N]) = N\mu_N$. We write $\mathcal{M}(\mathcal{X})$ for the set of finite measures on a set $\mathcal{X}$; the measurable structure on $\mathcal{X}$ will be clear from the context.

We will study the empirical measure of the collection $(\eta_{\mathbf{x}}(\mathcal{C}_j))_j$, i.e., the statistics of how many times a given type-configuration appears as the type-configuration of a component of $\mathcal{G}_N$. We will pay particular attention to the scale of the size of the component, more precisely, whether it is finite or it has a size $\asymp N$. We will call the first scale *microscopic* and the second scale *macroscopic* (macroscopic components are usually called *giant components*). Hence, the quantities of interest in our study are the *microscopic* and the *macroscopic empirical measures* of the type-configurations of the components, which we define as follows:

$$\mathrm{Mi}_N = \mathrm{Mi}_N(\mathbf{x}) = \frac{1}{N} \sum_j \delta_{\eta_{\mathbf{x}}(\mathcal{C}_j)} \quad \text{and} \quad \mathrm{Ma}_N = \mathrm{Ma}_N(\mathbf{x}) = \sum_j \delta_{\frac{1}{N}\eta_{\mathbf{x}}(\mathcal{C}_j)}. \tag{1.2}$$

It is clear that both $\mathrm{Mi}_N$ and $\mathrm{Ma}_N$ only depend on $\mathbf{x}$ through its empirical measure $\mu_N$. Both $\mathrm{Mi}_N$ and $\mathrm{Ma}_N$ are random measures on $\mathcal{M}(\mathcal{S})$, i.e., they are elements of $\mathcal{M}(\mathcal{M}(\mathcal{S}))$. More precisely, $\mathrm{Mi}_N$ is a measure on the set $\mathcal{M}_{\mathbb{N}_0}(\mathcal{S})$ of measures on $\mathcal{S}$ with values in $\mathbb{N}_0$, and $\mathrm{Ma}_N$ is a measure that takes values in $\mathbb{N}_0$ and is defined on the set $\mathcal{M}(\mathcal{S}) \setminus \{0\}$ of non-trivial measures on $\mathcal{S}$. Both $\mathrm{Mi}_N$ and $\mathrm{Ma}_N$ contain precisely the same information for fixed $N$, but in the limit $N \to \infty$, $\mathrm{Mi}_N$ asymptotically registers only the microscopic components and $\mathrm{Ma}_N$ only the macroscopic ones. Indeed, the non-microscopic components leave the state space $\mathcal{M}_{\mathbb{N}_0}(\mathcal{S})$ via the set of measures with unbounded total mass, and the non-macroscopic ones (with prefactor $1/N$) leave $\mathcal{M}(\mathcal{S})$ via the measures with vanishing total mass, i.e., via $\{0\}$. The topologies that we will introduce below reflect this effect; it lies at the heart of the phase transition of the emergence of a giant cluster, which we are also interested in here. On the other side, they are very natural, as they reduce to the pointwise respectively to the usual vague topology for a finite set $\mathcal{S}$.

The effect of a diverging or vanishing total mass can also be observed in terms of integrated versions of $\mathrm{Mi}_N$ respectively $\mathrm{Ma}_N$. Indeed, observe that for any measurable $R \subset \mathcal{S}$ and any $N \in \mathbb{N}$

$$\int_{\mathcal{M}_{\mathbb{N}_0}(\mathcal{S})} k(R) \, \mathrm{Mi}_N(\mathrm{d}k) = \frac{1}{N} \sum_j \int_{\mathcal{M}_{\mathbb{N}_0}(\mathcal{S})} k(R) \, \delta_{\eta_{\mathbf{x}}(\mathcal{C}_j)}(\mathrm{d}k)$$

$$= \frac{1}{N} \sum_j \eta_{\mathbf{x}}(\mathcal{C}_j)(R) = \mu_N(R). \tag{1.3}$$

According to our assumption that $\mu_N$ weakly converges towards $\mu$, the right-hand side of (1.3) converges to $\mu(R)$, as $N \to \infty$, if $R$ is a $\mu$-continuity set, i.e., if $\mu(\partial R) = 0$. However, the topology on the state space for $\mathrm{Mi}_N$ will be chosen as the vague one, and in this topology any accumulation point $\lambda$ of $(\mathrm{Mi}_N)_{N \in \mathbb{N}}$ satisfies *a priori* only $\int k(R) \lambda(\mathrm{d}k) \leq \mu(R)$, since the map $k \mapsto k(R)$ is unbounded in general. The same holds for $(\mathrm{Ma}_N)_{N \in \mathbb{N}}$.

To take into account the possibility of a loss of mass we introduce for any $\lambda \in \mathcal{M}(\mathcal{M}_{\mathbb{N}_0}(\mathcal{S}))$ and any $\alpha \in \mathcal{M}_{\mathbb{N}_0}(\mathcal{M}(\mathcal{S}) \setminus \{0\})$ the following measures on $\mathcal{S}$

$$c_\lambda(R) = \int_{\mathcal{M}_{\mathbb{N}_0}(\mathcal{S})} k(R) \lambda(\mathrm{d}k), \qquad R \subset \mathcal{S} \text{ measurable}, \tag{1.4}$$

$$c_\alpha(R) = \int_{\mathcal{M}(\mathcal{S}) \setminus \{0\}} y(R) \alpha(\mathrm{d}y), \qquad R \subset \mathcal{S} \text{ measurable}. \tag{1.5}$$

We call $c_\lambda$ and $c_\alpha$ the *integrated type-configurations* of $\lambda$, respectively of $\alpha$. If one sees $\lambda$ as a (not normalized) 'distribution' of finite point measures on $\mathcal{S}$, then $c_\lambda$ registers the 'expected' total mass of particles in a given subset of $\mathcal{S}$ that appear in this distribution $\lambda$; an analogous statement holds for $\alpha$. The total masses of $c_\lambda$ and $c_\alpha$ are equal to the integrals of $k \mapsto k(\mathcal{S})$ under $\lambda$ respectively under $\alpha$; they are $\leq 1$ for any accumulation point of $(\mathrm{Mi}_N)_{N \in \mathbb{N}}$ respectively of $(\mathrm{Ma}_N)_{N \in \mathbb{N}}$, according to the above.

The natural state space containing $\mathrm{Mi}_N$ for any $N \in \mathbb{N}$, is the set

$$\mathcal{L} = \{\lambda \in \mathcal{M}(\mathcal{M}_{\mathbb{N}_0}(\mathcal{S})) \colon c_\lambda \leq \mu \text{ or } c_\lambda \leq \mu_N \text{ for some } N \in \mathbb{N}, \lambda(\{0\}) = 0\}. \tag{1.6}$$

The condition $\lambda(\{0\}) = 0$ is clearly satisfied by any empirical measure $\mathrm{Mi}_N$, we could have identified it indeed as an element of $\mathcal{M}(\mathcal{M}_{\mathbb{N}_0}(\mathcal{S}) \setminus \{0\})$. However for later notational convenience we do not exclude $\{0\}$ but we add the constraint $\lambda(\{0\}) = 0$. Notice that with this constraint and the conditions on $c(\lambda)$ any $\lambda \in \mathcal{L}$ is a subprobability measure. Any $k$ in the support of $\lambda$ is a finite and non-zero point measure $k = \sum_i \delta_{z_i}$ with $z_i \in \mathcal{S}$ (with possible repetitions) and stands for the empirical measure of the types of a vertex set of a component, its type-configuration. Informally, the event $\{\mathrm{Mi}_N = \lambda\}$ is the event that, for every $k \in \mathcal{M}_{\mathbb{N}_0}(\mathcal{S})$, $\mathcal{G}_N$ has $N\lambda(\{k\})$ components with type-configuration $k$.

The natural state space containing $\mathrm{Ma}_N$ for any $N \in \mathbb{N}$, is the set

$$\mathcal{A} = \{\alpha \in \mathcal{M}_{\mathbb{N}_0}(\mathcal{M}(\mathcal{S}) \setminus \{0\}) \colon c_\alpha \leq \mu \text{ or } c_\alpha \leq \mu_N \text{ for some } N \in \mathbb{N}\}. \tag{1.7}$$

One can write $\alpha \in \mathcal{A}$ as a finite or at most countable point measure $\alpha = \sum_n \delta_{y_n}$ on measures $y_n$ on $\mathcal{S}$ (with possible repetitions). Note that $c_\alpha = \sum_n y_n$ and since $c_\alpha(\mathcal{S}) \leq 1$, the total masses $y_n(\mathcal{S})$ of the measures $y_n$ can accumulate only at zero. Another consequence is that every $\alpha \in \mathcal{A}$ is concentrated on the set of sub-probability measures. For each $y_n$ we interpret $N y_n$ as the type-configuration of a giant component. Informally, the event $\{\mathrm{Ma}_N = \alpha\}$ is the event that $\mathcal{G}_N$ has, for any $n$, a macroscopic component with type-configuration $N y_n$.

We will be working only with two particular choices of the type set: $\mathcal{S}$ as a finite set (equipped with the power set as topology and as sigma-field) and $\mathcal{S}$ equal to a compact metric space (equipped with the topology induced by the metric and the corresponding Borel sigma-field). We equip $\mathcal{M}(\mathcal{S})$ with the weak topology that is generated by all the test integrals against continuous bounded functions $\mathcal{S} \to \mathbb{R}$. However, on the sets $\mathcal{L}$ and $\mathcal{A}$, the appropriate topologies for our purposes are the *vague topologies*, the ones that are induced by all the test integrals against compactly supported continuous test functions $\mathcal{M}_{\mathbb{N}_0}(\mathcal{S}) \to \mathbb{R}$, respectively $\mathcal{M}(\mathcal{S}) \backslash \{0\} \to \mathbb{R}$. If $|\mathcal{S}| < \infty$, then $\mathcal{M}_{\mathbb{N}_0}(\mathcal{S})$ can be identified with $\mathbb{N}_0^{\mathcal{S}}$ and vague convergence in $\mathcal{M}(\mathcal{M}_{\mathbb{N}_0}(\mathcal{S}))$ is the same as pointwise convergence on $\mathcal{M}(\mathbb{N}_0^{\mathcal{S}})$. On $\mathcal{L} \times \mathcal{A}$ we use the product topology. We will show in Lemma 5.2 that both $\mathcal{L}$ and $\mathcal{A}$ are compact, hence also $\mathcal{L} \times \mathcal{A}$ is.

The convergence in this topology is the natural one that reflects the possible loss of mass that we are interested in in view of the phase transition. The crucial point is that mass of $\mathrm{Mi}_N$ can leak out only via the unboundedness of $k \mapsto k(\mathcal{S})$, i.e., via having larger and larger connected components, while mass of $\mathrm{Ma}_N$ can leak out due to the fact that $y \mapsto y(\mathcal{S})$ is not bounded away from zero. With other words, mass of $\mathrm{Ma}_N$ leaks out only via the zero measure, where every non-giant component leaves. See Sect. 5 for details. By the definitions of $\mathcal{L}$ respectively $\mathcal{A}$, the integrated type-configurations $c_\lambda$ and $c_\alpha$ for $\lambda \in \mathcal{L}$ and $\alpha \in \mathcal{A}$ are sub-probability measures, while for fixed $N \in \mathbb{N}$ both $c_{\mathrm{Mi}_N}$ and $c_{\mathrm{Ma}_N}$ are even probability measures. The total mass one of $c_{\mathrm{Mi}_N}$ can partially get lost and $c_{\mathrm{Ma}_N}$ may not lose all its mass in the limit $N \to \infty$. This is precisely the phase transition that we are after.

### 1.2 The large-deviations principle for the cluster statistics

In this section we formulate the main result of this paper. We assume that $\mathcal{S}$ is a compact metric space.

We need some notation. For any measure $\nu$ on $\mathcal{S}$ and any function $\kappa \colon \mathcal{S} \times \mathcal{S} \to [0, \infty)$, we write $\kappa \nu(r) = \int_{\mathcal{S}} \kappa(r, s) \nu(\mathrm{d}s)$. The total mass of a measure $\nu$ on a measure space $\mathcal{X}$ is denoted by $|\nu| = \nu(\mathcal{X})$. The relative entropy of two (possibly non-normalized) finite measures $\nu, \widetilde{\nu}$ on $\mathcal{X}$ is denoted by

$$\mathbb{H}(\nu | \widetilde{\nu}) = \begin{cases} |\widetilde{\nu}| - |\nu| + \int_{\mathcal{X}} \nu(\mathrm{d}x) \log \frac{\mathrm{d}\nu}{\mathrm{d}\widetilde{\nu}}(x), & \text{if } \frac{\mathrm{d}\nu}{\mathrm{d}\widetilde{\nu}} \text{ exists,} \\ +\infty & \text{otherwise.} \end{cases} \quad (1.8)$$

We write $\langle \nu, f \rangle$ for the integral of a function $f$ with respect to a measure $\nu$, and we write $f \nu$ for the measure that has the density $f$ with respect to a measure $\nu$.

An important reference measure is the distribution $\mathbb{Q}_\nu$ of a Poisson point process $\mathbb{X}$ on $\mathcal{S}$ with intensity measure $\nu \in \mathcal{M}(\mathcal{S})$, then $\mathbb{Q}_\nu$ is a measure on $\mathcal{M}_{\mathbb{N}_0}(\mathcal{S})$. Note that we do not assume that $\nu$ has a density, hence $\mathbb{X}$ is not necessarily simple.

We define a function $\tau$ by

$$\tau(k) := \sum_{T \in \mathcal{T}(k)} \prod_{\{i,j\} \in E(T)} \kappa(x_i, x_j), \quad k \in \mathcal{M}_{\mathbb{N}_0}(\mathcal{S}), \tag{1.9}$$

where $(x_1, \ldots, x_{|k|}) \in \mathcal{S}^{|k|}$ is any vector that is *compatible* with $k$, i.e., $k = \sum_{i=1}^{|k|} \delta_{x_i}$, and $\mathcal{T}(k)$ is the set of spanning trees on $[|k|]$. Notice that $\tau$ depends on $(x_1, \ldots, x_{|k|})$ only through $k$. We use the convention that $\mathcal{T}(0) = \emptyset$ and hence $\tau(0) = 0$, since the sum is empty. As we will see in Lemma 3.4, up to a factor of $N^{-|k|+1}$, $\tau(k)$ is equal to the large-$N$ asymptotics of the probability that the graph $\mathcal{G}(|k|, (x_1, \ldots, x_{|k|}), \frac{1}{N}\kappa)$, which can be seen as a subgraph of $\mathcal{G}(N, \mathbf{x}^{(N)}, \frac{1}{N}\kappa)$, is connected.

We say that $\kappa : \mathcal{S} \times \mathcal{S} \to [0, \infty)$ is *irreducible* with respect to a measure $\mu \in \mathcal{M}(\mathcal{S})$ if

$$A \subset \mathcal{S} \text{ and } \kappa = 0 \text{ a.e. on } A \times (\mathcal{S} \setminus A) \implies \mu(A) = 0 \text{ or } \mu(\mathcal{S} \setminus A) = 0. \tag{1.10}$$

Otherwise $\kappa$ is called reducible.

Here is the main result of this paper.

**Theorem 1.1** *(LDP for $(Mi_N, Ma_N)$). Fix a probability measure $\mu$ on a compact metric space $\mathcal{S}$ and a kernel $\kappa$ on $\mathcal{S} \times \mathcal{S}$ that is nonnegative, continuous and irreducible with respect to $\mu$. Assume that $\mathbf{x} = \mathbf{x}^{(N)} \in \mathcal{S}^N$ is such that its empirical measure $\mu_N$ converges weakly towards $\mu$ as $N \to \infty$. Assume that $\kappa_N$ is a nonnegative and continuous kernel for any $N \in \mathbb{N}$ such that $\kappa_N$ converges uniformly towards $\kappa$ as $N \to \infty$. Let $Mi_N$ and $Ma_N$ be, respectively, the microscopic and macroscopic empirical measure of the connected components of $\mathcal{G}_N = \mathcal{G}(N, \mathbf{x}^{(N)}, \frac{1}{N}\kappa_N)$, for any $N \in \mathbb{N}$, as defined in (1.2).*

*Then $(Mi_N, Ma_N)$ satisfies a large-deviations principle with speed $N$ and rate function $I$ defined by*

$$I(\lambda, \alpha) = \begin{cases} I_{Mi}(\lambda) + I_{Ma}(\alpha) + I_{Me}(\mu - c_\lambda - c_\alpha), & \text{if } c_\lambda + c_\alpha \leq \mu, \\ +\infty & \text{otherwise,} \end{cases}$$

*where, for $\lambda \in \mathcal{L}$, $\alpha \in \mathcal{A}$ and $\nu \in \mathcal{M}(\mathcal{S})$,*

$$I_{Mi}(\lambda) = \mathbb{H}(\lambda | \mathbb{Q}_\mu) - 1 - \langle \lambda, \log \tau \rangle + |c_\lambda| - |\lambda| + \frac{1}{2} \langle c_\lambda, \kappa \mu \rangle, \tag{1.11}$$

$$I_{Ma}(\alpha) = \int_{\mathcal{M}(\mathcal{S}) \setminus \{0\}} \alpha(\mathrm{d}y) \left[ \left\langle y, \log \frac{\mathrm{d}y}{(1 - e^{-\kappa y}) \mathrm{d}\mu} \right\rangle + \frac{1}{2} \langle y, \kappa(\mu - y) \rangle \right], \tag{1.12}$$

$$I_{Me}(\nu) = \left\langle \nu, \log \frac{\mathrm{d}\nu}{\kappa \nu \, \mathrm{d}\mu} \right\rangle + \frac{1}{2} \langle \nu, \kappa \mu \rangle. \tag{1.13}$$

As in the definition of $\mathbb{H}$ in (1.8), we define $I_{Ma}(\alpha) = \infty$ if it is not true that $\alpha$-almost everywhere $\frac{\mathrm{d}y}{(1 - e^{-\kappa y}) \mathrm{d}\mu}$ exists. Likewise we define $I_{Me}(\nu) = \infty$ if $\frac{\mathrm{d}\nu}{\kappa \nu \, \mathrm{d}\mu}$ does not exist. We use the convention that $\log 0 = -\infty$ and $0 \log 0 = 0$.

Let us recall the notion of an LDP: Theorem 1.1 says that $I(\cdot)$ is lower semi-continuous and, for any open set $G \subset \mathcal{L} \times \mathcal{A}$ and any closed set $F \subset \mathcal{L} \times \mathcal{A}$,

$$\liminf_{N \to \infty} \frac{1}{N} \log \mathbb{P}_N((\mathrm{Mi}_N, \mathrm{Ma}_N) \in G) \geq - \inf_G I(\cdot), \quad (1.14)$$

$$\limsup_{N \to \infty} \frac{1}{N} \log \mathbb{P}_N((\mathrm{Mi}_N, \mathrm{Ma}_N) \in F) \leq - \inf_F I(\cdot), \quad (1.15)$$

where we wrote $\mathbb{P}_N$ for the probability measure under the random graph $\mathcal{G}_N$. For the theory of large deviations, see e.g. [19].

An intuitive explanation of Theorem 1.1 is given in Sect. 1.3. The proof of Theorem 1.1 is in Sect. 5. It relies heavily on the special case of Theorem 1.1 for finite sets $\mathcal{S}$, see Theorem 3.1, whose proof we present first in Sect. 3. Our main strategy there is to identify the joint distribution of all the clusters in $\mathcal{G}(N, \mathbf{x}^{(N)}, \frac{1}{N}\kappa_N)$ in a combinatorial way and then to explicitly extract the exponential rates. The proof of Theorem 1.1 in Sect. 5 carries out an approximation procedure of $\mathcal{S}$ with finite state spaces in the spirit of the Dawson–Gärtner theorem.

Theorem 1.1 is an extension of [2, Theorem 1.1] from the special case $\mu = \delta_0$ and constant $\kappa$ (that is, from the standard Erdős–Rényi graph) to an inhomogeneous Erdős–Rényi graph. Note that this LDP is also highly non-trivial, interesting, and new in the case of an arbitrary $\mu$ and constant $\kappa$, to the best of our knowledge.

**Remark 1.2** (*Quenched and annealed LDPs*) One possible application of Theorem 1.1 is to the situation where the vertex types $x_1, \ldots, x_N$ are themselves random and independent with distribution $\mu$ each. Then Theorem 1.1 can be seen as a conditional LDP given $\mathbf{x}$, sometimes called a *quenched* LDP. The rate function turns out to be not random and depending only on $\mu$. One can then obtain an *annealed* version of the LDP, i.e., when the probabilities are also taken with respect to the vertices $x_1, \ldots, x_N$. The annealing follows from a standard mixture argument when $\mathcal{S}$ is a finite set of points; for general $\mathcal{S}$ the construction of a discretization suitable for use in our proof is a delicate matter that we do not explore here. One possible formulation of the annealed result would be that the triple, consisting of the empirical measures of the vertices, and $\mathrm{Mi}_N$ and $\mathrm{Ma}_N$ satisfies an LDP with rate function equal to $(\nu, \lambda, \alpha) \mapsto \mathbb{H}(\nu|\mu) + I_\nu(\lambda, \alpha)$, where we now wrote $I_\nu$ for the rate function $I$ of Theorem 1.1 with $\nu$ the limiting empirical measure of the type vector (instead of $\mu$).

**Remark 1.3** (*Detailedness*) We decided to register any component of the graph $\mathcal{G}_N$ only through its type-configuration, neglecting all the information about the internal connection structure. It is a natural wish to have also a more detailed analysis, for example by distinguishing each component as a sub*graph* instead of the type-configuration. From such a refined LDP, one could derive Theorem 1.1 via the contraction principle.

For the microscopic components it is indeed not too difficult to derive a refined version of the LDP of Theorem 1.1, since for each type-configuration $k \in \mathcal{M}_{\mathbb{N}_0}(\mathcal{S})$, the statistics of the $\approx N\lambda(k)$ components with type-configuration $k$ follow an explicit multinomial distribution. The form of the term $\tau(k)$ gives the hint that only spanning trees survive. The macroscopic components are much more involved and it is not clear which type of structure gives the decisive contribution in the limit.

Here is a standard corollary from the LDP in Theorem 1.1 about separate LDPs for $Mi_N$ and $Ma_N$.

**Corollary 1.4** (Separate LDPs for *$Mi_N$ and $Ma_N$*) *Under the assumptions of Theorem 1.1, the empirical measures $(Mi_N)_{N\in\mathbb{N}}$ and $(Ma_N)_{N\in\mathbb{N}}$ each satisfy an LDP on $\mathcal{L}$, respectively on $\mathcal{A}$, with rate functions*

$$\mathcal{I}_{Mi}(\lambda) = \inf_{\alpha\in\mathcal{A}} I(\lambda,\alpha) \quad \text{and} \quad \mathcal{I}_{Ma}(\alpha) = \inf_{\lambda\in\mathcal{L}} I(\lambda,\alpha).$$

The LDP assertion directly follows from the contraction principle (see [19]), since both projections $(\lambda,\alpha)\mapsto\lambda$ and $(\lambda,\alpha)\mapsto\alpha$ are continuous. The identification of the two contracted rate functions is formulated in Theorem 2.3 and discussed in Sect. 2.1.

***Remark 1.5*** (*LDP for the mesoscopic part*) Analogously to the corresponding result in [2], we could formulate and prove also a corollary about the *mesoscopic* part of the configuration $(\mathcal{C}_j)_j$ of the components of $\mathcal{G}_N$, i.e., about those components whose cardinalities satisfy $R < |\mathcal{C}_j| < \varepsilon N$ in the limit $N\to\infty$, followed by $R\to\infty$ and $\varepsilon\downarrow 0$. It is clear that we cannot consider the empirical measure of all these components anymore, but only the empirical measure of the total number of vertices of a given type in any of the mesoscopic components. Our conjecture is that (similarly to [2, Corollary 1.4]), this measure on $\mathcal{S}$ satisfies an LDP as $N\to\infty$ for fixed $R\in\mathbb{N}$ and $\varepsilon > 0$ with a rate that converges towards $I_{Me}$ defined in (1.13) as $R\to\infty$ and $\varepsilon\downarrow 0$. We abstained from writing out the details.

### 1.3 Interpretation of the LDP

Our main result, the LDP of Theorem 1.1, is highly compressed and contains a number of interesting results as special cases, so let us comment on the impact and draw some conclusions from it. We will restrict here to the large-deviations issues; the limiting issues and the consequences for the giant-cluster phase transition are deferred to Sect. 2.

We are examining the probability of the event $\{Mi_N \approx \lambda, Ma_N \approx \alpha\}$, asymptotically for large $N$, for any $\lambda\in\mathcal{L}$ and $\alpha = \sum_n \delta_{y_n} \in \mathcal{A}$. Indeed, we want to heuristically argue that, as $N\to\infty$, one has

$$\mathbb{P}_N(Mi_N \approx \lambda, Ma_N \approx \alpha) = e^{-NI(\lambda,\alpha)+o(N)}. \quad (1.16)$$

Recall from Sect. 1.1 that this is the event that $\mathcal{G}_N$ has $\sim N\lambda(k)$ components with type-configuration $k$, for any $k\in\mathcal{M}_{\mathbb{N}_0}(\mathcal{S})$, and a macroscopic component with type-configuration $\sim Ny_n$, for any $n$. We can clearly restrict to the case that $c_\lambda + c_\alpha \leq \mu$, since otherwise the number of vertices of some type in all the microscopic or macroscopic components together would be larger than the number of existing vertices of that type. However, it might be that the difference $\nu = \mu - c_\lambda - c_\alpha$ is a positive measure; this means that there are $\sim N\nu(R)$ of the vertices with type in $R$ in mesoscopic components for any $R \subset \mathcal{S}$, e.g., in components with $N$-dependent cardinalities like $\log N$ or $N^{1/3}$, or any mixture.

Recall that the types of all the vertices of $\mathcal{G}_N$ are approximately distributed as $\frac{1}{N}\sum_{i=1}^{N} \delta_{x_i} \approx \mu$. The key point in our proof is that the probability of $\{\mathrm{Mi}_N \approx \lambda, \mathrm{Ma}_N \approx \alpha\}$ consists of a number of terms that are more or less independent, i.e., lead to a sum of exponential rates. These terms are the following:

- a combinatorial term that expresses the number of decompositions of $[N]$ into the collection of subsets as above (respecting all the types),
- the probability that each of these subsets are connected (this depends on the type-configurations $k$ of the microscopic components and on the type-configuration $N y_n$ of the $n$-th macroscopic component, respectively),
- the probability that any two of all these vertex sets are not connected.

The above decomposition is the starting point of our combinatorial analysis in Lemma 3.3. Now, the crucial point in our LDP proof consists in the fact that we can get precise asymptotics for most of the terms in the decomposition of $\mathbb{P}_N(\mathrm{Mi}_N \approx \lambda, \mathrm{Ma}_N \approx \alpha)$. In particular, we rearrange the terms to isolate the contribution given by the microscopic components, respectively macroscopic, and collect in a remaining term all the rest, identifying this as the contribution of the mesoscopic part. In this rearrangement, we see that the probability that each subset is connected is the hardest term to handle. However, it turns out that it is enough to prove sharp asymptotics only in the case of finite size and order $N$ size components, for the remaining terms upper bounds are sufficient. Let us mention, in particular, that asymptotics of the probability of connectedness are not trivial at all for macroscopic components and obtaining such asymptotics is not only a crucial step in our LDP proof but also an interesting result on its own. We comment more on this in Sect. 1.4. Given this rearrangement, the large-deviations rate terms that we finally obtained are organized and interpreted in a slightly different fashion as follows.

Let us first consider the microscopic part. The first two terms in $I_{\mathrm{Mi}}(\lambda)$, $\mathbb{H}(\lambda|\mathbb{Q}_\mu)-1$, describe the number of labellings of $N|c_\lambda|$ vertices into microscopic subsets, according to $\lambda$ and respecting the type configurations. A priori, it is only notationally convenient to write this as an entropy, but this interpretation allowed us to make the step from discrete to continuous setting. The next three terms, $-\langle \lambda, \log \tau \rangle + |c_\lambda| - |\lambda|$, describe the probability that all the considered subsets are connected, see the comment below (1.9). The last term, $\frac{1}{2}\langle c_\lambda, \kappa\mu\rangle$, collects the costs of isolating each microscopic component from the rest of the system.

For the macroscopic part, given a macroscopic measure $\alpha = \sum_n \delta_{y_n}$, the term $\sum_n \langle y_n, \log \frac{dy_n}{d\mu}\rangle$ comes from the number of labellings of $N|c_\alpha|$ vertices into macroscopic subsets, according to $\alpha$ and to the type configurations. The term $-\sum_n \langle y_n, \log(1 - e^{-\kappa y_n})\rangle$ summarizes the connection probabilities of all the macroscopic components, see Theorem 3.6. Finally, for each $n$, the term $\frac{1}{2}\langle y_n, \kappa(\mu - y_n)\rangle$ is the cost of isolating the macroscopic component from the rest of the system.

In the mesoscopic part of the rate function, we see only the dependence on the measure $\nu$ of all the vertices of mesoscopic components, without any information about the components themselves. Again, the term $\langle \nu, \log \frac{d\nu}{d\mu}\rangle$ is a result of the relabelling of the $N|\nu|$ vertices according to the type configuration, and the integral of $-\log(\kappa\nu)$ with respect to $\nu$ describes in a summarizing way that each of the vertex sets is connected. The term $\frac{1}{2}\langle\nu,\kappa\mu\rangle$ represents the cost of isolating such vertices from the rest of the graph.

### 1.4 Connectivity of inhomogeneous graphs

On our way to a proof of Theorem 1.1, we derive some interesting formulas for quantities that are of general interest in the theory of multi-type Erdős–Rényi graphs. Indeed, in Lemma 3.3 we give a closed formula for the joint distribution of the entire collection of the vertex sets of all the components of $\mathcal{G}_N$. One important ingredient there is the probability for a given subset of vertices $\subset [N]$ to be connected. If the size of the subset is kept fixed, then it is straight-forward to get sharp estimates for the limit of this probability, as $N \to \infty$ (see Lemma 4.8). However, if the size of the subset is of order $N$, then it is much harder to derive sharp asymptotics for this probability. Our LDP relies precisely on this kind of asymptotics, which we derive in Theorem 3.6 in the case of a finite type space. We want to stress that finding Theorem 3.6 was crucial for proving the LDP and that we could not find any similar result in the literature that holds for inhomogeneous graphs. It requires a rather extensive proof, which is given in Sect. 4. Having established the LDP in its full generality one gets the same result about the connection probability in the more general setting of a compact type space.

**Corollary 1.6** (Connectivity probability of $\mathcal{G}_N$) *In the situation of Theorem 1.1 we have that*

$$\lim_{N\to\infty} \frac{1}{N} \log \mathbb{P}_N\big(\mathcal{G}(N, \boldsymbol{x}^{(N)}, \tfrac{1}{N}\kappa_N) \text{ is connected }\big) = \big\langle\mu, \log\big(1 - e^{-\kappa\mu}\big)\big\rangle. \quad (1.17)$$

The question about the connection probability of a random graph has attracted quite some interest. Let us mention [20], which studies for the inhomogeneous random graph the regime where the edge probability is of order $\frac{\log N}{N}$ and proves a phase transition: the probability of the graph to be connected either converges to 1 or to 0, depending on the parameters $\mu$ and $\kappa$ (in our notation). In our case, where the edge probability is of order $\frac{1}{N}$, we are in the case in which the probability of the graph being connected is always going to 0. Corollary 1.6 above identifies the exponential rate of its decay, which we think is of independent interest.

### 1.5 Related literature

This paper is a natural generalization to the inhomogeneous setting of [2], where we derived an LDP for all the cluster sizes of the Erdős–Rényi random graph in the sparse setting. Indeed, the classical sparse Erdős–Rényi graph corresponds to the case where $\mathcal{S}$ has only one element and $\kappa$ takes only one value, meaning that the result in [2] is a special case of the results in this paper. As mentioned in [2], the literature on the Erdős–Rényi graph is rich, but very few results on large deviations in the sparse regime are present. To the best of our knowledge, our paper is the first proving a large-deviations statement in the framework of the inhomogeneous graphs. Inhomogeneous random graphs have been introduced in [34] and the first mathematical treatment of the model has been presented in the seminal paper [7]. In [7] events that happen with high probability are studied, while the focus of the present paper is on rare events. In this section we list some earlier results concerning large deviations for random graphs and comment on their relation to our work.

In [2, Section 1.4], we gave a broad survey on known results on LDPs for sparse random graphs; summarizing, there are indeed some results on particular statistics of the graph $\mathcal{G}_N$, many of which are *a posteriori* contained in [2, Theorem 1.1] as special cases. For example in [32] two LDPs for the size of the largest component and for the number of isolated vertices have been derived. These quantities are continuous functionals of our measures $\mathrm{Ma}_N$, respectively of $\mathrm{Mi}_N$, and the contraction principle gives these results as consequences of ours.

It is important to mention the paper [8], where the authors proved an LDP for the empirical measure of the components rooted at each vertex in the sparse Erdős–Rényi graph. This object is a detailed size-biased version of our $\mathrm{Mi}_N$ and contains information about the edges that establish the connectedness of the components. However, in [8, Theorem 1.8], the authors show that their rate function is concentrated on trees, therefore any feasible microscopic component of size $k \in \mathbb{N}$ is indeed a spanning tree on $k$ vertices. This is also implicitly proven in our Lemma 3.4: the element $\tau$, defined in (1.9), ensures precisely that, when computing the probability that a microscopic component is formed on a certain finite set of vertices, the only important contribution is given by realizing a spanning tree on those vertices. This implies that a refinement of our proofs would give an LDP for a microscopic empirical measure on finite size components, as we mention in Remark 1.3. A size-biased version of such an empirical measure would correspond to a generalization to the inhomogeneous setting of the empirical neighborhood distribution in [8]. In this direction, in the very recent preprint [3], the authors deal with such an object in the case where the vertices of the graph have a type but the kernel $\kappa$ is constant. They prove an LDP using the techniques coming from [8] and relying on the notion of entropy for stochastic processes on marked rooted graphs introduced in [17]. Further investigation is needed (and desirable) to understand connections between this entropy, which comes out as the large-deviations rate function, and the rate function from our LDP.

In the framework of the sparse Erdős–Rényi graph, i.e., when the connection probability of $\mathcal{G}(N, p)$ satisfies $p \asymp N^{-1}$, recent progress has been made on the tails of triangle counts [12, 21], while we are not aware of the study of other rare events for the inhomogeneous graphs in such regime.

The case of large deviations in the dense Erdős–Rényi graph, i.e., for $\mathcal{G}(N, p)$ with fixed $p \in (0, 1)$, has been completely covered thanks to Chatterjee and Varadhan [14], see [13] for an overview. In [9, 28] extensions of the LDP to the framework of the dense inhomogeneous graphs are given. This regime is rather different from the sparse one and the LDP relies on the fact that each graph can be associated with a two dimensional symmetric function, called graphon. The limit of any sequence of sparse graphs is the graphon which is identically zero, showing that the space of graphons is not the right space in the sparse setting. Recently extensions to the sparse setting of the concept of graphons have been introduced, see [10], but connections to the sparse Erdős–Rényi graph to our knowledge are not known yet. Let us mention that the literature on LDPs for various statistics of different types of random graphs outside the sparse regime has seen a growth since the seminal work [14], see for example [15, 18]. However these graphs are by their nature so different from the setting of the present paper that we do not go into further details here.

For inhomogeneous graphs in the sparse regime, there are only a few results in the literature, starting with the seminal paper [7], which introduced the model and investigated the giant-cluster phase transition in detail. Furthermore, clusters of critical sizes of order $N^\alpha$ with some $\alpha \in (0, 1)$ around the phase transition have been studied for some types of inhomogeneous random graphs in [5], [6] and [36] under certain moment assumptions on the (scalar) types. Let us finally mention some results just outside the sparse setting: [11] analyzes the eigenvalues of the adjacency matrix of an inhomogeneous Erdős–Rényi random graph with vanishing edge probabilities of order $\gg \frac{1}{N}$. In [20] the authors study the probability of the graph to be connected when the edge probabilities are of order $\frac{\log N}{N}$. This compares to our result in Corollary 1.6, where we obtain precise asymptotics for such a probability in the sparse regime.

Finally, we have recently become aware of the preprint [27], where the authors study the fixed point equation (2.4) using combinatorial identities that we also rely on in Sect. 4. The focus of the paper is indeed to prove that a multi-type version of the Marcus–Lushnikov coagulation model with multiplicative kernel converges to the solution of the multi-type Flory equation (2.14). The particle masses of such a coagulation system are in one-to-one correspondence with the connected components of the inhomogeneous random graph, linking their convergence result to our Lemma 2.7.

## 2 Limiting consequences

In this section we present and discuss the second part of our main results, some consequences of the LDP of Theorem 1.1 that imply detailed and comprehensive limiting assertions about the inhomogeneous random graph. These results rely on involved variational analysis, using recursive formulas and elements of combinatorial power series analysis as methods to explicitly construct minimizers.

Like many large-deviations principles, also Theorem 1.1 implies a law of large numbers. This is particularly interesting here, since it implies and illustrates the well-known phase transition about the emergence of a giant cluster that was established in [7], and which we record and discuss in Sect. 2.1. There we also reveal our identification and interpretation of that phase transition in terms of the minimizer(s) of the rate function of our LDP in Theorem 1.1. In Sect. 2.2, we compare to the description of this phase transition that was given in [7] in terms of a strongly related multi-type branching process. Furthermore, in Sect. 2.3 we comment on the irreducibility of $\kappa$ and explain what the LDP looks like if this assumption is dropped. One of our main motivations for the present work is explained in Sect. 2.4 where we map the inhomogeneous Erdős–Rényi graph on a particular particle process with a random coagulating mechanism and discuss the consequences of our results for this process, in particular the phase transition of gelation type, i.e., the emergence of a gel, a macroscopic particle.

### 2.1 The phase transition

We give now comprehensive information about the well-known phase transition of the emergence of a giant component for the inhomogeneous Erdős–Rényi graph using the LDP of Theorem 1.1. Indeed, we derive a detailed picture of all the limiting microscopic and macroscopic clusters, according to their type-configurations. In this way, we go substantially beyond the work [7], to which we compare in Sect. 2.2 below. Unlike [7], our point of departure is not a multitype branching process, but the variational analysis of the rate function. This leads us in a natural way to deal with a transformed Poisson point process, which indeed shows deep connections with the multitype branching process.

We introduce first the minimizing microcluster distribution. Let us fix a kernel $\kappa$ as in Theorem 1.1. Recall that, for any measure $\nu$ on $\mathcal{S}$, we denote by $\mathbb{Q}_\nu$ the distribution of a Poisson point process on $\mathcal{S}$ with intensity measure $\nu$. For $c \in \mathcal{M}(\mathcal{S})$, we introduce the measure $\lambda_c$ on $\mathcal{M}_{\mathbb{N}_0}(\mathcal{S})$ by

$$\lambda_c(\mathrm{d}k) = \mathrm{e}^{\theta_c(\mathcal{S})} \tau(k) \mathbb{Q}_{\theta_c}(\mathrm{d}k), \qquad \text{where } \theta_c(\mathrm{d}r) = \mathrm{e}^{-\kappa c(r)} c(\mathrm{d}r). \tag{2.1}$$

In words, $\lambda_c$ is obtained by transforming the Poisson point process with intensity measure $\mathrm{e}^{-\kappa c(r)} c(\mathrm{d}r)$ with the function $\tau$. It will turn out in the subcritical case (and is implicit in the following theorem) that this measure possesses the integrated type configuration $c$, that is, the choice of the measure $\theta_c$ implies the crucial property that $c_{\lambda_c} = c$, where we recall the notation $c_\lambda(\mathrm{d}r) = \int \lambda(\mathrm{d}k)\, k(\mathrm{d}r)$.

We now introduce an important quantity (which was shown to be crucial in [7]) that we use for separating the sub- and supercritical regimes. For any measure $\nu$ on $\mathcal{S}$, we denote by $L^2(\nu)$ the usual $L^2$-space of functions $\mathcal{S} \to \mathbb{R}$ with respect to the measure $\nu$. We introduce the operator

$$T_{\kappa,\nu} \colon L^2(\nu) \to L^2(\nu), \qquad T_{\kappa,\nu} f(x) = \int_{\mathcal{S}} \kappa(x,y) f(y)\, \nu(\mathrm{d}y), \tag{2.2}$$

and its operator norm

$$\Sigma(\kappa, \nu) = \|T_{\kappa,\nu}\|_{L^2(\nu)} = \sup_{f \in L^2(\nu)\colon \|f\|_{L^2(\nu)}=1} \|T_{\kappa,\nu} f\|_{L^2(\nu)}. \tag{2.3}$$

Informally (an argument will follow in Sect. 4.1), in the special case that the support of $\nu$ is finite (i.e., the case of a finite set $\mathcal{S}$), then $T_{\kappa,\nu}$ can be identified with the matrix $(\kappa(r,s)\nu(s))_{r,s \in \mathcal{S}}$, and $\Sigma(\kappa, \nu)$ is its spectral radius.

Recall the rate function $I$ and the reference measure $\mu$ from Theorem 1.1, now we describe its minimizer.

**Theorem 2.1** (Minimizers of the rate function) *Suppose that $\kappa$ and $\mu$ are as in Theorem 1.1, then the following hold.*

(i) *If $\Sigma(\kappa, \mu) \leq 1$, then the unique minimizer of $I$ is equal to $(\lambda_\mu, 0)$.*
(ii) *If $\Sigma(\kappa, \mu) > 1$, then the unique minimizer of $I$ is equal to $(\lambda_{c^*}, \delta_{\mu-c^*})$, where the subprobability measure $c^*$ is the only solution to the characteristic equation*

$$e^{-\kappa c(r)} c(dr) = e^{-\kappa \mu(r)} \mu(dr) \quad \text{on } \mathcal{S}, \tag{2.4}$$

*that satisfies both $c^* \leq \mu$ and $c^* \neq \mu$. It further holds that $\Sigma(\kappa, c^*) < 1$.*

*In particular,*

$$\left(\text{Mi}_N, \text{Ma}_N\right) \stackrel{N \to \infty}{\Longrightarrow} \begin{cases} (\lambda_\mu, 0) & \text{if } \Sigma(\kappa, \mu) \leq 1, \\ (\lambda_{c^*}, \delta_{\mu-c^*}) & \text{if } \Sigma(\kappa, \mu) > 1. \end{cases} \tag{2.5}$$

The proof of Theorem 2.1 is in Sect. 7. Most of it is original research, but take from [7] the discussion of the solutions of the fixed point equation (2.4), see Lemma 4.1 where we summarize it.

The law of large numbers in (2.5) is a standard consequence of an LDP with a unique minimizer for the rate function. This is a very precise and detailed formulation of the famous giant-cluster phase transition in the graph $\mathcal{G}_N$. Indeed, the following happens with probability tending to one exponentially fast:

(i) In the subcritical phase $\Sigma(\kappa, \mu) < 1$, all vertices (meaning all up to $o(N)$) are in microscopic components, more precisely in the unique optimal configuration encoded by $\lambda_\mu$. That is, for any $k \in \mathcal{M}_{\mathbb{N}_0}(\mathcal{S})$, the number of components with vertex set given by $k$ is asymptotically $N e^{\theta_\mu(\mathcal{S})} \tau(k) \mathbb{Q}_{\theta_\mu}(dk)$, with $\theta_\mu(dr) = e^{-\kappa \mu(r)} \mu(dr)$. We have $c_{\lambda_\mu} = \mu$, no giant component appears, and the number of vertices in mesoscopic components is $o(N)$.
(ii) In the case $\Sigma(\kappa, \mu) > 1$, a unique giant cluster appears with $\sim N(1 - c^*(\mathcal{S}))$ vertices and type-configuration asymptotically equal to $N(\mu - c^*)$ with $c^*$ characterized by (2.4), since $\theta_\mu = \theta_{c^*}$. The microscopic components are distributed according to the optimal distribution $\lambda_{c^*}$, and the number of vertices in mesoscopic components is $o(N)$. Note that this microscopic distribution is not saturated, that is $\Sigma(\kappa, c^*) < 1$, as in the one-type setting [2]. That is, we encounter a phase transition of explosion type, rather than of saturation type, see Remark 2.2 and [2, Section 1.6].

**Remark 2.2** (*Phase transition: saturation versus explosion*) Here is an explanation of the phase transition in terms of a dynamical process. Consider a process of Erdős–Rényi graphs in increasing connection probability, i.e., by adding more and more bonds between the vertices, such that components grow. A suitable growth parameter is $\Sigma(\kappa, c)$, where $c$ stands for the rescaled type-configuration of all the vertices; however, we consider $\Sigma(\kappa, c)$ as a growing function of $\kappa$, the bond density. As we explained in [2, Section 1.6], the well-known giant-cluster phase transition (see also the discussion below Theorem 1.1 below) is an explosion phase transition in the sense that, when crossing the threshold one, a positive fraction of finite-size clusters merges rapidly into one giant cluster and, at any time, every cluster keeps participating in merge events. In particular, the total microscopic mass starts decreasing at the phase transition. In contrast, in condensation phase transitions like the famous Bose–Einstein condensation first all microscopic components reach their maximal size (the saturated state), before a macroscopic component, the condensate, appears, and then the microsocopic ones do not change anymore, but all of the additional mass goes exclusively into the condensate.

Now we give the description of the two rate functions for the contracted LDPs of $(\text{Mi}_N)_{N \in \mathbb{N}}$ and $(\text{Ma}_N)_{N \in \mathbb{N}}$, respectively, from Corollary 1.4.

**Theorem 2.3** (Minimizers of the contracted rate function) *Suppose that $\kappa$ and $\mu$ are as in Theorem 1.1. Then the following hold for $\lambda \in \mathcal{L}$ and for $\alpha \in \mathcal{A}$, respectively.*

$$\mathcal{I}_{Mi}(\lambda) = \inf_{\alpha \in \mathcal{A}} I(\lambda, \alpha) = \begin{cases} I_{Mi}(\lambda) + I_{Ma}(\delta_{\mu - c_\lambda}) & \text{if } c_\lambda \leq \mu, \\ \infty & \text{otherwise,} \end{cases} \quad (2.6)$$

*and*

$$\mathcal{I}_{Ma}(\alpha) = \inf_{\lambda \in \mathcal{L}} I(\lambda, \alpha) = \begin{cases} I_{Ma}(\alpha) + J(\mu - c_\alpha) & \text{if } c_\alpha \leq \mu, \\ \infty & \text{otherwise,} \end{cases} \quad (2.7)$$

*where for $c \in \mathcal{M}(\mathcal{S})$*

$$J(c) = \begin{cases} I_{Mi}(\lambda_c) = \langle c, \log \frac{dc}{d\mu} \rangle + \frac{1}{2} \langle c, \kappa(\mu - c) \rangle & \text{if } \Sigma(\kappa, c) \leq 1, \\ I_{Mi}(\lambda_{b^*}) + I_{Me}(c - b^*) & \text{if } \Sigma(\kappa, c) > 1, \end{cases} \quad (2.8)$$

*where $b^* = b^*(c) \in \mathcal{M}(\mathcal{S})$ is the minimal non-trivial (i.e., not equal to $c$) solution to*

$$\kappa(c - b^*)(r) \, b^*(dr) = (c - b^*)(dr), \quad b^* \leq c, \quad (2.9)$$

*and it holds that $\Sigma(\kappa, b^*) = 1$.*

The proof of Theorem 2.3 is in Sect. 6 for $\mathcal{I}_{Mi}$ and in Sect. 7 for $\mathcal{I}_{Ma}$. The above theorem suggests us a conditional law of large numbers. Informally, if we fix $\alpha \in \mathcal{A}$ and a sequence of $\alpha_N \in \mathcal{A}$ such that $\alpha_N \to \alpha$. Then, under the probability $\mathbb{P}_N(\cdot | \text{Ma}_N = \alpha_N)$ we have

$$\text{Mi}_N \overset{N \to \infty}{\Longrightarrow} \begin{cases} \lambda_{\mu - c_\alpha} & \text{if } \Sigma(\kappa, \mu - c_\alpha) \leq 1, \\ \lambda_{b^*} & \text{if } \Sigma(\kappa, \mu - c_\alpha) > 1. \end{cases} \quad (2.10)$$

Notice that when $\Sigma(\kappa, \mu - c_\alpha) > 1$, then $b^*$ is not equal to $\mu - c_\alpha$, therefore $b^* + c_\alpha \neq \mu$ and the missing mass is interpreted as being mesoscopic.

**Remark 2.4** (*Conditional limit with saturation phase transition*) In formula (2.8) we encounter a phase transition of saturation type in a conditional limit, in contrast to a transition of explosion type of the unconditional one, see Remark 2.2. We refer back to Sect. 1.3 for the interpretation. Recall that $\mathcal{I}_{\text{Ma}}(\alpha)$ is the negative exponential rate of the probability of the event $\{\text{Ma}_N \approx \alpha\}$. When $\Sigma(\kappa, \mu - c_\alpha) \leq 1$, (2.8) shows that $\mathcal{I}_{\text{Ma}}(\alpha) = I_{\text{Ma}}(\alpha) + I_{\text{Mi}}(\lambda_{\mu-c_\alpha})$, implying that $\mathbb{P}_N(\text{Ma}_N \approx \alpha) = \mathbb{P}_N(\text{Mi}_N \approx \lambda_{\mu-c_\alpha}, \text{Ma}_N \approx \alpha)e^{o(N)}$. The interpretation of this is that, conditionally on the event $\{\text{Ma}_N \approx \alpha\}$, the non-macroscopic mass optimally organizes according to the microscopic measure $\lambda_{\mu-c_\alpha}$. In contrast, in the case $\Sigma(\kappa, \mu - c_\alpha) > 1$, from (2.8) we see that $\mathcal{I}_{\text{Ma}}(\alpha) = I_{\text{Ma}}(\alpha) + I_{\text{Mi}}(\lambda_{b^*}) + I_{\text{Me}}(\mu - c_\alpha - b^*)$. This means that, conditionally on the event $\{\text{Ma}_N \approx \alpha\}$, the non-macroscopic mass cannot be organized fully in microscopic clusters, but it is organized microscopically according to $\lambda_{b^*}$ and the remaining vertices, with type-configuration $N(\mu - c_\alpha - b^*)$, are put in mesoscopic components. The particular rescaled type-configuration $b^*$ is saturated in the sense that $\Sigma(\kappa, b^*) = 1$. This means, given a fixed macroscopic type-configuration, if more bonds are thrown into the graph, then first all microscopic clusters grow until they reach the saturated state $\lambda_{b^*}$, and then this is frozen and only mesoscopic clusters grow. The latter effect is present in literature under the name of frozen percolation, see for example [16, 29, 33, 37] and the difference between the two phase transitions is reflected in substantial differences of the hydrodynamic limit, as we summarize in Sect. 2.4.

### 2.2 Comparison to [7]: branching-process interpretation

Our description of the limiting quantities that we presented in Sect. 2.1 is based on and derived from our analysis of the minimizer of $I$. Therefore we found it most suitable to present them in terms of transformed Poisson point processes. However, in the analysis of finite-size components of random graphs, it is common and was often successful to employ well-adapted branching processes for the description. The main idea is that a component can be efficiently (sampled and) analyzed by exploring it via such a branching algorithm. This idea was also a cornerstone in the seminal paper [7], and it produced a description of the limiting macroscopic component in terms of the extinction probability of a crucial branching process. In this section, we recall this description and compare it to our Poisson point process description, also including the microscopic components.

The main tool that is utilized in [7] is a *multitype branching process* with type space $\mathcal{S}$, in which each particle of type $r \in \mathcal{S}$ has offspring with distribution that is a Poisson process with intensity measure $\kappa(r, s)\mu(ds)$. We define $\rho(r) \in [0, 1]$ as the probability of non-extinction of the branching process, if it starts with precisely one particle that has type $r \in \mathcal{S}$. We summarize the most important facts from [7] that have relevance for our comparison as follows (see [7, Th. 3.1, Th. 3.12, Th. 6.1, Th. 9.10]).

**Theorem 2.5** (Existence of a giant component, [7]) *Suppose the situation of Theorem 1.1 is given. Abbreviate $\mathcal{G}_N = \mathcal{G}(N, \mathbf{x}, \frac{1}{N}\kappa_N)$, then the following hold.*

(i) $\rho \colon \mathcal{S} \to [0, \infty)$ *is the maximal solution of*

$$\rho = 1 - e^{-T_{\kappa,\mu}\rho}. \tag{2.11}$$

(ii) *If* $\Sigma(\kappa, \mu) \leq 1$, *then the largest component of* $\mathcal{G}_N$ *has size* $O(\log N)$ *as* $N \to \infty$ *with high probability.*

(iii) *If* $\Sigma(\kappa, \mu) > 1$, *then the largest component* $\mathcal{C}_1$ *of* $\mathcal{G}_N$ *has size* $\asymp N$. *More precisely, its normalized empirical measure* $\frac{1}{N}\eta_{\boldsymbol{x}}(\mathcal{C}_1)$ *(recall* (1.1)) *converges weakly towards the measure* $\rho(r)\,\mu(\mathrm{d}r)$, *and* $\rho$ *is positive* $\mu$-*almost everywhere in* $\mathcal{S}$.

Part (iii) identifies the limiting type-configuration of the giant component as $N$ times the measure with density $\rho$ with respect to $\mu$, and part (i) characterizes $\rho$ via the functional identity (2.11). It is easily seen to be equivalent to the characteristic equation (2.4) that we use via the substitution $\rho(r)\,\mu(\mathrm{d}r) = \mu(\mathrm{d}r) - c^*(\mathrm{d}r)$ or $c^*(\mathrm{d}r) = (1 - \rho(r))\,\mu(\mathrm{d}r)$. In our analysis of the minimizer of $I$, (2.4) arose via the Euler–Lagrange equations, while (2.11) emerged in [7, Lemma 5.4] via a standard formula for mixed moments of the offspring of the branching process (which itself uses standard Poisson point process theory).

The statement in part (ii) about the order of the largest component is out of reach of our large-deviations ansatz, which implies that all but $o(N)$ vertices are in components of finite size.

About the distribution of the microscopic clusters of $\mathcal{G}_N$, however, there is no explicit result contained in [7]. However, we can give a description in terms of the above branching process as well. We derive this description now from our form of the minimizer $\lambda_\mu$ defined in (2.1). Let $\Xi(\mathrm{d}r)$ be the total progeny of type $r$ of the branching process; then $\Xi$ is a random measure on $\mathcal{S}$. Then, if $\mathrm{P}_r$ denotes the measure if the process starts from one individual of type $r$ at time 0, we have

$$\mu(\mathrm{d}r)\mathrm{P}_r(\Xi \in \mathrm{d}k) = \lambda_\mu(\mathrm{d}k)\,k(\mathrm{d}r), \qquad k \in \mathcal{M}_{\mathbb{N}_0}(\mathcal{S}) \setminus \{0\}, r \in \mathcal{S}. \tag{2.12}$$

In Remark 4.6 we explain this relation in the setting where $\mathcal{S}$ is a finite set. In words, the empirical statistics of the microscopic components in $\mathcal{G}_N$ in the subcritical case approximate the distribution of the total offspring of the characteristic branching process.

### 2.3 The reducible case

Let us briefly comment on the case where the kernel $\kappa$ is reducible with respect to $\mu$, i.e., $\mathcal{S}$ is composed of at least two irreducible classes (maximal irreducible subsets). Then the graph $\mathcal{G}_N$ decomposes into disconnected subgraphs with types in only one of these classes. Accordingly, the collection of all the connected components can be decomposed into collections for each subgraph. In principle one can apply the LDP of Theorem 1.1 to each of the micro/macro empirical measures of the subgraphs. However, one might have the wish to have a joint LDP for the entire collection. Here one might expect that the same LDP holds true, and the decomposition into the subgraphs reappears in the rate function in a natural way.

It turns out that this expectation is not disappointed, as it concerns the microscopic part, but is disappointed for the macroscopic part. Actually, the formulation of the LDP slightly changes. The main point is that two macroscopic components can very cheaply be connected to form a significantly larger macroscopic component, just by throwing in one connecting edge, which cannot be seen on the exponential scale. Hence, the macroscopic part is very unstable on the exponential scale under adding edges. This argument fails for the microscopic part, since here we are talking now about $\asymp N$ independent copies of a component of a finite size; if one wants to connect them such that the microscopic statistics change, then one needs to change $\asymp N$ edges, whose probability is clearly seen on an exponential scale.

As a consequence of this effect, we will see that the rate function is finite only if the macroscopic measure $\alpha = \sum_n \delta_{y_n}$ is such that for each $n$ the rescaled type-configuration $y_n$ is supported in one of the irreducible classes. To be precise, we say that a measure $y \in \mathcal{M}(\mathcal{S})$ is *connectable* (with respect to $\kappa$ and $\mu$) if its support is contained in an irreducible class. Furthermore, a measure $\alpha \in \mathcal{A}$ is called *connectable* if each of its atoms is connectable.

For the microscopic configurations $\lambda$ connectability is implicitly ensured by the fact that the measure $\lambda$ has to be absolutely continuous with respect to $\tau \mathbb{Q}_\mu$ in order to have a finite value of the rate function; notice that $\tau(k) = 0$ if $\text{supp}(k)$ is not concentrated on a irreducible class of $\mathcal{S}$. We also have to restrain to a sequence of random graphs that are defined with respect to the same kernel $\kappa$, rather than an approximating sequence $\kappa_N$, since each $\kappa_N$ might be irreducible.

**Theorem 2.6** (LDP in the reducible case) *Suppose the setting as in Theorem* 1.1, *with the exception that $\kappa_N = \kappa$ for all $N \in \mathbb{N}$ and the kernel $\kappa$ on $\mathcal{S} \times \mathcal{S}$ is now assumed to be reducible, and suppose that $\mathcal{S}$ is equal to the support of $\mu$. Then $(\text{Mi}_N, \text{Ma}_N)$ satisfies an LDP with rate function $\widetilde{I}$ defined by*

$$\widetilde{I}(\lambda, \alpha) = \begin{cases} I(\lambda, \alpha), & \text{if } \alpha \text{ is connectable,} \\ +\infty & \text{otherwise,} \end{cases}$$

*where $I$ is as in Theorem* 1.1.

The proof follows in a straightforward way from our results, see Remark 3.11 for the finite type case. The intuitive reason is that, for $\alpha$ that is not connectable, on the event $\{\text{Ma}_N \approx \alpha\}$, there is a macroscopic component who has two non-trivial parts (i.e., each with $\asymp N$ types) in two different irreducible classes, even though there cannot be any edge between these sets. Hence this event has the probability zero.

We made the choice to state the theorem under the additional assumption that $\mathcal{S}$ corresponds to the support of $\mu$. If this was not the case, one could still approach a not connectable measure $y \in \mathcal{M}(\mathcal{S})$ with a finite exponential cost if, for each finite $N$, the support of $y$ is contained in an irreducible class of $\kappa$ with respect to $\mu^{(N)}$. The proof would anyway be an extension of our finite type case results, but it is out of our scope to cover this particular framework.

### 2.4 Motivation: an inhomogeneous coagulation process

With the present work, we actually continue our study of random particle models with coagulation in the light of large-deviation arguments initiated in [2]. Indeed, we make here the first step towards a *spatial* model.

The model that we are interested in is the following. Fix $N$ atoms $1, \ldots, N$ at the locations $x_1, \ldots, x_N$ in a compact metric space $\mathcal{S}$. We consider a Markov process in continuous time on the set of partitions of $[N] = \{1, \ldots, N\}$. Starting with the monodispersed configuration $M(0) = (\{i\})_{i \in [N]}$, at any time any two subsets $A$, $B$ in the current partition are replaced by their union $A \cup B$ after an exponentially distributed random time with parameter

$$\frac{1}{N} \sum_{i \in A, j \in B} \kappa(x_i, x_j), \tag{2.13}$$

where a symmetric $\kappa \colon \mathcal{S} \times \mathcal{S} \to [0, \infty)$ is given. All these random times are supposed to be independent. If $M(t)$ denotes the partition at time $t$, then $M(s)$ is a refinement of $M(t)$ for any $s < t$. Hence, the number of elements of $M(t)$ is a non-increasing (random) process starting at $N$. The special case of a singleton $\mathcal{S}$ and $\kappa \equiv 1$ (the homogeneous case) is the case of the *Marcus–Lushnikov model* that we studied in [2]. There we also explained how the Erdős–Rényi model can be mapped onto the Marcus–Lushnikov model, and this works also in the inhomogeneous setting. Indeed, to any unordered pair $\{i, j\} \subset [N]$ with $i \neq j$ we associate an exponential random time $e(i, j)$ with parameter $\frac{1}{N} \kappa(x_i, x_j)$. These random times are independent and we put a bond between $i$ and $j$ as soon as $e(i, j)$ elapses. At a fixed time $t \in (0, \infty)$, this graph has the distribution of the inhomogeneous random graph $\mathcal{G}_{t,N} = \mathcal{G}([N], (x_1, \ldots, x_N), \frac{1}{N} \kappa_{t,N})$ with type space $\mathcal{S}$ and

$$\kappa_{t,N}(r, s) = N\big(1 - e^{-\frac{1}{N} t \kappa(r,s)}\big), \qquad r, s \in \mathcal{S}.$$

Notice that the random partition $M(t)$ of the above coagulation model is equal in distribution to the collection $(\mathcal{C}_j)_j$ of the vertex sets of the components of $\mathcal{G}_{t,N}$. The two main reasons for this fact are the memorylessness of the exponential distribution and the property that the minimum of independent exponential times is also exponential with a parameter that is the sum of all the parameters. The only difference between the two models is that the graph model registers all the bonds that arrive within each of the components (and do not change anything in the connectedness), while the coagulation model just registers that a given set is connected.

We are interested in an LDP for the micro and the macro empirical measure of the partition sets of $M(t)$ in the limit $N \to \infty$, assuming the initial locations of particles are such that $\frac{1}{N} \sum_{i=1}^{N} \delta_{x_i} \to \mu$ for some measure $\mu$ on $\mathcal{S}$. Since $\kappa_{t,N} \to t\kappa$, Theorem 1.1 applies also to the above inhomogeneous coagulation process under the appropriate assumptions at a fixed time $t$. Furthermore, from Theorem 2.1, we obtain that the process $(M(t))_{t\in[0,\infty)}$ has a phase transition at the time

$$t_c = \frac{1}{\Sigma(\kappa,\mu)},$$

and we have a limiting distribution of the empirical micro and macro measures. This phase transition is of explosion type, as described in Remark 2.2, and in the coagulation literature is usually called *gelation*. Further consequences for the limiting distribution of $M(t)$ as $N \to \infty$ follow in a natural way, but we refrain from writing them down.

Interestingly, we can deduce that the minimizing process of microscopic cluster sizes satisfies the multitype version of the *Flory equation*, which is a modification of the well-known *Smoluchowski equation*. The classical (single-type) Smoluchowski and Flory equation are ubiquitous in the literature concerning coagulation processes, see for example [1]. A multi-type extension to the Flory equation can be formulated as follows. We think of an inhomogeneous deterministic coagulation process $(\lambda_t)_{t\in[0,\infty)}$, conceived as a process in $\mathcal{L}$. Each particle $k \in \mathcal{M}_{\mathbb{N}_0}(\mathcal{S})$ consists of $k(\mathcal{S})$ atoms, $k(\{r\})$ of which have the type $r$ for any $r \in \mathcal{S}$. Coagulation is nothing but addition of measures in this formulation, i.e., two particles $k$ and $\widetilde{k}$ coagulate to a particle $k + \widetilde{k}$. The kernel of this process is given as

$$K(k,\widetilde{k}) = \langle \widetilde{k}, \kappa k \rangle = \int_{\mathcal{S}^2} \kappa(r,s)\, k(\mathrm{d}r)\widetilde{k}(\mathrm{d}s), \qquad k, \widetilde{k} \in \mathcal{M}_{\mathbb{N}_0}(\mathcal{S}).$$

Then the weak formulation of the Flory equation is, for any test function $f \in \mathcal{C}_c(\mathcal{M}_{\mathbb{N}_0}(\mathcal{S}))$,

$$\begin{aligned}\frac{\mathrm{d}}{\mathrm{d}t}\int_{\mathcal{M}_{\mathbb{N}_0}(\mathcal{S})} f(k)\lambda_t(\mathrm{d}k) &= \frac{1}{2}\int_{\mathcal{M}_{\mathbb{N}_0}(\mathcal{S})^2} f(k+\widetilde{k})K(k,\widetilde{k})\lambda_t(\mathrm{d}k)\lambda_t(\mathrm{d}\widetilde{k}) \\ &\quad - \int_{\mathcal{M}_{\mathbb{N}_0}(\mathcal{S})^2} f(k)K(k,\widetilde{k})\lambda_t(\mathrm{d}k)\lambda_t(\mathrm{d}\widetilde{k}) \\ &\quad - \int_{\mathcal{M}_{\mathbb{N}_0}(\mathcal{S})^2} f(k)K(k,\widetilde{k})\lambda_t(\mathrm{d}k)(\lambda_0-\lambda_t)(\mathrm{d}\widetilde{k}), \quad (2.14)\end{aligned}$$

where $\lambda_0(\mathrm{d}k) = \int_{\mathcal{S}} \mu(\mathrm{d}r)\, \delta_{\delta_r}(\mathrm{d}k)$ is the initial condition, which expresses that $\mu$ is the atom type distribution. In words, the time-evolution of $(\lambda_t)_{t\in[0,\infty)}$ is described by saying that any coagulation of two particles $k$ and $\widetilde{k}$ (i.e., replacement of $k$ and $\widetilde{k}$ by $k + \widetilde{k}$) happens with rate $K(k,\widetilde{k})$. In our model the last term in the right-hand side can be rewritten as $-\int_{\mathcal{M}_{\mathbb{N}_0}(\mathcal{S})} f(k)\langle k, \kappa(\mu - c_{\lambda_t})\rangle \lambda_t(\mathrm{d}k)$. It captures the interaction between the microscopic particles and the gel (the macroscopic mass) once it forms.

See [30, Section 3] for a mathematical discussion of the (well-known) homogeneous version of the Flory equation and [31, Section 2] for an introduction of the inhomogeneous version of the equation. We now identify a solution $(\lambda_t)_{t\in[0,\infty)}$ to (2.14).

**Lemma 2.7** (Solution to the Flory equation) *Assume that $\mathcal{S}$ is a finite state space and $\kappa$ an irreducible nonnegative symmetric matrix on $\mathcal{S}$. Let $\lambda_0(k) = \sum_{r \in \mathcal{S}} \mu_r \delta_{\delta_r}(k)$ and for $t \in (0, \infty)$, define $\lambda_t$ to be the first component of the minimizer appearing in Theorem* 2.1 *with $\kappa$ replaced by $t\kappa$.*

*Then $t \mapsto \lambda_t$ is a solution to the Flory equation* (2.14) *on $[0, \infty)$.*

The proof of Lemma 2.7, as well as an explicit expression for $(\lambda_t)_{t \geq 0}$, is given in Sect. 7.4. We are confident that Lemma 2.7 is also true in the general setting of Theorem 2.1.

The Flory equation is closely related to the *Smoluchowski equation*, which we write in its multitype version:

$$\frac{d}{dt} \int_{\mathcal{M}_{\mathbb{N}_0}(\mathcal{S})} f(k) \lambda_t(dk) = \frac{1}{2} \int_{\mathcal{M}_{\mathbb{N}_0}(\mathcal{S})^2} f(k + \widetilde{k}) K(k, \widetilde{k}) \lambda_t(dk) \lambda_t(d\widetilde{k})$$
$$- \int_{\mathcal{M}_{\mathbb{N}_0}(\mathcal{S})^2} f(k) K(k, \widetilde{k}) \lambda_t(dk) \lambda_t(d\widetilde{k}). \quad (2.15)$$

The Smoluchowski equation only considers the microscopic clusters, that is, it excludes any interaction with a possible gel, which we see in the third line of the Flory equation (2.14). The solutions of equations (2.14) and (2.15) coincide until the gelation time $t_c = 1/\Sigma(\kappa, \mu)$, after which differences appear.

At the level of the underlying stochastic microscopic models this difference is seen in terms of the type of the phase transitions, as cited in Remark 2.2. The microscopic models of *frozen percolation*-type, as in [16, 29, 33, 37], correspond to (2.15), while models like ours correspond to (2.14).

## 3 Proof of Theorem 1.1 for a finite type set

In this section we assume that $\mathcal{S}$ is a finite set and derive the large-deviations principle (LDP) of Theorem 1.1 for this case, Theorem 3.1. This is not only an important special case that is worth being formulated and studied on its own, but it will be the first step in the proof of Theorem 1.1 that is completed in Sect. 5. The formulation in the discrete case is notationally pretty different from the formulation in the general setting and many objects simplify because of the finiteness of $\mathcal{S}$. Therefore we are going to formulate the setting and the LDP from scratch in Sect. 3.1.

The organization of this section is as follows. In Sect. 3.2 we derive a formula for the distribution of $\mathrm{Mi}_N$. The more involved terms that appear in our formula are certain connection probabilities whose asymptotics are stated in Sect. 3.3. In Sect. 3.4 we will decompose the distribution of $\mathrm{Mi}_N$ into a micro–, meso– and macroscopic part and derive the exponential rates for each part. The proof of the LDP of Theorem 3.1 is finally finished in Sect. 3.5.

### 3.1 Formulation of the LDP

Let us recall the objects that we need to formulate the LDP for a finte type space $\mathcal{S}$. Fix a probability measure $\mu$ on $\mathcal{S}$, which we will denote as a vector $\mu = (\mu_s)_{s\in\mathcal{S}}$. For any $N \in \mathbb{N}$ let $\mathbf{x}^{(N)} = (x_1^{(N)}, \ldots, x_N^{(N)}) \in \mathcal{S}^N$ be a type vector such that the normalized empirical measure, $\mu^{(N)} = \frac{1}{N}\sum_{i=1}^N \delta_{x_i^{(N)}}$, converges to $\mu$ as $N \to \infty$. Let $\kappa = (\kappa(r,s))_{r,s\in\mathcal{S}} \in [0,\infty)^{\mathcal{S}\times\mathcal{S}}$ be a nonnegative and symmetric matrix. For any $N \in \mathbb{N}$, let $\kappa_N$ be another such matrix, and assume that the sequence $\kappa_N$, $N \in \mathbb{N}$, converges pointwise to $\kappa$ as $N \to \infty$. Throughout the section we will work under the assumptions that we just stated for $\kappa_N$, $N \in \mathbb{N}$, and $\mu^{(N)}$, $N \in \mathbb{N}$.

Recall that the random graph $\mathcal{G}_N = \mathcal{G}(N, \mathbf{x}^{(N)}, \frac{1}{N}\kappa_N)$ consists of $N$ vertices and that the vertex $i \in [N]$ has $x_i^{(N)} \in \mathcal{S}$. The undirected edges of the graph are set independently for each pair of vertices and two vertices of type $r$ and $s$ are connected via an edge with probability $1 \wedge \frac{1}{N}\kappa_N(r,s)$.

We denote by $\{\mathcal{C}_j\}_j$ the collection of the vertex sets of all the connected components of $\mathcal{G}_N = \mathcal{G}(N, \mathbf{x}^{(N)}, \frac{1}{N}\kappa_N)$. We want to study empirical measures depending on the random collection $\{\mathcal{C}_j\}_j$. For this we introduce the type-registering mapping $\eta\colon \mathcal{P}([N]) \to \mathbb{N}_0^{\mathcal{S}}$ that gives the (type) composition of an arbitrary vertex set $A \subset [N]$, i.e., $\eta(A) = (\eta_s(A))_{s\in\mathcal{S}}$, and $\eta_s(A) = \#\{i \in A\colon x_i^{(N)} = s\}$ is the number of vertices in $A$ with type $s$. Note that the mapping $\eta$ depends on the entire type vector $\mathbf{x}^{(N)}$, not only on its normalized empirical measure, $\mu^{(N)} = \frac{1}{N}\sum_{i=1}^N \delta_{x_i^{(N)}}$.

We identify $\mathcal{M}_{\mathbb{N}_0}(\mathcal{S})$ with $\mathbb{N}_0^{\mathcal{S}}$ and will work in $[0,1]^{\mathcal{S}}$ instead of $\mathcal{M}(\mathcal{S})$. Now we recall the definition of the main objects, the empirical measures of the connected components of $\mathcal{G}_N$, in microscopic, respectively macroscopic, registration. The *microscopic empirical measure* $\mathrm{Mi}_N$ is defined as a measure on $\mathbb{N}_0^{\mathcal{S}}$ via

$$\mathrm{Mi}_N(\mathrm{d}k) := \frac{1}{N}\sum_j \delta_{\eta(\mathcal{C}_j)}(\mathrm{d}k). \tag{3.1}$$

Since $\mathbb{N}_0^{\mathcal{S}}$ is a discrete space, we will abbreviate $\mathrm{Mi}_N(k) = \mathrm{Mi}_N(\{k\})$ for any $k \in \mathbb{N}_0^{\mathcal{S}}$. The *macroscopic empirical measure* is defined as a measure on $[0,1]^{\mathcal{S}}\setminus\{0\}$ via

$$\mathrm{Ma}_N(\mathrm{d}y) := \sum_j \delta_{\frac{1}{N}\eta(\mathcal{C}_j)}(\mathrm{d}y). \tag{3.2}$$

Therefore our state spaces for $\mathrm{Mi}_N$ and $\mathrm{Ma}_N$ are now

$$\mathcal{L} := \{\lambda = (\lambda_k)_{k\in\mathbb{N}_0^{\mathcal{S}}} \in [0,\infty)^{\mathbb{N}_0^{\mathcal{S}}} : c(\lambda) \leq \mu \text{ or } c(\lambda) \leq \mu^{(N)} \text{ for some } N \in \mathbb{N}, \lambda_0 = 0\} \tag{3.3}$$

and

$$\mathcal{A} := \{\alpha \in \mathcal{M}_{\mathbb{N}_0}([0,1]^{\mathcal{S}}\setminus\{0\}) : c(\alpha) \leq \mu \text{ or } c(\alpha) \leq \mu^{(N)} \text{ for some } N \in \mathbb{N}\}, \tag{3.4}$$

respectively, where

$$c_r(\lambda) = \sum_{k \in \mathbb{N}_0^{\mathcal{S}}} \lambda_k k_r \quad \text{and} \quad c_r(\alpha) = \int_{[0,1]^{\mathcal{S}} \setminus \{0\}} \alpha(\mathrm{d}y) \, y_r, \quad r \in \mathcal{S}. \tag{3.5}$$

One can easily verify that for any fixed $N$ we have $c_r(\text{Mi}_N) = \mu_r^{(N)}$ as well as $c_r(\text{Ma}_N) = \mu_r^{(N)}$ for any $r \in \mathcal{S}$, so indeed $\text{Mi}_N \in \mathcal{L}$ and $\text{Ma}_N \in \mathcal{A}$. However, the idea is that due to the topologies that we choose, some of the (rescaled) vertex mass specified by $c(\lim \text{Mi}_N)$ may get lost when we take the limit for $N \to \infty$.

We equip $\mathcal{L}$ and $\mathcal{A}$ with the vague topologies that we introduced in Sect. 1.2. On $\mathcal{L}$, this is identical with the topology of pointwise convergence (i.e., $\lim_{N \to \infty} \lambda^{(N)} = \lambda$ if and only if $\lim_{N \to \infty} \lambda_k^{(N)} = \lambda_k$ for any $k \in \mathbb{N}_0^{\mathcal{S}}$). The vague topology on $\mathcal{A}$ is formulated by saying that $\lim_{N \to \infty} \alpha^{(N)} = \alpha$ if and only if $\lim_{N \to \infty} \int \alpha^{(N)}(\mathrm{d}y) f(y) = \int \alpha(\mathrm{d}y) f(y)$ for any continuous and compactly supported function $f : [0, 1]^{\mathcal{S}} \setminus \{0\} \to \mathbb{R}$; note that for every such function $f$ there is an $\varepsilon > 0$ such that $f = 0$ on $\{x \in [0, 1]^{\mathcal{S}} : |x| \le \varepsilon\}$.

Recall that we write $\langle a, f \rangle = \sum_r a_r f_r$ for the integral of a function $f$ with respect to a measure $a$ on $\mathcal{S}$ and also recall the notation $|a| = \sum_{s \in \mathcal{S}} a_s$. Further, recall the combinatorial quantity $\tau(k)$ that collects the weight of all spanning trees on a vertex set with type configuration $k \in \mathbb{N}_0^{\mathcal{S}}$, i.e.,

$$\tau(k) = \sum_{T \in \mathcal{T}(k)} \prod_{\{i,j\} \in E(T)} \kappa(x_i, x_j), \tag{3.6}$$

where $\mathbf{x} \in \mathcal{S}^{|k|}$ is a type vector compatible with $k$, i.e., $\sum_{i=1}^{|k|} \delta_{x_i} = k$, and $\mathcal{T}(k)$ is the set of spanning trees on $[|k|]$. We use the convention that $\mathcal{T}(0) = \emptyset$ and hence $\tau(0) = 0$.

Here is the main result of Sect. 3:

**Theorem 3.1** (LDP for $(\text{Mi}_N, \text{Ma}_N)$ with finitely many types) *Assume that the empirical measure $\mu^{(N)}$ of the type sequence $(x_1^{(N)}, \ldots, x_N^{(N)})$ converges weakly towards a positive probability vector $\mu \in (0, 1]^{\mathcal{S}}$ as $N \to \infty$ and that the kernel $\kappa_N$ converges on $\mathcal{S} \times \mathcal{S}$ towards a $\mu$-irreducible kernel $\kappa \in [0, \infty)^{\mathcal{S} \times \mathcal{S}}$.*

*Then $(\text{Mi}_N, \text{Ma}_N)$ satisfies a large deviations principle (LDP) with speed $N$ and rate function $(\lambda, \alpha) \mapsto I(\lambda, \alpha)$ defined by*

$$I(\lambda, \alpha) = \begin{cases} I_{Mi}(\lambda) + I_{Ma}(\alpha) + I_{Me}(\mu - c(\lambda) - c(\alpha)), & \text{if } c(\lambda) + c(\alpha) \le \mu, \\ +\infty & \text{otherwise,} \end{cases}$$

*where*

$$I_{Mi}(\lambda) = \sum_{k \in \mathbb{N}_0^{\mathcal{S}}} \lambda_k \log \frac{\lambda_k}{\tau(k) \prod_{s \in \mathcal{S}} \frac{\mu_s^{k_s}}{k_s!}} + \sum_{k \in \mathbb{N}_0^{\mathcal{S}}} \lambda_k(|k| - 1) + \frac{1}{2} \langle c(\lambda), \kappa \mu \rangle, \tag{3.7}$$

$$I_{Ma}(\alpha) = \int_{\mathcal{M}(\mathcal{S})\setminus\{0\}} \alpha(\mathrm{d}y) \left(\left\langle y, \log \frac{y}{(1-\mathrm{e}^{-\kappa y})\mu}\right\rangle + \frac{1}{2}\langle y, \kappa(\mu-y)\rangle\right), \quad (3.8)$$

$$I_{Me}(\nu) = \left\langle \nu, \log \frac{\nu}{(\kappa\nu)\mu}\right\rangle + \frac{1}{2}\langle \nu, \kappa\mu\rangle \tag{3.9}$$

*and where we always use the convention that* $\log 0 = -\infty$ *and* $0 \log 0 = 0$.

Theorem 3.1 is indeed nothing but the special case of Theorem 1.1 for a finite set $\mathcal{S}$. Indeed, it is clear that the setting and the two rate functions $I_{Ma}$ and $I_{Me}$ are the discrete-space versions, but (3.7) looks a bit different from the formula for $I_{Mi}$ in Theorem 1.1. But from substituting the notation of the entropy in (3.7) and noting that the distribution of a Poisson point process with intensity measure $\mu$ can here be identified as

$$\mathbb{Q}_\mu(k) = \prod_{r\in\mathcal{S}} \left[\mathrm{e}^{-\mu_r} \frac{\mu_r^{k_r}}{k_r!}\right], \quad k \in \mathbb{N}_0^\mathcal{S},$$

one sees that (3.7) is indeed a discrete analog of (1.11). Furthermore, one can also write

$$I_{Mi}(\lambda) = \sum_{k\in\mathbb{N}_0^\mathcal{S}} \lambda_k \log \frac{\lambda_k}{\lambda_k(\mu)} + \sum_{k\in\mathbb{N}_0^\mathcal{S}} \lambda_k(|k|-1) - \frac{1}{2}\langle c(\lambda), \kappa\mu\rangle \tag{3.10}$$

where $\lambda(\mu)$ is defined as the discrete analog of (2.1) with $c = \mu$, i.e.,

$$\lambda_k(\mu) = \tau(k) \prod_{s\in\mathcal{S}} \frac{(\mu_s \mathrm{e}^{-(\kappa\mu)_s})^{k_s}}{k_s!}, \quad k \in \mathbb{N}_0^\mathcal{S}. \tag{3.11}$$

This formula will be helpful in Sect. 6 when we will identify minimizers of $I_{Mi}$.

We will now give an extension of Theorem 3.1 for kernels $\kappa$ that are not $\mu$-irreducible. Recall the notion of connectability for $\alpha \in \mathcal{A}$ that was introduced in Sect. 2.3.

**Theorem 3.2** (Finite-type LDP for $(Mi_N, Ma_N)$ without irreducibility) *For* $(\lambda, \alpha) \in \mathcal{L} \times \mathcal{A}$ *define*

$$\widetilde{I}(\lambda, \alpha) = \begin{cases} I(\lambda, \alpha), & \text{if } \alpha \text{ is connectable,} \\ +\infty & \text{otherwise.} \end{cases} \tag{3.12}$$

*(i) Given all the assumptions from Theorem 3.1, except the assumption that $\kappa$ is $\mu$-irreducible, the pair $(Mi_N, Ma_N)$ satisfies the lower large-deviations bound (1.14) with speed $N$ and rate function $\widetilde{I}$.*

*(ii) Given all the assumptions from Theorem 3.1, except the assumption that $\kappa$ is $\mu$-irreducible and with the additional assumption that $\kappa_N = \kappa$ for all but finitely many $N \in \mathbb{N}$, the pair $(Mi_N, Ma_N)$ satisfies an LDP with speed $N$ and rate function $\widetilde{I}$.*

The proof is given in Remark 3.11.

We omit to restate the finite-$\mathcal{S}$ analogs of Theorems 2.3, 2.1 and all the related corollaries, as they can be deduced as special cases. For the critical quantity $\Sigma(\kappa, \mu)$, we refer to (4.2), and we recall that it is in this setting equal to the spectral radius of the matrix $(\kappa(r, s)\mu_s)_{(r,s)\in\mathcal{S}^2}$.

### 3.2 The distribution of $\text{Mi}_N$

Let us identify the distribution of $\text{Mi}_N$ for any $N$ in explicit terms. Note that as long as $N \in \mathbb{N}$ is fixed, the measure $\text{Ma}_N$ contains exactly the same information as $\text{Mi}_N$, hence, we are also deriving its distribution. We start by noting that $N\text{Mi}_N$ takes values in

$$\mathcal{L}_N := \left\{ \ell = (\ell_k)_{k \in \mathbb{N}_0^\mathcal{S}} : \ell_k \in \mathbb{N}_0 \text{ for all } k, \ \ell_0 = 0 \text{ and} \right.$$
$$\left. \sum_{k \in \mathbb{N}_0^\mathcal{S}} \ell_k k_r = N\mu_r^{(N)} \text{ for all } r \in \mathcal{S} \right\}. \tag{3.13}$$

Let $k \in \mathbb{N}_0^\mathcal{S}$ and let $\mathbf{x} = (x_1, \ldots, x_{|k|}) \in \mathcal{S}^{|k|}$ be a type vector which is compatible with $k$, meaning that $\sum_{i=1}^{|k|} \delta_{x_i} = k$. We define the connection probability of the graph $\mathcal{G}(|k|, \mathbf{x}, \frac{1}{N}\kappa_N)$ by

$$p_N(k) = \mathbb{P}\big(\mathcal{G}(|k|, \mathbf{x}, \tfrac{1}{N}\kappa_N) \text{ is connected}\big), \quad \text{for } k \neq 0, \tag{3.14}$$

and $p_N(0) = 0$. In the following lemma we write down the distribution of $\text{Mi}_N$ in terms of the quantities $p_N(k)$, $k \in \mathbb{N}_0^\mathcal{S}$.

**Lemma 3.3** (The distribution of $Mi_N$) *Let $N \in \mathbb{N}$ and assume that $\kappa_N(r, s) \leq N$ for all $r, s \in \mathcal{S}$. Then for any $\ell \in \mathcal{L}_N$ we have that*

$$\mathbb{P}\left(N Mi_N(k) = \ell_k \ \forall k \in \mathbb{N}_0^\mathcal{S}\right) = \left(\prod_{r \in \mathcal{S}} \left(N\mu_r^{(N)}\right)!\right) \times \prod_{k \in \mathbb{N}^\mathcal{S}} \zeta^{(N)}(\ell, k), \tag{3.15}$$

*where*

$$\zeta^{(N)}(\ell, k) = \frac{p_N(k)^{\ell_k}}{\ell_k! \prod_{r \in \mathcal{S}} (k_r!)^{\ell_k}} \left( \prod_{r,s \in \mathcal{S}} \left(1 - \frac{\kappa_N(r, s)}{N}\right)^{\frac{1}{2}k_r[N\mu_s^{(N)} - k_s]} \right)^{\ell_k}, \tag{3.16}$$

*and $p_N(k)$ is defined in (3.14).*

**Proof** This is proved in an analogous way to [2, Corollary 2.2]; we omit the details. □

The formula in (3.15) is easy to understand. Indeed, the combinatorial term on the right (with an additional factor $1/N!$) is equal to the inverse of the number of possible labelings of all the vertices; the event $\{N\mathrm{Mi}_N(k) = \ell_k\ \forall k\}$ means that $\ell_k$ is equal to the number of clusters in $[N]$ whose vertex set has the type configuration $k$, for any multi-index $k$. The product of $p_N(k)^{\ell_k}$ over $k$ is the probability that all these clusters are connected, the product over the powers of $1 - \kappa_N(r,s)/N$ is equal to the probability that each two of them are not connected, and the product of the two remaining combinatorial terms (with an additional factor $N!$) is equal to the number of ways to decompose all the types into clusters having the prescribed vertex structure.

### 3.3 Asymptotics for the connection probabilities

Recall the connection probability $p_N(k)$, for $k \in \mathbb{N}_0^{\mathcal{S}}$, that we defined in (3.14). When going to the limit for $N \to \infty$, it will be crucial to distinguish different cases depending on the asymptotic behaviour of $k$. In the first case, we keep $k \in \mathbb{N}_0^{\mathcal{S}}$ fixed, whereas $N \to \infty$ and refer to $p_N(k)$ as the connection probability of a microscopic cluster. In the second case we consider a sequence $k^{(N)} \in \mathbb{N}_0^{\mathcal{S}}$ where $|k^{(N)}|$ is of order $N$. More precisely, we assume that for any $s \in \mathcal{S}$ the limit $\lim_{N\to\infty} \frac{k_s^{(N)}}{N} = y_s$ exists and that the vector $y = (y_s)_{s\in\mathcal{S}}$ is non-trivial. In that case we refer to $p_N(k^{(N)})$ as the connection probability of a macroscopic cluster. In between the microscopic and the macroscopic regime there are the cases, in which the sequence $|k^{(N)}|$ diverges, but is in $o(N)$. Those are summarized under the notion of mesoscopic clusters.

In this section the results for the different cases will only be stated. Their proofs are collected in Sect. 4, since the derivation of the asymptotics for the macroscopic case is rather cumbersome. Recall the definition (3.6) for $\tau$. By $\tau_N$ we will denote the same quantity, but defined with respect to the kernel $\kappa_N$. Further, recall that we are working under the assumption that $\kappa_N$ converges pointwise to $\kappa$, which will be used in Lemma 3.4 and Theorem 3.6.

**Lemma 3.4** (Asymptotics for the connection probability of microscopic clusters) *Fix $k \in \mathbb{N}_0^{\mathcal{S}}$. Then, as $N \to \infty$,*

$$p_N(k) = N^{-(|k|-1)} \tau_N(k)(1 + o(1)). \tag{3.17}$$

**Lemma 3.5** (Estimate for the connection probability of mesoscopic clusters) *Fix $k \in \mathbb{N}_0^{\mathcal{S}}$ and choose any $r \in \mathcal{S}$ such that $k_r > 0$. Then*

$$p_N(k) \le (\|\kappa_N\|_\infty |S_k|)^{|S_k|-1} N^{-(|k|-1)} \frac{1}{k_r^2} \Big( \prod_{s \in S_k} (\kappa_N k)_s^{k_s - 1} k_s \Big), \tag{3.18}$$

*where $S_k := supp(k)$.*

The following theorem concerns the connection probabilities of macroscopic clusters. This result is to the best of our knowledge, a new one and might be of independent interest in the theory of random graphs. Note that it is the multi-type version of a result from [35], see [2, Lemma 2.4].

**Theorem 3.6** (Asymptotics of the connection probability of macroscopic clusters) *Fix $y \in [0,1]^\mathcal{S}$, $y \neq 0$. Let $\{k^{(N)}\}_{N \in \mathbb{N}}$ be a sequence in $\mathbb{N}_0^\mathcal{S}$ such that $\lim_{N \to \infty} \frac{k_r^{(N)}}{N} = y_r$ for all $r \in \mathcal{S}$.*

*(i) Then it holds that*

$$\limsup_{N \to \infty} \frac{1}{N} \log p_N(k^{(N)}) \leq \sum_{r \in \mathcal{S}} y_r \log\left(1 - e^{-(\kappa y)_r}\right), \quad (3.19)$$

*where the right-hand side takes the value $-\infty$ when $y \not\ll \kappa y$.*

*(ii) Assume that $\tau(k^{(N)}) > 0$ for all but finitely many $N \in \mathbb{N}$ and that $\{k_r^{(N)}\}_N$ is bounded in $N$ for all $r \notin \mathrm{supp}(y)$. Then*

$$\lim_{N \to \infty} \frac{1}{N} \log p_N(k^{(N)}) = \sum_{r \in \mathcal{S}} y_r \log\left(1 - e^{-(\kappa y)_r}\right) \in [-\infty, 0], \quad (3.20)$$

*where the right-hand side takes the value $-\infty$ when $y \not\ll \kappa y$.*

The additional assumption $\tau(k^{(N)}) > 0$ ensures that the connection probability is indeed strictly positive since otherwise the left-hand side of (3.20) is $-\infty$. The assumption about the boundedness of $\{k_r^{(N)}\}_N$ for $r \notin \mathrm{supp}(y)$ might be weakened. However, for our purposes the statement will suffice in this form.

The proof can be found in Sect. 4. The main idea is to construct a sequence of graphs with the same connection parameter $\kappa_N$, but a different number of vertices in such a way that it contains with high probability a macroscopic component $k^{(N)}$.

### 3.4 Exponential rates for micro-, meso- and macro parts

The proof of the LDP in Theorem 3.1 is carried out in the same way as in [2, Section 3]. The main idea is to split the distribution that we obtained in Lemma 3.3 into three parts, which we will call micro-, meso- and macroscopic part. These parts roughly give the terms $e^{-NI_{\mathrm{Mi}}(\lambda)}$, $e^{-NI_{\mathrm{Me}}(\mu - c(\lambda) - c(\alpha))}$ and $e^{-NI_{\mathrm{Ma}}(\alpha)}$, if a properly rescaled version of $\ell = \ell^{(N)} \in \mathcal{L}_N$ is close to $(\lambda, \alpha)$. In the next lemma we give the decomposition into the three parts. Afterwards, we derive the exponential asymptotics of them in Lemmas 3.8–3.10.

**Lemma 3.7** (Decomposition into three contributions) *Fix $\ell \in \mathcal{L}_N$. For $k \in \mathbb{N}_0^\mathcal{S}$ define*

$$z^{(N)}(\ell_k, k) := \frac{p_N(k)^{\ell_k} \prod_{s \in \mathcal{S}}\left(\frac{N\mu_s}{e}\right)^{k_s \ell_k}}{\ell_k! \prod_{s \in \mathcal{S}}(k_s!)^{\ell_k}} \prod_{r,s \in \mathcal{S}}\left(1 - \frac{\kappa_N(r,s)}{N}\right)^{\frac{1}{2}k_s(N\mu_r - k_r)\ell_k}, \quad (3.21)$$

*and for any two numbers $A, B \in [0, \infty)$ write*

$$z_{A,B}^{(N)}(\ell) := \prod_{k \in \mathbb{N}_0^\mathcal{S} : A < |k| \leq B} z^{(N)}(\ell_k, k).$$

*Then, for any fixed $R \in \mathbb{N}$ and $\varepsilon > 0$ we have that, as $N \to \infty$*

$$\mathbb{P}_N(NMi_N = \ell) = e^{o(N)} z_{0,R}^{(N)}(\ell) z_{R,\varepsilon N}^{(N)}(\ell) z_{\varepsilon N,N}^{(N)}(\ell).$$

**Proof** On the right-hand side of (3.15) we apply Stirling's formula $n = e^{o(n)}(n/e)^n$ to the terms $(N\mu_r^{(N)})!$ and use that $N\mu_r^{(N)} = \sum_k \ell_k k_r$. Note that we also used that $\mu_r^{(N)} \to \mu_r$. □

**Lemma 3.8** (Asymptotics of the micro part) *Fix $R \in \mathbb{N}$ and let $\lambda \in \mathcal{L}$. Define*

$$I_{Mi}^{(R)}(\lambda) := \sum_{k \in \mathbb{N}_0^S \colon |k| \leq R} \lambda_k \log \frac{\lambda_k e^{|k|-1}}{\tau(k) \prod_s \frac{\mu_s^{k_s}}{k_s!}} + \frac{1}{2} \langle c^{(R)}(\lambda), \kappa\mu \rangle, \qquad (3.22)$$

*where*

$$c_r^{(R)}(\lambda) = \sum_{k \in \mathbb{N}_0^S \colon |k| \leq R} \lambda_k k_r, \quad r \in \mathcal{S}, \qquad (3.23)$$

*and where we use the convention that $\log 0 = -\infty$ and $0 \log 0 = 0$. In particular, the right-hand side of (3.22) is equal to $+\infty$ if there is some $k \in \mathbb{N}_0^S$ with $|k| \leq R$ such that $\tau(k) = 0$, but $\lambda_k > 0$. Otherwise, the right-hand side is finite.*

*Then for all $\ell^{(N)} = (\ell_k^{(N)})_{k \in \mathbb{N}_0^S} \in \mathcal{L}_N$ satisfying $\lambda_k = \lim_{N \to \infty} \frac{1}{N} \ell_k^{(N)}$ for all $k$ with $|k| \leq R$ we have*

$$\lim_{N \to \infty} \frac{1}{N} \log z_{0,R}^{(N)}(\ell^{(N)}) = -I_{Mi}^{(R)}(\lambda). \qquad (3.24)$$

**Proof** We use Stirling's formula for the terms $\ell_k^{(N)}!$ as well as the fact that for any $k \in \mathbb{N}_0^S$ with $|k| \leq R$ we have that $p_N(k) = e^{o(1)} N^{1-|k|} \tau_N(k)$ as $N \to \infty$, by Lemma 3.4. Therefore, as $N \to \infty$, we have

$$\prod_{k \colon |k| \leq R} \frac{p_N(k)^{\ell_k^{(N)}} \prod_{s \in \mathcal{S}} (\frac{N\mu_s}{e})^{k_s \ell_k^{(N)}}}{\ell_k^{(N)}! \prod_{s \in \mathcal{S}} (k_s!)^{\ell_k^{(N)}}}$$

$$= e^{o(N)} e^{-\frac{1}{2} \sum_{k \colon |k| \leq R} \log \ell_k^{(N)}} \prod_{|k| \leq R} \left( \frac{\tau_N(k) \prod_s \frac{\mu_s^{k_s}}{k_s!}}{(\ell_k^{(N)}/N) e^{|k|-1}} \right)^{\ell_k^{(N)}}$$

$$= e^{o(N)} \exp\left( N \sum_{|k| \leq R} \lambda_k \log \frac{\tau(k) \prod_s \frac{\mu_s^{k_s}}{k_s!}}{\lambda_k e^{|k|-1}} \right), \qquad (3.25)$$

where, we have used that $0 \leq \sum_{k \colon |k| \leq R} \log \ell_k^{(N)} \leq |\{k \colon |k| \leq R\}| \log\left(\frac{\sum_{k \colon |k| \leq R} \ell_k^{(N)}}{|\{k \colon |k| \leq R\}|}\right) = o(N)$. Further, we use that $\lim_{N \to \infty}(1 + \frac{x}{N})^N = e^x$ to

get that, as $N \to \infty$,

$$\prod_{|k| \leq R} \prod_{r,s \in \mathcal{S}} \left(1 - \frac{\kappa_N(r,s)}{N}\right)^{\frac{1}{2}k_s(N\mu_r - k_r)\ell_k^{(N)}} = \exp\left(-\frac{N}{2}\left\langle \sum_{|k| \leq R} \frac{\ell_k^{(N)}}{N} k, \kappa(\mu - \frac{k}{N})\right\rangle\right)$$
$$= e^{o(N)} e^{-\frac{N}{2}\langle c^{(R)}(\lambda), \kappa\mu\rangle}, \qquad (3.26)$$

where in the last step we used that

$$0 \leq \frac{1}{2} \sum_{|k| \leq R} \frac{\ell_k^{(N)}}{N} \sum_{r,s} \frac{k_r}{N} \kappa_N(r,s) k_s \leq \frac{1}{2} \|\kappa_N\|_\infty \frac{R}{N} |c^{(R)}(\frac{\ell^{(N)}}{N})|,$$

together with the fact that the right-hand side converges to 0 as $N \to \infty$. Combining the asymptotics from (3.25) and (3.26) gives the claim. $\square$

**Lemma 3.9** (Asymptotics of the macro part) *Fix $\alpha \in \mathcal{A}$, and note that $\alpha$ can be written as $\alpha = \sum_{j \in J} \delta_{y^{(j)}}$ where $y^{(j)} \in [0,1]^{\mathcal{S}} \setminus \{0\}$ for all $j \in J$ and $J$ is a countable set. Fix any $\varepsilon > 0$ with $\varepsilon \notin \{|y^{(j)}| : j \in J\}$. Define $J_\varepsilon(\alpha) := \{j \in J : |y^{(j)}| > \varepsilon\}$, which is a finite set, and*

$$I_{Ma}^{(\varepsilon)}(\alpha) := \int \alpha(\mathrm{d}y) \mathbb{1}_{\{|y|>\varepsilon\}} \left\langle y, \log \frac{y}{(1-e^{-\kappa y})\mu}\right\rangle + \frac{1}{2} \int \alpha(\mathrm{d}y) \mathbb{1}_{\{|y|>\varepsilon\}} \langle y, \kappa(\mu - y)\rangle, \qquad (3.27)$$

*where we use the convention that $\log 0 = -\infty$ and $0 \log 0 = 0$. In particular, the right-hand side of (3.27) is equal to $+\infty$, if there is some $i \in J_\varepsilon(\alpha)$ such that the condition $y^{(i)} \ll \kappa y^{(i)}$ fails. Then we have the following.*

(i) *For any sequence $\ell^{(N)} \in \mathcal{L}_N$ denote $\alpha^{(N)} = \sum_k \ell_k^{(N)} \delta_{\frac{k}{N}}$ and assume that $\alpha^{(N)}$ restricted to $\{y : |y| > \varepsilon\}$ converges to $\alpha$ restricted to $\{y : |y| > \varepsilon\}$, as $N \to \infty$. Then it holds that*

$$\limsup_{N \to \infty} \frac{1}{N} \log z_{\varepsilon N, N}^{(N)}(\ell^{(N)}) \leq -I_{Ma}^{(\varepsilon)}(\alpha). \qquad (3.28)$$

(ii) *For all $j \in J_\varepsilon(\alpha)$, let $\{k^{(j,N)}\}_{N \in \mathbb{N}}$ be a sequence in $\mathbb{N}_0^{\mathcal{S}}$ such that $\tau(k^{(j,N)}) > 0$ for all $N \in \mathbb{N}$, $\lim_{N \to \infty} \frac{k^{(j,N)}}{N} = y^{(j)}$ and $\{k_s^{(j,N)}\}_{N \in \mathbb{N}}$ is bounded for all $s \notin \mathrm{supp}(y^{(j)})$. Let $\ell^{(N)}$ be an element of $\mathcal{L}_N$ such that $\ell_k^{(N)} = \#\{j \in J : k^{(j,N)} = k\}$ for $|k| > \varepsilon N$. Denote $\alpha^{(N)} = \sum_k \ell_k^{(N)} \delta_{\frac{k}{N}}$, then $\alpha^{(N)}$ restricted to $\{y : |y| > \varepsilon\}$ converges to $\alpha$ restricted to $\{y : |y| > \varepsilon\}$, and*

$$\lim_{N \to \infty} \frac{1}{N} \log z_{\varepsilon N, N}^{(N)}(\ell^{(N)}) = -I_{Ma}^{(\varepsilon)}(\alpha). \qquad (3.29)$$

**Proof** We start with the first statement. Let us first turn to the first term on the right of (3.21). We apply Stirling's bound $n! \geq n^n e^{-n}$ to each of the terms $k_s!$ and the simple bound $\ell_k^{(N)}! \geq 1$ for all $k \in \mathbb{N}_0^{\mathcal{S}}$. Using the upper bound (3.19) for $p_N(k)$ from Theorem 3.6 we obtain, as $N \to \infty$,

$$\prod_{k:\,|k|>\varepsilon N} \frac{p_N(k)^{\ell_k^{(N)}} \prod_{s\in\mathcal{S}}(\frac{N\mu_s}{e})^{k_s\ell_k^{(N)}}}{\ell_k^{(N)}! \prod_{s\in\mathcal{S}}(k_s!)^{\ell_k^{(N)}}}$$
$$\leq e^{o(N)} \exp\left(N \int_{\{|y|>\varepsilon\}} \left\langle y, \log\frac{(1-e^{-\kappa y})\mu}{y}\right\rangle \alpha^{(N)}(\mathrm{d}y)\right), \quad (3.30)$$

where we used that $\sum_{|k|>\varepsilon N} \ell_k^{(N)} \leq \frac{1}{\varepsilon}$. The second term on the right of (3.21) is estimated using $1 - x \leq e^{-x}$ for $x = \kappa_N(r,s)/N$ as follows for $N \to \infty$,

$$\prod_{|k|>\varepsilon N}\prod_{r,s\in\mathcal{S}} \left(1 - \frac{\kappa_N(r,s)}{N}\right)^{\frac{1}{2}k_s(N\mu_r-k_r)\ell_k^{(N)}}$$
$$\leq e^{o(N)} \exp\left(-\frac{N}{2}\int_{\{|y|>\varepsilon\}} \langle y, \kappa(\mu-y)\rangle \alpha^{(N)}(\mathrm{d}y)\right). \quad (3.31)$$

Note that the product of the right-hand sides of (3.30) and (3.31) is equal to $e^{o(N)}e^{-NI_{\mathrm{Ma}}^{(\varepsilon)}(\alpha^{(N)})}$. Using the convergence assumption on $\alpha^{(N)}$ and the fact that $\varepsilon < \inf_{j\in J_\varepsilon} |y^{(j)}|$ one can verify that $I_{\mathrm{Ma}}^{(\varepsilon)}(\alpha^{(N)}) \to I_{\mathrm{Ma}}^{(\varepsilon)}(\alpha)$ as $N \to \infty$. This gives the result.

To show the second statement, let us first notice that it is clear from the definition that $\alpha^{(N)}$ restricted to $\{y: |y| > \varepsilon\}$ converges to $\alpha$ restricted to $\{y: |y| > \varepsilon\}$. Therefore we can apply the first assertion and we have the upper bound (3.28). In order to get also the lower bound, we lower estimate $z_{\varepsilon N,N}^{(N)}(\ell^{(N)})$ against the sum on $k^{(N,j)}$ over the finite set $j \in J_\varepsilon(\alpha)$ and note that these are the only summands $k$ with $|k| > \varepsilon N$ such that $\ell_k^{(N)} > 0$. For each such $j$ we apply the asymptotic (3.20) from Theorem 3.6 and obtain the corresponding lower bound, also noting that $\prod_{|k|>\varepsilon N}(\ell_k^{(N)}!)^{-1} \geq e^{-\frac{1}{\varepsilon}\log N} = e^{o(N)}$ and, by estimating $k_s \leq N$, we get that $(k_s)^{-\frac{\ell_k^{(N)}}{2}} \geq e^{-\frac{1}{2\varepsilon}\log N}$ for any $s \in \mathcal{S}$. To derive the lower bound equivalent of (3.31) we use that $(1 - \frac{c}{N})^N = e^{-c}(1+o(1))$. This proves the claim. □

In the second statement of Lemma 3.9 we restrict to sequences $\ell^{(N)}$ such that each type-configuration $k$ with $\ell_k^{(N)} > 0$ is connectable in the sense that $\tau(k) > 0$. We will see in the proof of Theorem 3.1 that for each $\alpha \in \mathcal{A}$ such that $I_{\mathrm{Ma}}(\alpha) < \infty$ and for each $\varepsilon > 0$, we will always be able to find such a sequence. Our second restriction $\varepsilon \notin \{|y^{(j)}|: j \in J\}$ is clearly only technical and gives no problem at all when we later take the limit as $\varepsilon \downarrow 0$; it frees us from unwanted terms.

**Lemma 3.10** (Asymptotics of the meso part) *Fix $R \in \mathbb{N}$, $\varepsilon > 0$ and $\nu \in [0,1]^{\mathcal{S}}$. For a sequence $\ell^{(N)} \in \mathcal{L}_N$, $N \in \mathbb{N}$, we write $\nu_s^{(N)} := \frac{1}{N}\sum_{R<|k|\leq \varepsilon N} \ell_k^{(N)} k_s$, for $s \in \mathcal{S}$,* *and assume that* $\lim_{N\to\infty} \nu_s^{(N)} = \nu_s$ *holds for all* $s \in \mathcal{S}$. *Then*

$$\limsup_{N\to\infty} \frac{1}{N} \log z_{R,\varepsilon N}^{(N)}(\ell^{(N)}) \leq -I_{Me}(\nu) + C(R, \varepsilon, \nu), \tag{3.32}$$

*where the term* $C(R, \varepsilon, \nu)$ *is continuous in* $\nu$ *and converges to* 0 *as* $R \to \infty$ *and* $\varepsilon \to 0$.

**Proof** For fixed $k \in \mathbb{N}_0^{\mathcal{S}}$ denote $S_k := \text{supp}(k)$. For $p_N(k)$ and $r \in S_k$ we use the upper bound (3.18) from Lemma 3.5. Also, we apply the Stirling bound $n! \geq n^n e^{-n}$ to the terms $\ell_k^{(N)}!$ as well as the terms $k_s!$ for all $s \in \mathcal{S}$. This gives that

$$\prod_{R<|k|\leq \varepsilon N} \frac{p_N(k)^{\ell_k^{(N)}} \prod_{s\in\mathcal{S}} (\frac{N\mu_s}{e})^{k_s \ell_k^{(N)}}}{\ell_k^{(N)}! \prod_{s\in\mathcal{S}} (k_s!)^{\ell_k^{(N)}}} \leq \prod_{R<|k|\leq \varepsilon N} \Big(\frac{NC \prod_{s\in S_k} (\kappa_N k)_s^{k_s-1} k_s}{k_r^2 \ell_k^{(N)} e^{-1} \prod_{s\in S_k} (\frac{k_s}{\mu_s})^{k_s}}\Big)^{\ell_k^{(N)}}$$

$$\leq e^{N\langle \nu^{(N)}, \log \mu \rangle} \times \prod_{R<|k|\leq \varepsilon N} \prod_{s\in S_k} \Big(\frac{(\kappa_N k)_s}{k_s}\Big)^{(k_s-1)\ell_k^{(N)}} \times \prod_{R<|k|\leq \varepsilon N} \Big(\frac{NeC}{\ell_k^{(N)}(k_r)^2}\Big)^{\ell_k^{(N)}}, \tag{3.33}$$

where $C = 1 \wedge ((\|\kappa\|_\infty + 1)|\mathcal{S}|)^{|\mathcal{S}|-1}$ and where we assumed that $N$ is large enough such that $\|\kappa_N\|_\infty \leq \|\kappa\|_\infty + 1$. The first term on the right-hand side of (3.33) is clearly equal to $e^{o(N)} e^{N\langle \nu, \log \mu \rangle}$. We now take a look at the second one. For $s \in \mathcal{S}$ we put

$$C_s(\ell^{(N)}) := \sum_{k:\, R<|k|\leq \varepsilon N,\, k_s>0} \ell_k^{(N)}(k_s - 1) = N\nu_s^{(N)} - \sum_{k:\, R<|k|\leq \varepsilon N,\, k_s>0} \ell_k^{(N)}. \tag{3.34}$$

Then $\frac{1}{N} C_s(\ell^{(N)}) \in [\nu_s^{(N)} - \frac{1}{R}, \nu_s^{(N)}]$ since $\sum_{k:\, R<|k|\leq \varepsilon N} \ell_k^{(N)} \leq \frac{1}{R} \sum_{k:\, R<|k|} |k| \ell_k^{(N)} \leq N/R$. Next, we apply Jensen's inequality and the fact that $x \mapsto \log x$ is concave to get that

$$\prod_{R<|k|\leq \varepsilon N} \prod_{s\in S_k} \Big(\frac{(\kappa_N k)_s}{k_s}\Big)^{(k_s-1)\ell_k^{(N)}}$$

$$= \exp\Big(\sum_{s\in\mathcal{S}} C_s(\ell^{(N)}) \sum_{R<|k|\leq \varepsilon N,\, k_s>0} \frac{\ell_k^{(N)}(k_s-1)}{C_s(\ell^{(N)})} \log \frac{(\kappa_N k)_s}{k_s}\Big)$$

$$\leq \exp\Big(\sum_{s\in\mathcal{S}} C_s(\ell^{(N)}) \log\Big(\sum_{R<|k|\leq \varepsilon N,\, k_s>0} \frac{\ell_k^{(N)}(k_s-1)}{C_s(\ell^{(N)})} \frac{(\kappa_N k)_s}{k_s}\Big)\Big)$$

$$\leq \exp\Big(N \sum_{s\in\mathcal{S}} \tfrac{1}{N} C_s(\ell^{(N)}) \log \frac{(\kappa_N \nu^{(N)})_s}{\tfrac{1}{N} C_s(\ell^{(N)})}\Big) \leq e^{o(N)} \exp\Big(N\Big\langle \nu, \log \frac{\kappa\nu}{\nu}\Big\rangle\Big) e^{N\delta_R(\nu)}$$

where

$$\delta_R(\nu) = \sum_{s \in \mathcal{S}} \nu_s \log \frac{\nu_s}{\nu_s - \frac{1}{R}} + \frac{1}{R} \sum_{s \in \mathcal{S}} \left( \log \frac{\nu_s}{(\kappa \nu)_s} \right) \wedge 0$$

which is continuous in $\nu$ and converges to 0 as $R \to \infty$.

It remains to argue that the large-$N$ exponential scale of the last term on the right-hand side of (3.33) vanishes when taking $R \to \infty$ afterwards. Recall that the choice of $r$ may depend on $k$; we will denote it as $r_k$. We first use that

$$\prod_{R < |k| \le \varepsilon N} \left( \frac{NeC}{\ell_k^{(N)}(k_{r_k})^2} \right)^{\ell_k^{(N)}} = \exp\left( N \sum_{R < |k| \le \varepsilon N} \frac{\ell_k^{(N)}}{N} \log \frac{NeC}{\ell_k^{(N)}(k_{r_k})^{2+|\mathcal{S}|}} \right)$$
$$\times \exp\left( N|\mathcal{S}| \sum_{R < |k| \le \varepsilon N} \frac{\ell_k^{(N)}}{N} \log |k| \right).$$

Abbreviating $D := \sum_{R < |k| < \varepsilon N} \ell_k^{(N)}/N$ and using Jensen's inequality we get that

$$\exp\left( N \sum_{R < |k| \le \varepsilon N} \frac{\ell_k^{(N)}}{N} \log \frac{NeC}{\ell_k^{(N)}(k_{r_k})^{2+|\mathcal{S}|}} \right) \le \exp\left( ND \log \left( \frac{eC}{D} \sum_{R < |k| \le \varepsilon N} k_{r_k}^{-2-|\mathcal{S}|} \right) \right). \tag{3.35}$$

It is easy to see that choosing $r_k$ such that $k_{r_k} = \max_{s \in \mathcal{S}} k_s$ ensures the convergence of $\sum_k k_{r_k}^{-2-|\mathcal{S}|}$ and the fact that $\sum_{|k| \ge R} k_{r_k}^{-2-|\mathcal{S}|}$ is polynomial in $R$. Further, as we remarked below (3.34), $D \le 1/R$ and therefore the right-hand side of (3.35) is bounded by $\exp(N\delta'_R)$ for some $\delta'_R$ that vanishes as $R \to \infty$. Next, we use that $R/\log R \le |k|/\log |k|$ holds on our summation area of $k$ if $R$ is large enough and therefore

$$\exp\left( N|\mathcal{S}| \sum_{R < |k| \le \varepsilon N} \frac{\ell_k^{(N)}}{N} \log |k| \right) \le \exp\left( N|\mathcal{S}| \frac{\log R}{R} \sum_{R < |k| \le \varepsilon N} \frac{\ell_k^{(N)}}{N} |k| \right)$$
$$= \exp\left( N|\mathcal{S}| \frac{\log R}{R} \right). \tag{3.36}$$

Noting that $\frac{1}{R} \log R \to 0$ as $R \to \infty$, we have shown that the last term on the right-hand side of (3.33) can be bounded by $e^{N\delta''_R}$ for some $\delta''_R$ that vanishes as $R \to \infty$. So far, we have handled the first term in the definition of $z_{R,\varepsilon N}^{(N)}(\ell^{(N)})$, and we saw that its exponential rate is not larger than the first term in $-I_{\text{Me}}(\nu)$.

Let us now handle the second and last part of $z_{R,\varepsilon N}^{(N)}(\ell^{(N)})$. We use $1 - x \le e^{-x}$ for $x = \kappa_N(r,s)/N$ as well as $\frac{1}{2} \sum_{R < k \le \varepsilon N} \ell_k^{(N)} \langle k, \kappa k \rangle \le$ $\frac{1}{2}\|\kappa\|_\infty \varepsilon \sum_{R<|k|\leq \varepsilon N} \ell_k^{(N)} |k| \leq \frac{N}{2}\|\kappa\|_\infty \varepsilon$ to get that

$$\prod_{R<|k|\leq \varepsilon N} \prod_{r,s\in \mathcal{S}} \left(1 - \frac{\kappa_N(r,s)}{N}\right)^{\frac{1}{2}k_s(N\mu_r - k_r)\ell_k^{(N)}} \leq e^{o(N)} e^{-\frac{N}{2}\langle \nu, \kappa \mu\rangle} e^{\frac{N}{2}\|\kappa\|_\infty \varepsilon}.$$

Note that the exponential rate of the last factor on the right-hand side vanishes as $\varepsilon \to 0$. Collecting all our estimates we have shown that the estimate (3.32) holds with $C(R, \varepsilon, \nu) = \delta_R(\nu) + \delta_R' + \delta_R'' + \frac{1}{2}\|\kappa\|_\infty \varepsilon$. □

### 3.5 Proof of Theorem 3.1

Here we finish now the proof of the LDP of Theorem 3.1. Let $d_\mathcal{L}$ and $d_\mathcal{A}$ be, respectively, metrics that induce the vague topologies on the state spaces $\mathcal{L}$ and $\mathcal{A}$ (see (3.3) and (3.4)). Notice that, thanks to the constraints $|c(\lambda)| \leq 1$ and $|c(\alpha)| \leq 1$, the spaces $\mathcal{L}$ and $\mathcal{A}$ endowed with the vague topologies are compact by the Bolzano–Weierstrass theorem and Fatou's lemma. Moreover the rate function $I$ is lower semicontinuous on $\mathcal{L} \times \mathcal{A}$, which makes it a good rate function. Hence, a weak LDP implies the claim of Theorem 3.1 and it suffices to prove that

$$\lim_{\delta,\rho \to 0} \lim_{N\to\infty} \frac{1}{N} \log \mathbb{P}_N\big(\mathrm{Mi}_N \in B_\delta(\lambda), \mathrm{Ma}_N \in B_\rho(\alpha)\big) = -I(\lambda, \alpha), \quad \lambda \in \mathcal{L}, \alpha \in \mathcal{A}. \tag{3.37}$$

where $B_\delta(\lambda)$ and $B_\rho(\alpha)$ are closed balls centered at $\lambda$ and $\alpha$ with radii $\delta$ and $\rho$, respectively.

Fix $\lambda \in \mathcal{L}$ and $\alpha \in \mathcal{A}$.

*Step 1: Cardinality of $\mathcal{L}_N$.* Recall the definition of $\mathcal{L}_N$ that was given in (3.13). Each $\ell \in \mathcal{L}_N$ has a unique representation as a product measure $\ell = \bigotimes_{r\in \mathcal{S}} \ell^{(r)}$, where $\ell^{(r)} = (\ell_j^{(r)})_{j\in \mathbb{N}_0}$ and $\sum_j \ell_j^{(r)} j = N\mu_r^{(N)}$ for all $r \in \mathcal{S}$. By the same argument as in [2, Lemma 3.2] there are at most $e^{o(N\mu_r^{(N)})}$ ways to choose the marginal $\ell^{(r)}$. Consequently we have

$$|\mathcal{L}_N| = e^{o(N)}.$$

Fix $R \in \mathbb{N}$ and $\varepsilon > 0$. We denote by $d_\mathcal{L}^R$ and $d_\mathcal{A}^\varepsilon$ the distance of the projections of measures on $\{k: |k| \leq R\}$ and $\{y: |y| > \varepsilon\}$, respectively. Then $d_\mathcal{L} \geq d_\mathcal{L}^R$ and $d_\mathcal{A} \geq d_\mathcal{A}^R$. For $\delta > 0$ and $\rho > 0$ denote by $\mathcal{L}_N^{(R,\varepsilon)}(\delta, \rho)$ the set of all $\ell \in \mathcal{L}_N$ with $d_\mathcal{L}^R(\frac{1}{N}\ell, \lambda) < \delta$ and $d_\mathcal{A}^\varepsilon(\ell_{\lfloor N\cdot \rfloor}, \alpha) < \rho$. Note that

$$\mathbb{P}\big(\mathrm{Mi}_N \in B_\delta(\lambda), \mathrm{Ma}_N \in B_\rho(\alpha)\big) \leq \sum_{\ell \in \mathcal{L}_N^{(R,\varepsilon)}(\delta,\rho)} \mathbb{P}\left(N\mathrm{Mi}_N = \ell\right).$$

According to the preceding, also

$$|\mathcal{L}_N^{(R,\varepsilon)}(\delta,\rho)| \leq |\mathcal{L}_N| = e^{o(N)}.$$

*Step 2: Case* $c(\lambda) + c(\alpha) \not\leq \mu$. Assume that there is some $r \in \mathcal{S}$ such that $c_r(\lambda) + c_r(\alpha) > \mu_r$. Then it is easy to see that $\mathcal{L}_N(\delta,\rho) = \emptyset$ for sufficiently large $N$ and hence

$$\lim_{N\to\infty} \frac{1}{N} \log \mathbb{P}\big(\mathrm{Mi}_N \in B_\delta(\lambda), \mathrm{Ma}_N \in B_\rho(\alpha)\big) = -\infty = -I(\lambda,\alpha),$$

which is proved with the same argument as in [2, Lemma 3.7].

*Step 3: Proof of the upper bound in* (3.37) *for* $c(\lambda) + c(\alpha) \leq \mu$. For any $R \in \mathbb{N}$ and any $\varepsilon \in (0,1]$ we have by Step 1 and Lemma 3.7 that

$$\begin{aligned}
\mathbb{P}\big(\mathrm{Mi}_N \in B_\delta(\lambda), \mathrm{Ma}_N \in B_\rho(\alpha)\big) &\leq e^{o(N)} \sup_{\ell \in \mathcal{L}_N^{(R,\varepsilon)}(\delta,\rho)} \mathbb{P}(N\mathrm{Mi}_N = \ell) \\
&\leq e^{o(N)} \sup_{\ell \in \mathcal{L}_N^{(R,\varepsilon)}(\delta,\rho)} z_{0,R}^{(N)}(\ell) z_{R,\varepsilon N}^{(N)}(\ell) z_{\varepsilon N, N}^{(N)}(\ell).
\end{aligned} \tag{3.38}$$

For $\tilde{\lambda} \in \mathcal{L}$ and $\tilde{\alpha} \in \mathcal{A}$ we define

$$c_s^{(R)}(\tilde{\lambda}) := \sum_{k:\, |k| \leq R} \tilde{\lambda}_k k_s \quad \text{and} \quad c_s^{(\varepsilon)}(\tilde{\alpha}) := \int_{\{|y|>\varepsilon\}} y_s\, \tilde{\alpha}(\mathrm{d}y), \qquad s \in \mathcal{S}. \tag{3.39}$$

We require that $\varepsilon \in (0,1] \setminus \{|y|: y \in \mathrm{supp}(\alpha)\}$ (recall that $\mathrm{supp}(\alpha)$ is countable). This is a prerequisite to apply Lemma 3.9 to $\tilde{\alpha}$. Now, applying Lemmas 3.8, 3.9(1) and 3.10 to the right-hand side of (3.38) we get

$$\begin{aligned}
&\limsup_{N\to\infty} \frac{1}{N} \log \mathbb{P}\big(\mathrm{Mi}_N \in B_\delta(\lambda), \mathrm{Ma}_N \in B_\rho(\alpha)\big) \\
&\leq \sup_{\substack{\tilde{\lambda}:\, d_\mathcal{L}^R(\tilde{\lambda},\lambda)<\delta \\ \tilde{\alpha}:\, d_\mathcal{A}^\varepsilon(\tilde{\alpha},\alpha)<\rho}} \Big( - I_{\mathrm{Mi}}^{(R)}(\tilde{\lambda}) - I_{\mathrm{Ma}}^{(\varepsilon)}(\tilde{\alpha}) - I_{\mathrm{Me}}\big(\tilde{\nu}^{(R,\varepsilon)}\big) + C(R,\varepsilon, \tilde{\nu}^{(R,\varepsilon)})\Big),
\end{aligned}$$

where we used the cut-off versions of the rate functions defined in Lemmas 3.8, 3.9 and 3.10. Also, we abbreviated $\tilde{\nu}^{(R,\varepsilon)} := \mu - c^{(R)}(\tilde{\lambda}) - c^{(\varepsilon)}(\tilde{\alpha})$ and recall that $C(R,\varepsilon, \tilde{\nu}^{(R,\varepsilon)})$ is the term given in Lemma 3.10. Note that the functions $I_{\mathrm{Mi}}^{(R)}$ and $c^{(R)}$ are continuous (in any point) and the functions $I_{\mathrm{Ma}}^{(\varepsilon)}$ and $c^{(\varepsilon)}$ are continuous in $\alpha$ due to our requirement that $\varepsilon \notin \{|y|: y \in \mathrm{supp}(\alpha)\}$. This also implies that $\tilde{\nu}^{(R,\varepsilon)}$ converges to $\nu^{(R,\varepsilon)} := \mu - c^{(R)}(\lambda) - c^{(\varepsilon)}(\alpha)$, as $\delta, \rho \to 0$, and due to the continuity of $I_{\mathrm{Me}}$ and

$C(R, \varepsilon, \cdot)$ the respective terms converge. Altogether, we get that

$$\limsup_{\delta,\rho \to 0} \limsup_{N\to\infty} \frac{1}{N} \log \mathbb{P}\big(\mathrm{Mi}_N \in B_\delta(\lambda), \mathrm{Ma}_N \in B_\rho(\alpha)\big)$$
$$\leq -I_{\mathrm{Mi}}^{(R)}(\lambda) - I_{\mathrm{Ma}}^{(\varepsilon)}(\alpha) - I_{\mathrm{Me}}(\nu^{(R,\varepsilon)}) + C(R, \varepsilon, \nu^{(R,\varepsilon)})$$

Observe that the right-hand side converges to $-I(\lambda, \alpha)$, if we let $R \to \infty$ and $\varepsilon \to 0$, which proves the upper bound in (3.37). Notice that the requirement that $\varepsilon \notin \{|y|: y \in \mathrm{supp}(\alpha)\}$ is not a problem since $\mathrm{supp}(\alpha)$ is countable.

*Step 4: Construction of a recovery sequence.* In this step, we prepare for the proof of the lower bound in (3.37) (see Step 5) by constructing an almost optimal sequence of $\ell$'s. We handle here only the case that $\kappa$ is irreducible; see Remark 3.11 for hints how to handle the case of a reducible $\kappa$. We may assume that $c(\lambda) + c(\alpha) \leq \mu$, since the rate function is equal to $\infty$ otherwise. For the same reason, we also may assume that the mesoscopic mass $\nu := \mu - c(\lambda) - c(\alpha)$ satisfies $\nu \ll \kappa\nu$, since otherwise $I_{\mathrm{Me}}(\nu) = \infty$.

For $R \in \mathbb{N}$ and $\varepsilon \in (0, 1]$ we construct a suitable recovery sequence $\ell^{(N)} = \ell^{(N)}(R, \varepsilon)$ that will turn out in Step 5 as asymptotically optimal. To this end, we construct it in such a way that it will put all mesoscopic mass $\nu := \mu - c(\lambda) - c(\alpha)$ into several components that can all be described by the same configuration $k^{(\mathrm{Me},N)}$ and which are actually on the lower end of the macroscopic scale, such that (3.20) of Theorem 3.6 can be applied. Our construction also ensures that all components of macroscopic scale have a strictly positive probability of being connected. Denote by $\mathbb{1}$ the element of $\mathbb{N}_0^S$ that is equal to 1 in each entry. Define $k^{(\mathrm{Me},N)} := \lfloor N\varepsilon\nu \rfloor + \mathbb{1}$. Using the representation $\alpha = \sum_i \delta_{y^{(i)}}$ put $k^{(i,N)} := \lfloor Ny^{(i)}\rfloor + \mathbb{1}$ for $i \in J_\varepsilon(\alpha) := \{i: |y^{(i)}| > \varepsilon\}$.

Let us check now that we can actually apply Theorem 3.6 to $p_N(k^{(i,N)})$ for all $i \in J_\varepsilon(\alpha) \cup \{\mathrm{Me}\}$. To check that $\tau(k^{(i,N)}) > 0$ for all $i \in J_\varepsilon(\alpha) \cup \{\mathrm{Me}\}$ note the following: Depending on the support of $y^{(i)}$ we might not have that $\kappa$ is irreducible with respect to $y^{(i)}$. However, the vectors $k^{(i,N)}$ have support on the full type space $\mathcal{S}$ for all $N \in \mathbb{N}$ by construction. Consequently, $\tau(k^{(i,N)}) > 0$ is implied by the fact that $\kappa$ is irreducible with respect to $\mu$. Secondly, since $\nu \ll \kappa\nu$, we also have that $y^{(i)} \ll \kappa y^{(i)}$ for all $i \in J_\varepsilon(\alpha)$: if not, $I_{\mathrm{Ma}}(\alpha) = \infty$ by definition and the lower bound $-\infty$ in (3.37) trivially holds. By construction it holds for all $i \in J_\varepsilon(\alpha)$ that $k_s^{(i,N)} = 1$ for $s \notin \mathrm{supp}(y^{(i)})$ and for all $N \in \mathbb{N}$, which ensures the boundedness condition. The same holds for $k^{(\mathrm{Me},N)}$.

Define now

$$\ell_k^{(N)} := \begin{cases} \lfloor \lambda_k N \rfloor & \text{if } 2 \leq |k| \leq R \\ 0 & \text{if } R < |k| \leq \varepsilon N \text{ and } k \neq k^{(\mathrm{Me},N)} \\ \lfloor \varepsilon^{-1} \rfloor & \text{for } k = k^{(\mathrm{Me},N)} \\ \#\{i \in J_\varepsilon(\alpha): k = k^{(i,N)}\} & \text{if } \varepsilon N < |k| \\ N\mu_r^{(N)} - \sum_{m:2\leq|m|} \ell_m^{(N)} m_r & \text{if } k = \mathbf{e}_r \text{ for any } r \in \mathcal{S}, \end{cases} \quad (3.40)$$

where the last line ensures that $\ell^{(N)} \in \mathcal{L}_N$. Observe that

$$\lim_{N\to\infty} \frac{1}{N}\ell^{(N)}_{\mathbf{e}_r} - \lambda_{\mathbf{e}_r} = c_r(\lambda) - c_r^{(R)}(\lambda) + \nu_r - \lfloor \varepsilon^{-1}\rfloor(\varepsilon\nu_r)$$
$$+ c_r(\alpha) - c_r^{(\varepsilon)}(\alpha) \text{ for } r \in \mathcal{S}$$
$$\lim_{N\to\infty} \frac{1}{N}\ell^{(N)}_k - \lambda_k = 0 \qquad \text{for } 2 \leq |k| \leq R,$$

which implies that $\lim_{R\to\infty, \varepsilon\to 0} \lim_{N\to\infty} d_{\mathcal{L}}(\frac{1}{N}\ell^{(N)}, \lambda) = 0$. Similarly, $\lim_{\varepsilon\to 0} \lim_{N\to\infty} d_{\mathcal{A}}(\alpha^{(N)}, \alpha) = 0$, where $\alpha^{(N)} = \sum_k \ell^{(N)}_k \delta_{\frac{k}{N}}$.

*Step 5: Proof of the lower bound in* (3.37). Now we finish the proof of the lower bound by showing that the recovery sequence $(\ell^{(N)})_{N\in\mathbb{N}}$ that we constructed in Step 4 is giving the right asymptotics.

Fix $\delta, \rho > 0$. Then by choosing $R \in \mathbb{N}$ large enough and $\varepsilon > 0$ small enough, we have that $\frac{1}{N}\ell^{(N)} \in B_\delta(\lambda)$ and $\alpha^{(N)} \in B_\rho(\alpha)$ for all but finitely many $N \in \mathbb{N}$. Hence, Lemma 3.7 implies that

$$\mathbb{P}\left(\mathrm{Mi}_N \in B_\delta(\lambda), \mathrm{Ma}_N \in B_\rho(\alpha)\right) \geq \mathbb{P}(N\mathrm{Mi}_N = \ell^{(N)})$$
$$= e^{o(N)} z^{(N)}_{0,R}(\ell^{(N)}) z^{(N)}_{R,\varepsilon N}(\ell^{(N)}) z^{(N)}_{\varepsilon N, N}(\ell^{(N)}).$$

Next, we want to use Lemmas 3.8 and 3.9(2) to get the exponential rates for the $z$-terms. Note that we do not have $\lim_{N\to\infty} \frac{1}{N}\ell^{(N)} = \lambda$ on $\{k\colon |k| \leq R\}$, so we will instead apply Lemma 3.8 to the sequence $\lfloor \lambda N \rfloor$, which gives

$$\lim_{N\to\infty} \frac{1}{N} \log z^{(N)}_{0,R}(\lfloor\lambda N\rfloor) = -I^{(R)}_{\mathrm{Mi}}(\lambda)$$

Using the definition (3.21), it is easy to verify that

$$\lim_{N\to\infty} \frac{1}{N}\log z^{(N)}_{0,R}(\ell^{(N)}) = \lim_{N\to\infty} \frac{1}{N}\log z^{(N)}_{0,R}(\lfloor\lambda N\rfloor) + \lim_{N\to\infty} \frac{1}{N}\sum_{r\in\mathcal{S}} \log \frac{z^{(N)}(\ell^{(N)}_{\mathbf{e}_r}, \mathbf{e}_r)}{z^{(N)}(\lfloor\lambda_{\mathbf{e}_r} N\rfloor, \mathbf{e}_r)}$$
$$= -I^{(R)}_{\mathrm{Mi}}(\lambda) + \widetilde{C}(R,\varepsilon).$$

for some constant $\widetilde{C}(R,\varepsilon)$ vanishing as $R \to \infty$ and $\varepsilon \to 0$. For the mesoscopic part, we use the asymptotic formula (3.20) of Theorem 3.6 for $p_N(k^{(\mathrm{Me},N)})$. We proved in Step 4 that we can actually apply Theorem 3.6 to $p_N(k^{(\mathrm{Me},N)})$. This, together with Stirling's formula, gives us

$$\lim_{N\to\infty} \frac{1}{N}\log z^{(N)}_{R,\varepsilon N}(\ell^{(N)}) = \lfloor\varepsilon^{-1}\rfloor\left\langle\varepsilon\nu, \log\frac{(1-e^{-\kappa(\varepsilon\nu)})\mu}{\varepsilon\nu}\right\rangle - \frac{1}{2}\lfloor\varepsilon^{-1}\rfloor\langle\varepsilon\nu, \kappa(\mu - \varepsilon\nu)\rangle.$$

The right-hand side clearly converges to $-I_{\mathrm{Me}}(\nu)$ as $\varepsilon \to 0$. The analysis of the macroscopic part in Step 4 shows that we can directly use Lemma 3.9(2) to see that

$$\lim_{N \to \infty} \frac{1}{N} \log z_{\varepsilon N, N}^{(N)}(\ell^{(N)}) = -I_{\mathrm{Ma}}^{(\varepsilon)}(\alpha).$$

Taking the limits $R \to \infty$ and $\varepsilon \to 0$ we obtain the lower bound in (3.37). This finishes the proof of Theorem 3.1.

*Remark 3.11* (*The reducible case: proof of Theorem* 3.2) Recall from Sect. 2.3 that in the case that $\kappa$ is reducible, one can decompose the type space $\mathcal{S} = \bigcup_j \mathcal{S}_j$ in such a way that $\kappa$ restricted to $\mathcal{S}_j \times \mathcal{S}_j$ is irreducible and $\kappa|_{\mathcal{S}_i \times \mathcal{S}_j} = 0$ for $i \neq j$. Assume that $\alpha \in \mathcal{A}$ is connectable as defined in Sect. 2.3: for each $y \in \mathrm{supp}(\alpha)$ there is some $j$ such that $\mathrm{supp}(y) \subset \mathcal{S}_j$. Then we can argue that the lower bound in the proof of (3.37) holds. For any $j$ we can construct a recovery sequence as before, where we approximate the macroscopic components by defining $k^{(i,N)} := \lfloor N y^{(i)} \rfloor + \mathbb{1}_{\mathcal{S}_j}$ for any $y^{(i)} \in \mathrm{supp}(\alpha)$ to ensure that $\tau(k^{(i,N)}) > 0$. In the same way we approximate the mesoscopic mass by defining $\nu^{(j)} := \nu \mathbb{1}_{\mathcal{S}_j}$ and $k^{(\mathrm{Me},j,N)} := \lfloor N \varepsilon \nu^{(j)} \rfloor + \mathbb{1}_{\mathcal{S}_j}$. The rest will work out as in the proof above. The only additional observation needed is that for any $j$ we have that $(\kappa \nu^{(j)})_s = (\kappa \nu)_s$ if $s \in \mathcal{S}_j$.

On the other hand, if $\alpha \in \mathcal{A}$ is not connectable, then we can argue that the upper bound of the LDP holds with $I(\lambda, \alpha) = \infty$ by additionally requiring that $\kappa_N = \kappa$ for all but finitely many $N \in \mathbb{N}$: Indeed, there exists $y \in \mathrm{supp}(\alpha)$ such that for any sequence $k^{(N)}$ with $\lim_{N \to \infty} \frac{k^{(N)}}{N} = y$ we have that $\tau(k^{(N)}) = 0$ and by using that $\kappa_N = \kappa$ we get that $p_N(k^{(N)}) \leq \tau(k^{(N)}) N^{-(|k|-1)} = 0$, which gives that $I(\lambda, \alpha) = \infty$.

## 4 Proofs of the results from Sect. 3.3

The aim of this section is to prove the results about the connection probabilities that where formulated in Sect. 3.3, namely Lemmas 3.4, 3.5 and Theorem 3.6. In our presentation we focus on deriving Theorem 3.6, while the lemmas will be a by-product of the procedure. Their proofs can be found after Lemma 4.8.

The idea for the proof of Theorem 3.6 is to construct a sequence of graphs that contain with high probability a macroscopic component with the desired type configuration. In order to understand how to choose the parameters correctly we have to understand the characteristic equation (2.4) and how it emerges from the generating function of weighted trees. This will be done in Sect. 4.1. In Sect. 4.2 we will collect estimates that provide the link between the weighted trees and the connection probabilities of the graph and finally prove Theorem 3.6, which is a consequence of Lemmas 4.10 and 4.13.

### 4.1 The characteristic equation and tree combinatorics

In this section we discuss a power series representation of the solution of fixed point equation (2.4). This will be crucial both in the proof of Theorem 3.6 and in the analysis of the minimizers of the rate function in Sect. 6.

We rewrite the equation in the following way. Fix $\nu = (\nu_s)_{s \in \mathcal{S}} \in [0, \infty)^{\mathcal{S}}$. We say that $\nu^* = (\nu_s^*)_{s \in \mathcal{S}}$ is a solution to the characteristic equation with respect to $\nu$ if

$$\nu_s^* e^{-(\kappa \nu^*)_s} = \nu_s e^{-(\kappa \nu)_s}, \quad s \in \mathcal{S}. \tag{4.1}$$

Recall that this is equivalent to the characteristic equation (2.11) that was studied in [7] via the substitution $\nu = \mu$ and $\rho = 1 - \frac{d\nu^*}{d\nu}$. However, we will need Eq. (4.1) on different occasions and for several choices of $\nu$. Note that any solution $\nu^*$ to (4.1) is necessarily non-negative and for any $s \in \mathcal{S}$ we have that $\nu_s^* > 0$ if and only if $\nu_s > 0$. Also note that $\nu$ itself is always a solution. Whether there exists a non-trivial solution $\nu^*$ (i.e., $\nu^* \neq \nu$) or not will turn out to depend on the quantity $\Sigma(\kappa, \nu)$, introduced in (2.2) and (2.3):

$$\Sigma(\kappa, \nu) = \|T_{\kappa,\nu}\|_\nu = \sup\left\{\|T_{\kappa,\nu} f\|_\nu \colon f \in [0, \infty)^{\mathcal{S}}, \|f\|_\nu = 1\right\}$$
$$= \sup\left\{\left(\sum_{r \in \mathcal{S}} \nu_r \left(\sum_{s \in \mathcal{S}} \kappa(r,s) f(s) \nu_s\right)^2\right)^{1/2} \colon f \in [0, \infty)^{\mathcal{S}}, \sum_{r \in \mathcal{S}} \nu_r f(r)^2 \leq 1\right\}, \tag{4.2}$$

where we write $\|\cdot\|_\nu$ for the norm on $L^2(\mathbb{R}^{\mathcal{S}}, \nu)$ and also for the corresponding operator norm. We note that $\Sigma(\kappa, \nu)$ is equal to the spectral radius of the matrix $T_{\kappa,\nu} = (\kappa(r,s)\nu_s)_{r,s \in \mathcal{S}}$, as is seen from an elementary analysis of (4.2), also using Frobenius eigenvalue theory. Indeed, the variational equations for the maximizer $f$ in (4.2) (with $\sum_r \nu_r f(r)^2 = 1$ instead of $\leq 1$) say that $f$ is a positive eigenvector of the matrix $T_{\kappa,\nu}^2$. Since also the (up to positive multiples, unique) positive eigenvector of $T_{\kappa,\nu}$ is a positive eigenvector of $T_{\kappa,\nu}^2$, we have that $f$ is the Frobenius eigenvector of $T_{\kappa,\nu}$. The corresponding Frobenius eigenvalue, i.e., the spectral radius, is equal to $\Sigma(\kappa, \nu)$.

Another elementary application of Frobenius eigenvalue theory yields that the map $\nu \mapsto \Sigma(\kappa, \nu)$ is non-decreasing with respect to componentwise ordering.

Now we cite the results from [7] regarding the solutions of the characteristic equation.

**Lemma 4.1** *Let $\nu \in [0, \infty)^{\mathcal{S}}$.*

(i) *If $\Sigma(\kappa, \nu) \leq 1$, then the only solution $\nu^*$ to the characteristic equation (4.1) satisfying $\nu^* \leq \nu$ is given by $\nu^* = \nu$.*
(ii) *If $\Sigma(\kappa, \nu) > 1$, then there exists a solution $\nu^*$ to (4.1) that satisfies $\nu^* \leq \nu$ and $\nu^* \neq \nu$. If additionally $\kappa$ is irreducible, then $\nu^*$ is the only solution to (4.1) that satisfies $\nu^* \leq \nu$ and $\nu^* \neq \nu$. Further, $\Sigma(\kappa, \nu^*) < 1$.*

**Proof** See Theorem 6.2 and Theorem 6.7 in [7] and substitute $\rho = 1 - \frac{dv^*}{dv}$. □

Our aim is to verify the following.

**Proposition 4.2** *Let $v \in [0, \infty)^{\mathcal{S}}$. Then for any $r \in \mathcal{S}$*

$$\sum_{k \in \mathbb{N}_0^{\mathcal{S}}} \tau(k) k_r \prod_{s \in \mathcal{S}} \frac{(v_s e^{-(\kappa v)_s})^{k_s}}{k_s!} = v_r^* \qquad (4.3)$$

*where $v^*$ is the smallest solution to (4.1).*

The result of Proposition 4.2 will be used in Sect. 4.2 to derive the asymptotics of the probability that a macroscopic set of vertices is connected. Further, it will be used in Sect. 6 to optimize the microscopic part of the rate function. The left-hand side of Eq. (4.3), when divided by $v_r$, is equal to the extinction probability of the branching process that we mentioned in Sect. 2.2. This observation, together with additional results from [7], is already enough to prove Proposition 4.2. However, we will provide a different proof that mainly uses the structure of the power series and additionally gives us a refined control of the convergence of the series, thanks to the estimates given in Lemma 4.5. They will be used later to bound probabilistic terms that we derive from components of mesoscopic sizes, but also in Sect. 6 where we need uniform convergence for a certain family of power series.

We now prepare for the proof of Proposition 4.2. We define a function $\Gamma = (\Gamma_r)_{r \in \mathcal{S}} \colon [0, \infty)^{\mathcal{S}} \to [0, \infty)^{\mathcal{S}}$ by putting, for $\theta = (\theta_s)_{s \in \mathcal{S}} \in [0, \infty)^{\mathcal{S}}$,

$$\Gamma_r(\theta) := \sum_{k \in \mathbb{N}_0^{\mathcal{S}}} \tau(k) k_r \prod_{s \in \mathcal{S}} \frac{\theta_s^{k_s}}{k_s!} \in [0, \infty] \quad \text{for } r \in \mathcal{S}. \qquad (4.4)$$

The idea of the proof is to show that $\theta \mapsto \Gamma(\theta)$ is the inverse of the function $v \mapsto \theta(v) := v e^{-\kappa v}$ on the domain $\{v[0, \infty)^{\mathcal{S}} \colon \Sigma(\kappa, v) \leq 1\}$. It turns out that this is an easy example of a technique known as Lagrange inversion, where directed trees are used to extract the variables of a power series, see [25, 26] for more general results on this topic. The relation of our combinatorial quantity $\tau$ to directed trees is given in the following lemma. The second statement will be useful to derive a criterion for analyticity of the power series defined in (4.4).

**Lemma 4.3** *Let $k \in \mathbb{N}_0^{\mathcal{S}}$ and $r \in \mathcal{S}$. Then the following holds.*

1. *Let $\vec{\mathcal{T}}_r(k)$ denote the set of directed trees with vertex set $\{1, \ldots, |k|\}$, root of type $r$ and edges directed away from the root. We use the convention that $\vec{\mathcal{T}}_r(0) = \emptyset$. Then*

$$\tau(k) k_r = \sum_{T \in \vec{\mathcal{T}}_r(k)} \prod_{(i,j) \in E(T)} \kappa(x_i, x_j). \qquad (4.5)$$

2. Write $S_0 := \mathrm{supp}(k) \subset \mathcal{S}$ and $\overrightarrow{\mathcal{T}}_r(S_0)$ for the set of all directed trees with vertex set $S_0$, root vertex given by $r$ and edges directed away from the root. We set $\overrightarrow{\mathcal{T}}_r(\emptyset) = \emptyset$. Then

$$\tau(k)k_r = \left( \prod_{s \in S_0} (\kappa k)_s^{k_s - 1} \right) \times \Delta_r(k) \qquad (4.6)$$

where

$$\Delta_r(k) = \sum_{A \in \overrightarrow{\mathcal{T}}_r(S_0)} \prod_{(s,s') \in E(A)} \kappa(s, s')k_s. \qquad (4.7)$$

**Proof** (1) Given some vertex $i \in [|k|]$ of type $r$, we can turn any undirected tree from $\mathcal{T}(k)$ uniquely into a directed tree that has vertex $i$ as its root by giving to each edge the direction away from the root. For each undirected tree $T$ in $\mathcal{T}(k)$ there are $k_r$ ways to choose the root, hence the weight of $T$ appears $k_r$ times on the right-hand side of (4.5).

(2) We use the formula derived in [4, Theorem 2] with $x_{s,s',i} := \kappa(s,s') = \kappa(s',s)$. Note that in [4] the edges of the trees are directed towards the root. Adapted to our notation their formula reads as

$$\sum_{T \in \overrightarrow{\mathcal{T}}_r(k)} \prod_{s,s' \in S_0} \kappa(s, s')^{\#\{(i,j): \, x_i = s, x_j = s'\}} = \left( \prod_{s \in S_0} (\kappa k)_s^{k_s - 1} \right) \times \Delta_r(k). \qquad (4.8)$$

By (4.5) we have that

$$\tau(k)k_r = \sum_{T \in \overrightarrow{\mathcal{T}}_r(k)} \prod_{(i,j) \in E(T)} \kappa(x_i, x_j) = \sum_{T \in \overrightarrow{\mathcal{T}}_r(k)} \prod_{s,s' \in S_0} \kappa(s, s')^{\#\{(i,j): \, x_i = s, x_j = s'\}},$$

so together with (4.8) we have proven equation (4.6). Note that in the case $r \notin S_0$, we have that $\overrightarrow{\mathcal{T}}_r(S_0) = \emptyset$, so both sides of (4.8) are 0. □

**Lemma 4.4** *If $\Gamma$ is analytic in $\theta$ (i.e., $\Gamma_r$ is analytic in $\theta$ for all $r \in \mathcal{S}$), then for all $r \in \mathcal{S}$*

$$\Gamma_r(\theta) = \theta_r \exp\left((\kappa \Gamma(\theta))_r\right). \qquad (4.9)$$

**Proof** Using formula (4.5) we can rewrite the power series using directed trees. Given a tree $T$ in $\overrightarrow{\mathcal{T}}_r(k)$ and some $h \in \mathbb{N}_0^{\mathcal{S}}$ we say that a vertex $i \in [|k|]$ is of parent-type $(r, h)$ if $i$ is of type $r$ and points to exactly $h_s$ vertices of type $s$ for any $s \in \mathcal{S}$. In that case, its weight is defined by $w(i) = \prod_{s \in \mathcal{S}} \kappa(r, s)^{h_s} =: g_{r,h}$. The weight of the tree $T$ is defined by $w(T) := \prod_{i=1}^{|k|} w(i)$. Defining $g_r(x) = \sum_{h \in \mathbb{N}_0^{\mathcal{S}}} g_{r,h} \prod_{s \in \mathcal{S}} \frac{x_s^{h_s}}{h_s!}$,

for $x \in \mathbb{R}^{\mathcal{S}}$, we can apply the multinomial theorem to see that $g_r(x) = \exp\{(\kappa x)_r\}$. Althogether, we have that

$$\Gamma_r(\theta) = \sum_{k \in \mathbb{N}_0^{\mathcal{S}}} \left( \sum_{T \in \vec{\mathcal{T}}_r(k)} w(T) \right) \prod_{s \in \mathcal{S}} \frac{\theta_s^{k_s}}{k_s!},$$

which is an exponential generating function, evaluated in $\theta \in \mathbb{R}^{\mathcal{S}}$ and shall be understood as a formal power series. We now cite a fact from [23] (see the proof of Theorem 1):

$$\Gamma_r(\theta) = \theta_r \sum_{h \in \mathbb{N}_0^{\mathcal{S}}} g_{r,h} \prod_{s \in \mathcal{S}} \frac{(\Gamma_s(\theta))^{h_s}}{h_s!} = \theta_r g_r(\Gamma(\theta)), \quad r \in \mathcal{S}, \quad (4.10)$$

where all equations denote equality of formal power series, i.e., we have equality of the coefficients of $\theta_s^{k_s}$. In particular, $\Gamma_r(\theta)$ might be infinite. However, by our assumption that $\Gamma$ is analytic in $\theta$, Eq. (4.10) holds in the usual sense, i.e., all series in (4.10) converge absolutely. Recalling that $g_r(x) = \exp\{(\kappa x)_r\}$, we have verified that (4.9) holds.

The idea behind proving (4.10) is to decompose the set of trees with root of type $r$ into sets of trees with root of parent-type $(r, h)$ for each $h \in \mathbb{N}_0^{\mathcal{S}}$. For fixed $h \in \mathbb{N}_0^{\mathcal{S}}$, one then decomposes a tree with root of parent-type $(r, h)$ into a single vertex of type $r$ and exactly $h_s$ many trees with root of type $s$ for each $s \in \mathcal{S}$. By studying how this decomposition affects the exponential generating function, one obtains the formula. □

Next, we derive a criterion for analyticity of the map $\theta \mapsto \Gamma(\theta)$, i.e., about the domain of convergence of the corresponding power series. We need to introduce the quantity

$$\chi(\kappa, \theta) = \inf \left\{ \left\langle v, \log \frac{v}{(\kappa v)\theta} \right\rangle : v \in \mathcal{M}_1(\mathcal{S}), v \ll \theta \right\}, \quad (4.11)$$

where $\mathcal{M}_1(\mathcal{S})$ denotes the set of probability measures on $\mathcal{S}$. One can easily see that $\theta \mapsto \chi(\kappa, \theta)$ is lower semi-continuous.

**Lemma 4.5** *1. For any $\theta \in [0, \infty)^{\mathcal{S}}$ and $r \in \mathcal{S}$,*

$$\sum_{k \in \mathbb{N}_0^{\mathcal{S}}: |k|=n} \tau(k) k_r \prod_{s \in \mathcal{S}} \frac{\theta_s^{k_s}}{k_s!} = e^{o(n)} e^{-n[\chi(\kappa,\theta)-1]}, \quad n \to \infty. \quad (4.12)$$

2. *Fix $\theta \in [0, \infty)^{\mathcal{S}}$. If $\chi(\kappa, \theta) > 1$, then for any $r \in \mathcal{S}$ the series $\Gamma_r$ defined in (4.4) is analytic in $\theta$. If on the other hand $\chi(\kappa, \theta) < 1$, then $\Gamma_r(\theta)$ diverges.*
3. *For arbitrary $v \in (0, \infty)^{\mathcal{S}}$ and $\theta_s(\kappa, v) = v_s e^{-(\kappa v)_s}$, $s \in \mathcal{S}$, we have that*

$$\chi(\kappa, \theta(\kappa, v)) = \Sigma(\kappa, v) - \log \Sigma(\kappa, v) \geq 1 \quad (4.13)$$

*with equality if and only if $\Sigma(\kappa, \nu) = 1$.*

**Proof** (1) We will use the identity (4.6) from Lemma 4.3. Let $n \in \mathbb{N}$. We write $\mathcal{M}_1^{(n)}(\mathcal{S})$ for the set of all probability measures $v$ on $\mathcal{S}$ satisfying $nv \in \mathbb{N}_0^{\mathcal{S}}$. For any $k \in \mathbb{N}_0^{\mathcal{S}}$ with $|k| = n$ we substitute $k = nv$ to get

$$\sum_{k \in \mathbb{N}_0^{\mathcal{S}}: |k|=n} \tau(k) k_r \prod_{s \in \mathcal{S}} \frac{\theta_s^{k_s}}{k_s!} = \sum_{v \in \mathcal{M}_1^{(n)}(\mathcal{S})} \left( \prod_{s \in \text{supp}(v)} (n(\kappa v)_s)^{nv_s - 1} \right) \prod_{s \in \mathcal{S}} \frac{\theta_s^{nv_s}}{(nv_s)!} \Delta_r(nv)$$

$$= \sum_{v \in \mathcal{M}_1^{(n)}(\mathcal{S})} e^{o(n)} n^{-|\text{supp}(v)|} e^{n \left[ 1 - \sum_s v_s \log \frac{v_s}{(\kappa v)_s \theta_s} \right]} \Delta_r(nv),$$

where we used Stirlings formula $N! = N^N e^{-N} e^{o(N)}$ as $N \to \infty$ and assumed throughout that $k \ll \theta$ and $v \ll \theta$ respectively. Note that $\left| \mathcal{M}_1^{(n)}(\mathcal{S}) \right| = e^{o(n)}$ and that the terms $n^{-|\mathcal{S}|}$ and $\Delta_r(nv)$ are only polynomial in $n$, thus also of order $e^{o(n)}$. Collecting everything, we arrive at (4.12).

(2) As a consequence of (1) we get the following. If $\chi(\kappa, \theta) > 1$, then the series $\Gamma_r(\theta)$ converges and since the mapping $\theta \mapsto \chi(\kappa, \theta)$ is lower-semicontinuous, we get that $\Gamma_r$ is analytic in $\theta$. If $\chi(\kappa, \theta) < 1$, then $\Gamma_r(\theta)$ diverges.

(3) Put $\phi(x) = x \log x$ for $x \geq 0$. Inserting the definition $\theta_s(v) = v_s e^{-(\kappa v)_s}$, $s \in \mathcal{S}$, we can estimate for any $v \in \mathcal{M}_1(\mathcal{S})$

$$\left\langle v, \log \frac{v}{(\kappa v) \theta(v)} \right\rangle = \left\langle v, \log \frac{v}{(\kappa v) v} \right\rangle + \langle v, \kappa v \rangle$$

$$= \langle v, \kappa v \rangle \sum_{r \in \mathcal{S}} \frac{(\kappa v)_r v_r}{\langle v, \kappa v \rangle} \phi \left( \frac{v_r}{(\kappa v)_r v_r} \right) + \langle v, \kappa v \rangle$$

$$\geq \langle v, \kappa v \rangle \phi \left( \frac{\sum_r v_r}{\langle v, \kappa v \rangle} \right) + \langle v, \kappa v \rangle = -\log \langle v, \kappa v \rangle + \langle v, \kappa v \rangle$$

where the estimate is due to Jensen's inequality applied to the convex function $\phi$. Further, we used the fact that $v$ is a probability measure and the symmetry of $\kappa$. Since $\phi$ is even strictly convex, the application of Jensen's inequality gives an equality if and only if $v$ is such that the map $\mathcal{S} \ni r \mapsto v_r / ((\kappa v)_r v_r)$ is constant, i.e., if there is some $a \in \mathbb{R}$ such that $(\kappa v)_r = a \frac{v_r}{v_r}$ for any $r \in \mathcal{S}$. Clearly, $a > 0$. It follows that $\langle v, \kappa v \rangle = a$ and that $w = (v_r / v_r)_{r \in \mathcal{S}}$ is an eigenvector of the matrix $T_{\kappa, v} = (\kappa(r, s) v_s)_{r, s \in \mathcal{S}}$ with eigenvalue $a$. Our assumption $v > 0$ and the irreducibility of $\kappa$ imply that $T_{\kappa, v}$ is also irreducible. Hence the Perron–Frobenius theorem gives that, if $a$ would be smaller than the spectal radius $\Sigma(\kappa, v)$ of $T_{\kappa, v}$, then $w$ could not be non-negative, implying that $v = (w_r v_r)_r$ would not be in $\mathcal{M}_1(\mathcal{S})$. Hence, $a$ is necessarily equal to $\Sigma(\kappa, v)$ and Eq. (4.13) holds. Since $x - 1 \geq \log x$ holds for any $x \geq 0$ with equality if and only if $x = 1$, the rest of the statement in (3) holds. □

Now, we can give the proof of Proposition 4.2.

**Proof of Proposition 4.2** We will proceed in three steps and weaken the assumptions on $\nu$ gradually.

(1) Assume that $\nu \in (0, \infty)^{\mathcal{S}}$ and that $\kappa$ is irreducible with respect to $\nu$. Consider the case $\Sigma(\kappa, \nu) \neq 1$, then by Lemma 4.5 we have that $\chi(\kappa, \theta(\nu)) > 1$ and thus for any $r \in \mathcal{S}$ the power series $\Gamma_r$ defined as in (4.4) is analytic in $\theta(\nu)$. Applying Lemma 4.4 and using the definition of $\theta(\nu)$ we get that

$$\nu_r \exp(-(\kappa\nu)_r) = \theta_r(\nu) = \Gamma_r(\theta(\nu)) \exp(-(\kappa\Gamma(\theta(\nu)))_r) \qquad (4.14)$$

for any $r \in \mathcal{S}$. In other words, $\Gamma(\theta(\nu))$ is a solution of (4.1). Now, Lemma 4.1 gives that $\nu = \Gamma(\theta(\nu))$ if $\Sigma(\kappa, \nu) \leq 1$ and the claim follows. If $\Sigma(\kappa, \nu) > 1$, then by Lemma 4.1 there exists (a strictly positive) $\nu^* \leq \nu$, $\nu^* \neq \nu$, solving Eq. (4.1) and satisfying $\Sigma(\kappa, \nu^*) < 1$. So, $\Gamma(\theta(\nu)) = \Gamma(\theta(\nu^*)) = \nu^*$ follows by applying the previous case.

Now, consider the case $\Sigma(\kappa, \nu) = 1$, which is equivalent to $\chi(\theta(\nu)) = 1$. Let $\nu^{(n)} \nearrow \nu$, in particular $\Sigma(\kappa, \nu^{(n)}) < 1$ for all $n$. Then by Fatou's lemma and the first case we get that for any $r \in \mathcal{S}$

$$\Gamma_r(\theta(\nu)) = \sum_{k \in \mathbb{N}_0^{\mathcal{S}}} \tau(k) k_r \prod_s \frac{\theta_s(\nu)^{k_s}}{k_s!} \leq \liminf_{n \to \infty} \sum_{k \in \mathbb{N}_0^{\mathcal{S}}} \tau(k) k_r \prod_s \frac{\theta_s(\nu^{(n)})^{k_s}}{k_s!}$$
$$= \liminf_{n \to \infty} \nu_r^{(n)} = \nu_r. \qquad (4.15)$$

Put $\nu^* := \Gamma(\theta(\nu))$ and assume towards a contradiction that $\nu^* \neq \nu$. Then there exists $s \in \mathcal{S}$ such that $\nu_s^* < \nu_s$ and by irreducibility of $\kappa$ there exists at least one $s'$ such that $\kappa(s', s) > 0$, and thus $\kappa(s', s)\nu_s^* < \kappa(s', s)\nu_s$. Hence the Perron–Frobenius theorem implies that $\Sigma(\kappa, \nu^*) < 1$. Therefore the power series is analytic in $\theta(\nu)$, so by Lemma 4.4 and the definition of $\theta(\nu)$ we get that $\nu^* \exp(-\kappa\nu^*) = \nu \exp(-\kappa\nu)$. Applying Lemma 4.1 yields $\nu^* = \nu$ in contradiction to our assumption.

(2) Let $\nu \in (0, \infty)^{\mathcal{S}}$ and an arbitrary $\kappa$. Then we can decompose $\mathcal{S}$ into disjoint sets $S_j$, $j \in J$, such that $\kappa^{(j)} = \kappa|_{S_j \times S_j}$ is irreducible with respect to $\nu^{(j)} = \nu|_{S_j}$ for any $j \in J$. For any $j \in J$ we can apply (1) to get that

$$\Gamma_r^{(j)}(\theta(\nu^{(j)})) := \sum_{k \in \mathbb{N}_0^{S_j}} \tau(k) k_r \prod_{s \in S_j} \frac{(\nu^{(j)} e^{-(\kappa^{(j)}\nu^{(j)})_s})^{k_s})}{k_s!} = \nu_r^{(*, j)}, \quad r \in S_j,$$
$$\qquad (4.16)$$

where $\nu^{(j)}$ solves Eq. (4.1) on $S_j$. Observe that if $k \in \mathbb{N}_0^{\mathcal{S}}$ is such that $\mathrm{supp}(k) \not\subset S_j$ holds for any $j \in J$, then $\tau(k) = 0$. Consequently for fixed $r$ and $j$ such that $r \in S_j$ we get that

$$\Gamma_r(\theta(\nu)) = \Gamma_r^{(j)}(\theta(\nu^{(j)})) \qquad (4.17)$$

where we also used that $(\kappa v)_s = (\kappa^{(j)} v^{(j)})_s$ for $s \in S_j$ holds by construction. Define $v^* = (v_s^*)_{s \in \mathcal{S}}$ by putting $v_s^* := v_s^{(*,j)}$ if $s \in S_j$. Then it is easy to verify that $v^*$ is the smallest solution to (4.1).

(3) Let $v \in [0, \infty)^{\mathcal{S}}$. Observe that for $r \notin \text{supp}(v)$ we have $\Gamma_r(\theta(v)) = 0$. The rest follows by restricting to $\text{supp}(v)$ and applying (2). □

**Remark 4.6** (*Connection with branching processes*) For the reader who is familiar with the results and techniques used in [7], a natural question that arises is about the connection between the power series that we study in Proposition 4.2 and the multi-type branching process that is used in [7] to explore the clusters of the random graph, as we mention in Sect. 2.2. Indeed, both objects carry the same information, as is shown in the following lemma.

**Lemma 4.7** *Fix $v \in (0, \infty)^{\mathcal{S}}$. Let $\mathcal{X}$ be a multi-type branching process, where the individuals are equipped with types from $\mathcal{S}$ and an individual of type $r \in \mathcal{S}$ gives birth to a number of individuals of type $s \in \mathcal{S}$ that is Poisson distributed with parameter $\kappa(r,s)v_s$, independently for each $s \in \mathcal{S}$. Let $\Xi$ be the vector in $\mathbb{N}_0^{\mathcal{S}}$ that counts the total progeny of the process $\mathcal{X}$ according to their types. For $r \in \mathcal{S}$ let $\mathrm{P}_r$ be the probability measure under which $\mathcal{X}$ starts with one single individual of type $r$. Then*

$$\mathrm{P}_r(\Xi = k) = \tau(k) \frac{k_r}{v_r} \prod_{s \in \mathcal{S}} \frac{(v_s e^{-(\kappa v)_s})^{k_s}}{k_s!}, \quad k \in \mathbb{N}_0^{\mathcal{S}}. \qquad (4.18)$$

We omit the proof of the lemma, which comes from combinatorial manipulations and properties of the Poisson distribution of the offspring.

As a consequence, the result of Proposition 4.2 can be reformulated in terms of the branching process $\mathcal{X}$. If $\rho_r^*$ denotes the probability that the process $\mathcal{X}$ goes extinct under $\mathrm{P}_r$, then

$$\rho_r^* = \sum_{k \in \mathbb{N}_0^{\mathcal{S}}} \mathrm{P}_r(\Xi = k) = \frac{1}{v_r} \sum_{k \in \mathbb{N}_0^{\mathcal{S}}} \tau(k) k_r \prod_{s \in \mathcal{S}} \frac{(v_s e^{-(\kappa v)_s})^{k_s}}{k_s!}, \quad r \in \mathcal{S}. \quad (4.19)$$

Furthermore, by substituting $\rho^* := \frac{v^*}{v}$, the statement of Proposition 4.2 is equivalent to the fact that the survival probability $\rho = 1 - \rho^*$ is the maximal non-negative solution to

$$\rho = 1 - e^{-\kappa(\rho v)}$$

as stated in Theorem 2.5.

Recall that, in the single-type case $|\mathcal{S}| = 1$, the right-hand side of (4.18) is identical to the Borel distribution on $\mathbb{N}_0$. Hence, we call it also in the general case where $\mathcal{S}$ is a finite set the *multi-type Borel distribution* and denote it by $\mathrm{Bo}_{\kappa,v}$. It is a probability measure on $\mathbb{N}_0^{\mathcal{S}}$ if and only if $\rho^* = 1$.

### 4.2 Proofs for the asymptotics for the connection probabilities

Using the results from Sect. 4.1 we can now study the connection probabilities for the different cases, namely for microscopic, mesoscopic and macroscopic clusters.

The first result provides the link between the connection probabilities of microscopic clusters and the weight of spanning trees. Its proof is elementary and uses well-known arguments. As a consequence we can prove Lemmas 3.4 and 3.5. Recall the definition of $\tau(k)$ from (3.6) and that $\kappa \in [0, \infty)^{S \times S}$ is the limiting matrix of the sequence $\kappa_N$ as $N \to \infty$. We define $\tau_N(k)$ as in (3.6) with respect to $\kappa_N$ instead of $\kappa$. We will assume that $\kappa_N(r, s) \leq N$ for any $r, s \in S$, which holds in the finite type setting if $N \in \mathbb{N}$ is large enough.

**Lemma 4.8** (Bounds for the connection probability of microscopic clusters) *For any $N \in \mathbb{N}$ and $k \in \mathbb{N}_0^S \setminus \{0\}$,*

$$\left(1 - \frac{\|\kappa_N\|_\infty}{N}\right)^{|k|^2/2} \leq \frac{p_N(k)}{N^{1-|k|} \tau_N(k)} \leq 1. \tag{4.20}$$

**Proof** We start with the upper bound. For $T \in \mathcal{T}(k)$ we denote by $\Omega_T$ the event that the edge set $E(T)$ of $T$ is contained in the edge set of $\mathcal{G}(|k|, \mathbf{x}, \frac{1}{N}\kappa_N)$. Since $\mathcal{G}(|k|, \mathbf{x}, \frac{1}{N}\kappa_N)$ has to contain at least one spanning tree in order to be connected, we have

$$p_N(k) = \mathbb{P}\Big(\bigcup_{T \in \mathcal{T}(k)} \Omega_T\Big) \leq \sum_{T \in \mathcal{T}(k)} \mathbb{P}(\Omega_T)$$

$$= \sum_{T \in \mathcal{T}(k)} \prod_{\{i,j\} \in E(T)} \frac{\kappa_N(x_i, x_j)}{N} \leq \left(\frac{1}{N}\right)^{|k|-1} \tau_N(k),$$

where we used the fact that each $T \in \mathcal{T}(k)$ has exactly $|k| - 1$ edges.

Now we continue with the lower bound. For $T \in \mathcal{T}(k)$ we denote by $\widetilde{\Omega}_T$ the event that the edge set of $\mathcal{G}(|k|, \mathbf{x}, \frac{1}{N}\kappa_N)$ is equal to $E(T)$. Note that the events $\widetilde{\Omega}_T$, $T \in \mathcal{T}(k)$, are disjoint and therefore

$$p_N(k) \geq \sum_{T \in \mathcal{T}(k)} \mathbb{P}(\widetilde{\Omega}_T) = \sum_{T \in \mathcal{T}(k)} \Big(\prod_{\{i,j\} \in E(T)} \frac{\kappa_N(x_i, x_j)}{N}\Big)$$

$$\times \prod_{\{i,j\} \notin E(T)} \left(1 - \frac{\kappa_N(x_i, x_j)}{N}\right)$$

$$\geq \left(1 - \frac{\|\kappa_N\|_\infty}{N}\right)^{\binom{|k|}{2} - (|k|-1)} \left(\frac{1}{N}\right)^{|k|-1} \tau_N(k)$$

$$\geq \left(1 - \frac{\|\kappa_N\|_\infty}{N}\right)^{|k|^2/2} \left(\frac{1}{N}\right)^{|k|-1} \tau_N(k).$$

□

**Proof of Lemma 3.4** Equation (3.17) is a Corollary of the estimate (4.20) and the fact that the left-hand side of (4.20) converges to 1, as $N \to \infty$, since $|k|$ does not depend on $N$ and $\|\kappa_N\|_\infty$ is bounded in $N$. □

**Proof of Lemma 3.5** For fixed $k \in \mathbb{N}_0^S$ denote $S_k := \mathrm{supp}(k)$. For $p_N(k)$ and $r \in S_k$ we use the upper bound in (4.20) from Lemma 4.8 and also the formula for $\tau_N(k) k_r$ given in (4.6) to obtain

$$p_N(k) \leq N^{1-|k|} \tau_N(k) = N^{1-|k|} \frac{1}{k_r} \Big( \prod_{s \in S_k} (\kappa_N k)_s^{k_s-1} \Big) \Delta_{N,r}(k), \qquad (4.21)$$

where $\Delta_{N,r}(k)$ is defined as in (4.7) but with respect to the kernel $\kappa_N$. Now, observe that

$$\Delta_{N,r}(k) = \sum_{A \in \vec{\mathcal{T}}_r(S_k)} \prod_{\{s,s'\} \in E(A)} \kappa_N(s,s') k_{s'} \leq \|\kappa_N\|_\infty^{|S_k|-1} |\mathcal{T}_r(S_k)| \prod_{s \in S_k \setminus \{r\}} k_s, \qquad (4.22)$$

and by Cayley's formula we have that $|\mathcal{T}_r(S_k)| = |S_k|^{|S_k|-1}$. Hence, combining Eq. (4.21) with the inequality (4.22) gives (3.18). □

Now we turn to the proof of Theorem 3.6, which comes as a consequence of the following Lemmas 4.10–4.13. The intuitive idea of the proof is to embed $\mathcal{G}(|k^{(N)}|, \mathbf{x}, \frac{1}{N}\kappa_N)$ in a larger random graph with vertex set given by some $m^{(N)} \in \mathbb{N}_0^S$ such that the component $k^{(N)}$ appears with high probability as the typical giant component in this graph. We make heavy use of the expansion in Lemma 4.9, which is a straightforward multi-type generalization of equation (4) from [24].

**Lemma 4.9** *Fix $m \in \mathbb{N}_0^S$ and $r \in S$ such that $m_r \geq 1$. Then the following formula holds*

$$1 = \sum_{h \in \mathbb{N}_0^S : \, \mathbf{e}_r \leq h \leq m} \Big[ \prod_{s \in S} \binom{m_s - \delta_{r,s}}{h_s - \delta_{r,s}} \Big] p_N(h) \prod_{s, \tilde{s} \in S} \Big( 1 - \frac{\kappa_N(s, \tilde{s})}{N} \Big)^{h_s(m_{\tilde{s}} - h_{\tilde{s}})}. \qquad (4.23)$$

**Proof** Consider the graph $\mathcal{G}(|m|, x, \frac{1}{N}\kappa_N)$ where $x$ is a fixed vector compatible with $m$. Fix a vertex $i \in \{1, \ldots, |m|\}$ of type $r$. For fixed $h \in \mathbb{N}_0^S$ with $\mathbf{e}_r \leq h \leq m$ denote by $\Omega(h)$ the event that the connected component containing $i$ is given by some set $C \subset \{1, \ldots, |m|\}$ that contains exactly $h_s$ vertices of type $s$ for each $s \in S$. We claim that the summand on the right-hand side of (4.23) is the probability of $\Omega(h)$. Indeed, there are $\prod_{s \in S} \binom{m_s - \delta_{r,s}}{h_s - \delta_{r,s}}$ possibilities to choose a set $C \setminus \{i\}$ from $\{1, \ldots, |m|\} \setminus \{i\}$. The probability that $C$ is the connected component containing vertex $i$ is given as the product of $p_N(h)$, i.e., the probability that it is connected, and the probability that no edge exists between $C$ and its complement $\{1, \ldots, |m|\} \setminus C$, which is easily seen to be equal to $\prod_{s,\tilde{s} \in S}(1 - \frac{\kappa_N(s,\tilde{s})}{N})^{h_s(m_{\tilde{s}} - h_{\tilde{s}})}$. Equation (4.23) follows by the observation that the events $\Omega(h)$, $\mathbf{e}_r \leq h \leq m$, form a decomposition of the underlying probability space. □

The idea is now to pick $m = m^{(N)}$ in such a way that the summand for $h = k^{(N)}$ is maximal, such that, on the exponential scale, the right-hand side of (4.23) can be replaced by just this summand. It will turn out that the correct choice is $m^{(N)} \sim N\nu$ with $\nu_r = y_r/(1-e^{-(\kappa y)_r})$ for $r \in \mathcal{S}$. Intuitively, this choice of $\nu$ comes from inverting equation (2.11): indeed one can see that $y = \rho\nu$ where $\rho$ solves

$$\rho = 1 - e^{-T_{\kappa,\nu}\rho}.$$

Notice that the assumption $y \ll \kappa y$ is crucial for $\nu$ to be well-defined. While in the upper bound of Lemma 4.10 this is not important (in that case we see from (4.25) that the upper bound of $p_N(k^{(N)})$ converges to 0), we need it in order to get a non-trivial lower bound in Lemma 4.12.

**Lemma 4.10** (Upper bound in (3.20)) *Let $k \in \mathbb{N}_0^{\mathcal{S}}$ with $k_r \geq 1$ for some $r \in \mathcal{S}$. Let $m \in \mathbb{N}_0^{\mathcal{S}}$ such that $m \leq k$. Then*

$$p_N(k) \leq \left[\prod_{s\in\mathcal{S}} \binom{m_s - \delta_{r,s}}{k_s - \delta_{r,s}}\right]^{-1} \prod_{s,\tilde{s}\in\mathcal{S}} \left(1 - \frac{\kappa_N(s,\tilde{s})}{N}\right)^{-k_s(m_{\tilde{s}} - k_{\tilde{s}})}. \quad (4.24)$$

*Fix $y \in [0,1]^{\mathcal{S}}\setminus\{0\}$. Let $\{k^{(N)}\}_{N\in\mathbb{N}}$ be a sequence in $\mathbb{N}_0^{\mathcal{S}}$ such that $\lim_{N\to\infty} \frac{k_r^{(N)}}{N} = y_r$ for all $r \in \mathcal{S}$. Then*

$$\limsup_{N\to\infty} \frac{1}{N} \log p_N(k^{(N)}) \leq \sum_{r\in\mathcal{S}} y_r \log\left(1 - e^{-(\kappa y)_r}\right), \quad (4.25)$$

*where the right-hand side takes the value $-\infty$ when $y \ll \kappa y$.*

**Proof** The inequality (4.24) is a direct consequence of (4.23), where from the right-hand side we pick only the summand for $h = k$, which is present since $m \geq k$.

Let us focus now on (4.25). Let $m^{(N)} \in \mathbb{N}_0^{\mathcal{S}}$ be such that $m_r^{(N)} := \lfloor k_r^{(N)}(1 - e^{-\frac{1}{N}(\kappa k^{(N)})_r})^{-1}\rfloor$ for all $r \in \mathcal{S}$, therefore $m^{(N)} \geq k^{(N)}$. Fix $r \in \mathcal{S}$ such that $k_r^{(N)} \geq 1$ and consider (4.24).

For brevity we will write $k$ and $m$ instead of $k^{(N)}$ and $m^{(N)}$. With the help of Stirling's formula $n! = (n/e)^n e^{o(n)}$ and the exponential limit theorem $\lim_{n\to\infty}(1 + \frac{c}{n})^n = e^c$ and some elementary calculation, we get that, as $N \to \infty$,

$$\left[\prod_{s\in\mathcal{S}} \binom{m_s - \delta_{r,s}}{k_s - \delta_{r,s}}\right] \prod_{s,\tilde{s}\in\mathcal{S}} \left(1 - \frac{\kappa_N(s,\tilde{s})}{N}\right)^{k_s(m_{\tilde{s}} - k_{\tilde{s}})}$$
$$= \exp\left(-\langle k, \log\frac{k}{m}\rangle - \langle m-k, \log\frac{m-k}{m}\rangle\right) \times \exp\left(-\langle m-k, \kappa k\rangle\right) e^{o(N)}$$
$$= \exp\left(-\langle k, \log(1 - e^{-\frac{1}{N}\kappa k})\rangle\right) e^{o(N)},$$

where in the last step we used that $k/m = 1 - e^{-\frac{1}{N}\kappa k}$ (by the definition of $m$). Now, we insert the asymptotics in (4.24) and see that

$$p_N(k^{(N)}) \leq \exp\left(N \sum_{r \in \mathcal{S}} \frac{k_r^{(N)}}{N} \log\left(1 - e^{-\frac{1}{N}(\kappa k^{(N)})_r}\right)\right) e^{o(N)}.$$

Consequently, the claim in (4.25) follows, including the convention for $y \not\ll \kappa y$. □

In order to prove the lower bound we need the following auxiliary lemma.

**Lemma 4.11** *For $h', h \in \mathbb{N}_0^{\mathcal{S}}$ with $h' \leq h$ we have, for any $N \in \mathbb{N}$,*

$$p_N(h') \leq p_N(h) \prod_{r \in \mathcal{S}} \left(1 - e^{-\frac{1}{N}(\kappa_N h')_r}\right)^{-(h_r - h'_r)}. \tag{4.26}$$

*Proof* Let $I := [|h|]$ and let $x = (x_i)_{i \in I}$ be compatible with $h$, i.e., $\sum_{i \in I} \delta_{x_i} = h$. Consider the graph $G := \mathcal{G}(|h|, x, \frac{1}{N}\kappa_N)$ on the vertex set $I$. There exists $I' \subset I$ such that $x' = (x_i)_{i \in I'}$ is compatible with $h'$. The subgraph $G'$ of $G$ that is induced by $I'$ can be identified with $\mathcal{G}(|h'|, x', \frac{1}{N}\kappa'_N)$, where $\kappa'$ denotes the restriction of $\kappa$ to $I' \times I'$. If $G'$ is connected and for any $i \in I \setminus I'$ there is an edge $\{i, j\} \in E(G)$ with $j \in I'$, then also $G$ is connected. Therefore, using that $1 - x \leq e^{-x}$ for any $x \in \mathbb{R}$, we have

$$p_N(h) \geq p_N(h') \prod_{r \in \mathcal{S}} \left(1 - \prod_{s \in \mathcal{S}} \left(1 - \frac{\kappa_N(r, s)}{N}\right)^{h'_s}\right)^{h_r - h'_r}$$
$$\geq p_N(h') \prod_{r \in \mathcal{S}} \left(1 - e^{-\frac{1}{N}(\kappa_N h')_r}\right)^{h_r - h'_r}.$$

□

In the following lemmas we give the lower bound for the connection probabilities in the macroscopic setting. It is sufficient to restrict to the case $y \ll \kappa y$, since otherwise the limit is already ensured to be $-\infty$ by Lemma 4.10.

**Lemma 4.12** (Lower bound in (3.20) for irreducible $\kappa$) *Fix $y \in (0, 1]^{\mathcal{S}}$ satisfying $y \ll \kappa y$ and assume that $\kappa$ is irreducible with respect to $y$. Let $\{k^{(N)}\}_{N \in \mathbb{N}}$ be a sequence in $\mathbb{N}_0^{\mathcal{S}}$ such that $\lim_{N \to \infty} \frac{k_r^{(N)}}{N} = y_r$ for all $r \in \mathcal{S}$ and assume that $\tau(k^{(N)}) > 0$ holds for all $N \in \mathbb{N}$. Then*

$$\liminf_{N \to \infty} \frac{1}{N} \log p_N(k^{(N)}) \geq \sum_{r \in \mathcal{S}} y_r \log\left(1 - e^{-\kappa y_r}\right). \tag{4.27}$$

*Proof* For $\delta \geq 0$ satisfying $2\delta < \inf\{y_s : s \in \mathcal{S}\}$ define $y^{(\delta)} := y - \delta$ as well as

$$v_s^{(\delta)} := \frac{y_s^{(\delta)}}{1 - e^{-(\kappa y^{(\delta)})_s}}, \quad \text{for } s \in \mathcal{S}. \tag{4.28}$$

For purely technical reasons that will become apparent later, we cannot work directly with $\nu := \nu^{(0)}$. Note that for all $\delta \geq 0$ our construction ensures that with respect to $\nu^{(\delta)}$ the characteristic equation (4.1) has a non-trivial solution, which is equivalent to $\Sigma(\kappa, \nu^{(\delta)}) > 1$ by Lemma 4.1. Let $m^{(N,\delta)} \in \mathbb{N}_0^{\mathcal{S}}$ be such that $m_r^{(N,\delta)} = \lfloor (k^{(N)} - N\delta)(1 - e^{-\frac{1}{N}(\kappa(k-N\delta))_r})^{-1} \rfloor$. In particular, $\lim_{N \to \infty} \frac{m^{(N,\delta)}}{N} = \nu^{(\delta)}$. Clearly, we have that $\nu^{(\delta)} \to \nu$ as $\delta \to 0$. Note that $y \leq \nu$. Moreover, there exist $\delta^* > 0$ and $N_0$ such that $m^{(N,\delta)} \geq k^{(N)}$ for all $\delta \leq \delta^*$ and $N \geq N_0$, which we will assume from now on. For brevity we will write $k$ and $m^{(\delta)}$ instead of $k^{(N)}$ and $m^{(N,\delta)}$.

Fix $r \in \mathcal{S}$ with $y_r > 0$ and note that by the assumption $y \ll \kappa y$ this implies $(\kappa y)_r > 0$ and thus $y_r < \nu_r$. From here on, $h$ will always denote an element in the set $\{h \in \mathbb{N}_0^{\mathcal{S}} : \mathbf{e}_r \leq h \leq m^{(\delta)}\}$. For such $h$ abbreviate

$$a_N^{(\delta)}(h) := \left[\prod_{s \in \mathcal{S}} \binom{m_s^{(\delta)} - \delta_{r,s}}{h_s - \delta_{r,s}}\right] p_N(h) \prod_{s,\tilde{s} \in \mathcal{S}} \left(1 - \frac{\kappa_N(s, \tilde{s})}{N}\right)^{h_s(m_{\tilde{s}}^{(\delta)} - h_{\tilde{s}})}. \quad (4.29)$$

Recall from Lemma 4.9 that the sum of $a_N^{(\delta)}(h)$ over the mentioned $h$'s is equal to one. In the following, we will split this sum into the three parts where $|h| \leq R$, and $R < |h| \leq \varepsilon N$ and $|h| > \varepsilon N$ for some large $R \in \mathbb{N}$ and some small $\varepsilon > 0$. We will show that the first part converges towards something in $(0, 1)$, the second part vanishes, and we will identify the exponential rate of the other one explicitly via some Laplace approximation. Explicitly, we will prove the following three claims:

*Claim 1:* $\quad \limsup_{\delta \downarrow 0} \limsup_{R \to \infty} \limsup_{N \to \infty} \sum_{h: |h| \leq R} a_N^{(\delta)}(h) \leq 1 - \frac{y_r}{\nu_r} \in [0, 1).$ (4.30)

*Claim 2:* If $\varepsilon$ is small enough, for some sufficiently small $\delta_0 > 0$ and some sufficiently large $N_0 \in \mathbb{N}$,

$$\limsup_{R \to \infty} \sup_{N \geq N_0} \sup_{\delta \in (0, \delta_0]} \sum_{h: R < |h| \leq \varepsilon N} a_N^{(\delta)}(h) = 0. \quad (4.31)$$

*Claim 3:* For any $\delta > 0$,

$$\sum_{h: |h| > \varepsilon N} a_N^{(\delta)}(h) \leq e^{o(N)} e^{NC(\delta)} e^{-N\tilde{f}(y,\nu)} p_N(k) + \varepsilon_N, \quad (4.32)$$

for some sequence $(\varepsilon_N)_N = (\varepsilon_N(\delta))_N$ that converges to 0, where $\tilde{f}(y, \nu) = \langle y, \log(1 - e^{-\kappa y})\rangle$, and $C(\delta)$ is a constant only depending on $\delta$ and vanishing as $\delta \downarrow 0$.

Let us first explain how the assertion of the lemma follows from these three claims. We start by using Lemma 4.9 and Claim 3 to obtain, in the limit as $N \to \infty$,

$$1 - \sum_{h: |h| \leq R} a_N^{(\delta)}(h) - \sum_{h: R < |h| < \varepsilon N} a_N^{(\delta)}(h) \leq e^{o(N)} e^{NC(\delta)} e^{-N\tilde{f}(y,\nu)} p_N(k) + \varepsilon_N.$$

By (4.30) and (4.31), the left hand side is not smaller than $y_r/v_r$, when taking the limits as $N \to \infty$, followed by $R \to \infty$ and $\delta \downarrow 0$, if $\varepsilon$ is small enough. In particular, the left-hand side is bounded away from zero when taking these limits. Hence, the exponential rate of the right-hand side as $N \to \infty$ is nonnegative. This implies that $\liminf_{N \to \infty} \frac{1}{N} \log p_n(k) \geq \widetilde{f}(y, \nu)$, which is the assertion of the lemma. It remains to prove the three claims.

*Proof of Claim 3:* As we explained in the proof of Lemma 4.10, the exponential rate of the two products in the definition of $a_N^{(\delta)}(h)$ can easily be identified with the help of Stirling's formula and the exponential limit theorem $\lim_{n \to \infty}(1 + \frac{c}{n})^n = e^c$ and elementary calculations. Indeed, for any sequence $h^{(N)}$ such that $\lim_{N \to \infty} \frac{h^{(N)}}{N} = x$ exists and $x \in [0, 1]^{\mathcal{S}} \setminus \{0\}$, and for any $\delta > 0$,

$$a_N^{(\delta)}(h^{(N)}) = e^{o(N)} e^{-Nf(x, \nu^{(\delta)})} p_N(h^{(N)}), \qquad N \to \infty, \tag{4.33}$$

where we introduce $f(x, \tilde{x}) = \langle x, \log \frac{x}{\tilde{x}} \rangle + \langle \tilde{x} - x, \log \frac{\tilde{x}-x}{\tilde{x}} \rangle + \langle \tilde{x} - x, \kappa x \rangle$ for $x, \tilde{x} \in [0, 1]^{\mathcal{S}} \setminus \{0\}$ and $x \leq \tilde{x}$. We can write

$$\begin{aligned} f(x, \tilde{x}) &= \sum_{s \in \mathcal{S}} \tilde{x}_s \left[ \frac{x_s}{\tilde{x}_s} \log \frac{\frac{x_s}{\tilde{x}_s}}{1 - e^{-(\kappa x)_s}} + \left(1 - \frac{x_s}{\tilde{x}_s}\right) \log \frac{1 - \frac{x_s}{\tilde{x}_s}}{e^{-(\kappa x)_s}} \right] \\ &\quad + \sum_{s \in \mathcal{S}} x_s \log(1 - e^{-(\kappa x)_s}) \\ &= \langle \tilde{x}, H(x; \tilde{x}) \rangle + \langle x, \log(1 - e^{-\kappa x}) \rangle, \end{aligned} \tag{4.34}$$

where we write $H = (H_s)_{s \in \mathcal{S}}$ and for $s \in \mathcal{S}$ we write $H_s(x; \tilde{x})$ for the entropy of the Bernoulli distribution with parameter $x_s/\tilde{x}_s$ with respect to the one with parameter $1 - e^{-(\kappa x)_s}$.

The second term in the second line of (4.34) will be handled jointly with the exponential rate of $p_N(h^{(N)})$, so let us discuss here the minimum of the first. Since every component of $H$ is an entropy between probability measures, we have that $H(x, \nu^{(\delta)}) \geq 0$ pointwise for all $x \leq \nu^{(\delta)}$ with equality if and only if $\frac{x}{\nu^{(\delta)}} = 1 - e^{-\kappa x}$. As $\kappa$ is irreducible with respect to $\nu^{(\delta)}$ the latter condition is true if and only if $x = y^{(\delta)}$ or $x = 0$, see Lemma 4.1. Hence, we have seen that, for $\varepsilon$ sufficiently small, $\min_{x: |x| \geq \varepsilon} \langle \nu^{(\delta)}, H(x; \nu^{(\delta)}) \rangle = \langle \nu^{(\delta)}, H(y^{(\delta)}; \nu^{(\delta)}) \rangle$.

Let us now prove (4.32). The work we still have to do is to combine (4.33) with estimates for the $p_N$-term from Lemmas 4.10 and 4.11, and we have to distinguish the cases that $h$ is left of $k$ but close to $k$ (then Lemma 4.11 applies), or bounded away from $k$ (then Lemma 4.10 suffices). Let $h^{(N)}$ be such that $|h^{(N)}| > \varepsilon N$ and assume that $x := \lim_{N \to \infty} \frac{h^{(N)}}{N}$ exists. Note that $x \in [0, 1] \setminus \{0\}$. We first examine the case where $h^{(N)} \in [k - 2N\delta, k]$ for large $N$ and hence $|x - y| \leq 2\delta$. With the help of (4.33) and using the estimate from Lemma 4.11 we get that

$$a_N^{(\delta)}(h^{(N)}) \leq e^{o(N)} e^{N(f(y,\nu) - f(x,\nu^{(\delta)}) - \langle y-x, \log(1-e^{-\kappa x}) \rangle)} e^{-Nf(y,\nu)} p_N(k)$$

$$\leq e^{o(N)} e^{NC(\delta)} e^{-Nf(y,\nu)} p_N(k) \qquad (4.35)$$

where

$$C(\delta) = \sup_{x:\,|x-y|\leq 2\delta} |f(y,\nu) - f(x,\nu^{(\delta)}) - \langle y - x, \log(1 - e^{-\kappa x})\rangle| \qquad (4.36)$$

and it can be easily verified that $C(\delta) \to 0$ as $\delta \to 0$. Now, let us examine the case $h^{(N)} \notin [k - 2N\delta, k]$ for large $N$, i.e., $|x - y^{(\delta)}| \geq \delta$. We start with (4.33) and use the estimate (4.25) from Lemma 4.10 to get

$$a_N^{(\delta)}(h^{(N)}) \leq e^{o(N)} e^{-N\langle \nu^{(\delta)}, H(x,\nu^{(\delta)})\rangle} \leq e^{o(N)} e^{-N\widetilde{C}(\delta)} \qquad (4.37)$$

where the number

$$\widetilde{C}(\delta) := \inf_{x:\,|x|>\varepsilon,\,|x-y^{(\delta)}|\geq \delta} \langle \nu^{(\delta)}, H(x,\nu^{(\delta)})\rangle \qquad (4.38)$$

is strictly positive since the function $x \mapsto H(x, \nu^{(\delta)})$ is continuous and its only zeros are at $x = 0$ and $x = y^{(\delta)}$. Now, it is not hard to see that $\#\{h \in \mathbb{N}_0^\mathcal{S}: h \leq m^{(\delta)}\} \leq (\nu N)^{|\mathcal{S}|} = e^{o(N)}$, hence a simple Laplace approximation argument implies that (4.32) holds with $\varepsilon_N := e^{o(N)} e^{-N\widetilde{C}(\delta)}$, which vanishes as $N \to \infty$ since $\widetilde{C}(\delta) > 0$.

*Proof of Claim 1*: Applying Lemma 4.8 to $p_N(h)$ we get, for any $N, R \in \mathbb{N}$ and $\delta > 0$,

$$\begin{aligned}
&\sum_{h:\,|h|\leq R} a_N^{(\delta)}(h) \\
&\leq \frac{N}{m_r^{(\delta)}} \sum_{h:\,|h|\leq R} \tau_N(h) h_r \prod_{s\in\mathcal{S}} \frac{\frac{m_s^{(\delta)}}{N} \frac{m_s^{(\delta)}-1}{N} \cdots \frac{m_s^{(\delta)}-h_s+1}{N} \left(\prod_{\tilde{s}\in\mathcal{S}}(1 - \frac{\kappa_N(s,\tilde{s})}{N})^{m_{\tilde{s}}^{(\delta)}-h_{\tilde{s}}}\right)^{h_s}}{h_s!} \\
&\leq \frac{N}{m_r^{(\delta)}} \sum_{h:\,|h|\leq R} \tau_N(h) h_r \prod_{s\in\mathcal{S}} \left[\frac{1}{h_s!} \left(\frac{m_s^{(\delta)}}{N} e^{-\kappa_N(\frac{m^{(\delta)}-h}{N})_s}\right)^{h_s}\right],
\end{aligned} \qquad (4.39)$$

where we used that $1 - x \leq e^{-x}$ for any $x \in \mathbb{R}$ in the second line. Hence,

$$\limsup_{N\to\infty} \sum_{h:\,|h|\leq R} a_N^{(\delta)}(h) \leq \frac{1}{\nu_r^{(\delta)}} \sum_{h:\,|h|\leq R} \tau(h) h_r \prod_{s\in\mathcal{S}} \left[\frac{1}{h_s!} \left(\nu_s^{(\delta)} e^{-(\kappa\nu^{(\delta)})_s}\right)^{h_s}\right] =: \Lambda_R^{(\delta)} \qquad (4.40)$$

Observe that by Proposition 4.2 and by the definition of $\nu^{(\delta)}$ given in (4.28), we have that $\lim_{R\to\infty} \Lambda_R^{(\delta)} = \frac{\nu_r^{(\delta)} - y_r^{(\delta)}}{\nu_r^{(\delta)}} \in [0, 1)$ and clearly $\lim_{\delta\to 0} \lim_{R\to\infty} \Lambda_R^{(\delta)} = \frac{\nu_r - y_r}{\nu_r} \in (0, 1)$.

*Proof of Claim 2:* For the sum on $h$ satisfying $R < |h| \leq \varepsilon N$, we can use basically the same estimates as in (4.39), but we modify the second line by now estimating

$$\prod_{s,\tilde{s}\in\mathcal{S}} \left(1 - \tfrac{\kappa_N(s,\tilde{s})}{N}\right)^{(m_{\tilde{s}}^{(\delta)} - h_{\tilde{s}})h_s} \leq \left(\prod_{s\in\mathcal{S}} e^{-\frac{1}{N}(\kappa_N m^{(\delta)})_s h_s}\right) e^{\frac{1}{N}\|\kappa_N\|_\infty |h|^2}, \quad (4.41)$$

since we do not have that $m^{(\delta)} \geq h$ on this sum. This gives

$$\sum_{h:\, R<|h|\leq \varepsilon N} a_N^{(\delta)}(h) \leq \frac{N}{m_r^{(\delta)}} \sum_{n=R+1}^{\varepsilon N} e^{n\varepsilon \|\kappa_N\|_\infty}$$

$$\sum_{h:\, |h|=n} \tau_N(h) h_r \prod_{s\in\mathcal{S}} \left[\frac{1}{h_s!} \left(\frac{m_s^{(\delta)}}{N} e^{-(\frac{1}{N}\kappa_N m^{(\delta)})_s}\right)^{h_s}\right]. \quad (4.42)$$

This is the main part of a power series with coefficients $\tau$, which we studied in Lemma 4.5. Note that the term in round brackets converges to $\nu_s e^{-(\kappa\nu)_s}$ as $N \to \infty$, followed by $\delta \downarrow 0$. Furthermore, this vector $\theta(\kappa, \nu) = \nu e^{-\kappa\nu}$ satisfies $\chi(\kappa, \theta(\kappa, \nu)) > 1$, since $\Sigma(\kappa, \nu) > 1$ by construction (see Lemma 4.5(3)). It will turn out that this is sufficient to estimate the sum on $n$ against some convergent geometric series.

We pick a small threshold $\varepsilon' > 0$. We may pick $\delta_0$ so small and then $N_0 \in \mathbb{N}$ so large that

$$\frac{m_s^{(\delta)}}{N} e^{-\left(\frac{1}{N}\kappa_N m^{(\delta)}\right)_s} \leq \nu_s e^{-(\kappa\nu)_s} e^{\varepsilon'}, \qquad \delta \in (0, \delta_0], N \geq N_0, s \in \mathcal{S}.$$

Additionally, we can also assume that $N/m_r^{(\delta)} \leq 1 + 1/\nu_r$ for these $\delta$ and $N$. By $\tau$ we denote the function defined as in (3.6) for the matrix $\kappa$, then, since $\lim_{N\to\infty} \kappa_N = \kappa$, we may also pick $N_0$ so large that additionally

$$\tau_N(h) \leq \tau(h) e^{\varepsilon'|h|}, \qquad h \in \mathbb{N}_0^{\mathcal{S}}, N \geq N_0.$$

We also may and will assume that $\|\kappa_N\|_\infty \leq 1 + \|\kappa\|_\infty$ for $N \geq N_0$. Therefore, we can now estimate, for $\delta \in (0, \delta_0]$ and $N \geq N_0$, and any $R \in \mathbb{N}$,

$$\sum_{h:\, R<|h|\leq\varepsilon N} a_N^{(\delta)}(h) \leq \left(1 + \tfrac{1}{\nu_r}\right) \sum_{n=R+1}^{\infty} e^{n\varepsilon(1+\|\kappa\|_\infty)} e^{2\varepsilon' n}$$

$$\sum_{h:\, |h|=n} \tau(h) h_r \prod_{s\in\mathcal{S}} \left[\frac{1}{h_s!}(\nu_s e^{-(\kappa\nu)_s})^{h_s}\right].$$

Now we can apply (4.12) in Lemma 4.5 to estimate the sum on $h$ against $e^{o(n)}e^{-n[\chi(\kappa,\theta(\kappa,\nu))-1]}$. We recall that $\chi(\kappa, \theta(\kappa, \nu)) - 1 > 0$ and pick now $\varepsilon$ and $\varepsilon'$ so small that the sum on $n$ can be estimated against a convergent geometric series. Hence, the entire sum vanishes as $R \to \infty$, which ends the proof. □

In the following we want to verify the statement of Lemma 4.12 without assuming that $\kappa$ is irreducible with respect to $y$ but still under the assumption that $\tau(k^{(N)}) > 0$. What we have in mind are components $k^{(N)}$ where for a certain set of types $\tilde{S} \subset S$ only $o(N)$ vertices of those types are available, i.e., $k_s^{(N)} \in o(N)$ for all $s \in \tilde{S}$. However, these types might be necessary in order to connect the vertices with types in $S \setminus \tilde{S}$. For that case we prove in the following that on the exponential scale the asymptotics of the connection probability $p_N(k^{(N)})$ are the same as in the previous lemma.

**Lemma 4.13** (Lower bound in (3.20) under generalized assumptions) *Fix $y \in [0, 1]^S \setminus \{0\}$ satisfying $y \ll \kappa y$. Let $\{k^{(N)}\}_{N \in \mathbb{N}}$ be a sequence in $\mathbb{N}_0^S$ such that $\lim_{N \to \infty} \frac{k_r^{(N)}}{N} = y_r$ for all $r \in S$ and assume that $\{k_s^{(N)}\}_N$ is bounded for all $s \notin \mathrm{supp}(y)$. Assume further that $\tau(k^{(N)}) > 0$ holds for all but finitely many $N \in \mathbb{N}$. Then*

$$\liminf_{N \to \infty} \frac{1}{N} \log p_N(k^{(N)}) \geq \sum_{r \in S} y_r \log\left(1 - e^{-\kappa y_r}\right). \tag{4.43}$$

*Proof* Throughout the proof we will assume without loss of generality that $S = \bigcap_N \mathrm{supp}(k^{(N)})$. If $y > 0$ on $S$ and $\kappa$ is irreducible with respect to $y$, then we are in the setting of Lemma 4.12, which yields the claim. In the other case, it always holds that $\tilde{S} := S \setminus \mathrm{supp}(y) \neq \emptyset$. We decompose $\mathrm{supp}(y)$ into disjoint sets $S_j$, $j \in J$, satisfying that the restriction of $\kappa$ to $S_j \times S_j$ is irreducible with respect to $y$ restricted to $S_j$.

We will start with the special case where $k_s^{(N)} = 1$ for all $s \in \tilde{S}$ and $N \in \mathbb{N}$. It is easily seen that the assumption $\tau(k^{(N)}) > 0$ for all but finitely many $N \in \mathbb{N}$ implies that there exists a tree $T$ on the type set $S$ such that $\prod_{\{r,s\} \in E(T)} \kappa(r, s) > 0$. Define

$$\tilde{E} := \{\{r, s\} \in E(T) \colon r \in \tilde{S} \text{ or } s \in \tilde{S}\}. \tag{4.44}$$

For every $r \in S$ we now fix a vertex $i_r \in [|k^{(N)}|]$ that is of type $r$. Note that $\{i_r \colon r \in \tilde{S}\}$ already contains all the vertices with types in $\tilde{S}$. We will abbreviate $\mathcal{G}_N = \mathcal{G}(|k^{(N)}|, x, \frac{1}{N}\kappa_N)$. Note that the event $\{\{i_r, i_s\} \in E(\mathcal{G}_N) \colon \{r, s\} \in \tilde{E}\}$ is independent of the existence of edges $\{i, i'\}$, if the types of $i$ and $i'$ are in $S_j$ for any $j \in J$. Therefore,

$$\begin{aligned} p_N(k^{(N)}) &\geq \Big[\prod_{j \in J} p_N(k^{(N)} \mathbb{1}_{S_j})\Big] \times \prod_{\{r,s\} \in \tilde{E}} \mathbb{P}(\{i_r, i_s\} \in E(\mathcal{G}_N)) \\ &\geq \Big[\prod_{j \in J} p_N(k^{(N)} \mathbb{1}_{S_j})\Big] \times \prod_{\{r,s\} \in \tilde{E}} \frac{1}{N} \kappa_N(r, s). \end{aligned}$$

Clearly the second product is $\geq e^{o(N)}$, since $\kappa_N$ is bounded away from zero, uniformly in all large $N$. Applying Lemma 4.12 to $k^{(N)}\mathbb{1}_{S_j}$ for every $j \in J$, we get that

$$\frac{1}{N} \log p_N(k^{(N)}) \geq o(1) + \sum_{j \in J} \sum_{r \in S_j} y_r \log(1 - e^{-(\kappa(y\mathbb{1}_{S_j}))_r})$$
$$= o(1) + \sum_{r \in \mathcal{S}} y_r \log(1 - e^{-(\kappa y)_r}), \quad (4.45)$$

where in the last equation we used that with fixed $j \in J$, for any $r \in S_j$ and $s \in \mathcal{S} \setminus S_j$ we have either $\kappa(r,s) = 0$ or $y_s = 0$. Thus, for $r \in S_j$ we get

$$(\kappa(y\mathbb{1}_{S_j}))_r = \sum_{s \in S_j} \kappa(r,s)y_s = \sum_{s \in \mathcal{S}} \kappa(r,s)y_s = (\kappa y)_r. \quad (4.46)$$

This implies that Eq. (4.43) holds in the special case.

Now, consider the general case, where $\{k_s^{(N)}\}_N$ is bounded in $N$ for all $s \in \tilde{S}$. Then by Lemma 4.11

$$p_N(k^{(N)}) \geq p_N(k^{(N)}\mathbb{1}_{\mathcal{S}\setminus\tilde{S}} + \mathbb{1}_{\tilde{S}}) \times \prod_{r \in \tilde{S}}(1 - e^{-(\kappa y)_r})^{k_r^{(N)}-1}.$$

Observe that the last product is $e^{o(N)}$ by assumption, hence the claim is implied by the treatment of the first case. □

## 5 Proof of Theorem 1.1 for a general type set

The aim of this section is to finally prove our main result, Theorem 1.1, for the general case, where the type space $\mathcal{S}$ is some compact metric space. The main idea is to approximate the general graph model introduced in Sect. 1.1 by a discretized model with a finite type space and discretized kernels. To derive the LDP we will use the LDP of Theorem 3.1 for finite type spaces, together with the Dawson–Gärtner theorem, which we will slightly modify for our purposes.

In Sect. 5.1 we introduce the approximation scheme and derive lower and upper approximations for the distribution of the empirical measures of the general graph via a comparison with certain discretized graph models. In order to lift statements about the distribution of the discretized model to the general case we introduce in Sect. 5.2 a projective system. There, we also identify the projective limit spaces with the state spaces of our empirical measures, i.e., with $\mathcal{L} \times \mathcal{A}$, and verify that the projective limit topology is strong enough to imply an LDP result also with respect to our chosen topology. In Sect. 5.3 we conclude with the derivation of our main result by using the Dawson–Gärtner approach, adapted to our purposes: we have to use different distributions for the lower and the upper bound and will have to deal with some additional technical difficulty concerning the lower bound. Finally, Theorem 1.1 is implied by combining the results of Lemma 5.16, where we prove the upper bound, and Lemma 5.18, where we prove the lower bound.

The following objects will be fixed for the remainder of the section. Let $\mathcal{S}$ be a compact metric space. Fix a probability measure $\mu$ on $\mathcal{S}$. For any $N$ fix a type vector $\mathbf{x} = \mathbf{x}^{(N)} = (x_1, \ldots, x_N) \in \mathcal{S}^N$ such that its empirical measure $\mu_N = \frac{1}{N}\sum_{i=1}^N \delta_{x_i}$ weakly converges to $\mu$ as $N \to \infty$. Furthermore, let a continuous kernel $\kappa \colon \mathcal{S} \times \mathcal{S} \to [0, \infty)$ be given that is irreducible w.r.t. $\mu$, as well as a sequence of continuous kernels $\kappa_N$ that converges to $\kappa$ uniformly on $\mathcal{S} \times \mathcal{S}$ as $N \to \infty$. By $\mathcal{G}_N = \mathcal{G}([N], \mathbf{x}, \frac{1}{N}\kappa_N)$ we denote the inhomogeneous random graph introduced in Sect. 1.1 and by $(\mathcal{C}_i)_i$ we denote the vertex sets of the connected components of $\mathcal{G}_N$. We will denote the probability measure corresponding to the graph $\mathcal{G}_N$ by $\mathbb{P}_N$. Recall the definition of the microscopic and macroscopic empirical measures $\mathrm{Mi}_N$ and $\mathrm{Ma}_N$ given in (1.2). The goal is to derive the LDP for the pair $(\mathrm{Mi}_N, \mathrm{Ma}_N)$ that is formulated in Theorem 1.1.

Before chosing a discretization scheme, we would like to collect some properties of our state space $\mathcal{L} \times \mathcal{A}$. Recall that

$$\mathcal{L} := \{\lambda \in \mathcal{M}(\mathcal{M}_{\mathbb{N}_0}(\mathcal{S})) \colon c_\lambda \leq \mu \text{ or } c_\lambda \leq \mu_N \text{ for some } N \in \mathbb{N}, \lambda(\{0\}) = 0\}, \tag{5.1}$$

$$\mathcal{A} := \{\alpha \in \mathcal{M}_{\mathbb{N}_0}(\mathcal{M}(\mathcal{S})\setminus\{0\}) \colon c_\alpha \leq \mu \text{ or } c_\alpha \leq \mu_N \text{ for some } N \in \mathbb{N}\}, \tag{5.2}$$

and that both are equipped with vague topologies, i.e., a sequence $(\lambda_n)_{n\in\mathbb{N}}$ in $\mathcal{L}$ converges to $\lambda$, if for any continuous, compactly supported test function $g \colon \mathcal{M}_{\mathbb{N}_0}(\mathcal{S}) \to \mathbb{R}$ the integrals $\int g \, d\lambda_n$ converge to $\int g \, d\lambda$, as $n \to \infty$. In the same way, a sequence $(\alpha_n)_{n\in\mathbb{N}}$ in $\mathcal{A}$ converges to $\alpha$, if for any continuous, compactly supported test function $g \colon \mathcal{M}(\mathcal{S})\setminus\{0\} \to \mathbb{R}$ the integrals $\int g \, d\alpha_n$ converge to $\int g \, d\alpha$, as $n \to \infty$. Both on $\mathcal{M}_{\mathbb{N}_0}(\mathcal{S})$ and $\mathcal{M}(\mathcal{S})\setminus\{0\}$ we consider the topologies of weak convergence. In the next lemma, we give a short characterization of compactness that is implied by this choice; a verification is left to the reader.

**Lemma 5.1** *The following assertions hold:*

1. *A subset $\mathcal{N} \subset \mathcal{M}_{\mathbb{N}_0}(\mathcal{S})$ is compact if and only if it is closed and $\sup\{\nu(\mathcal{S}) \colon \nu \in \mathcal{N}\} < \infty$.*
2. *A subset $\mathcal{N} \subset \mathcal{M}_{\leq 1}(\mathcal{S})\setminus\{0\}$ is compact if and only if it is closed and $\inf\{\nu(\mathcal{S}) \colon \nu \in \mathcal{N}\} > 0$.*

The following lemma implies that any vague accumulation points of $\mathrm{Mi}_N$ and $\mathrm{Ma}_N$ are indeed elements of $\mathcal{L}$ and $\mathcal{A}$ if they exist.

**Lemma 5.2** *Both $\mathcal{L}$ and $\mathcal{A}$ are compact spaces.*

***Proof*** Denote $M := \{\mu\} \cup \{\mu_N \colon N \in \mathbb{N}\}$.

(1) Compactness of $\mathcal{L}$: For any $\lambda \in \mathcal{L}$ note that $|\lambda| \leq c_\lambda(\mathcal{S}) \leq \sup_{\tilde{\mu}\in M} \tilde{\mu}(\mathcal{S}) = 1$. The set $B_1 := \{\lambda \in \mathcal{M}(\mathcal{M}_{\mathbb{N}_0}(\mathcal{S})) \colon |\lambda| \leq 1\}$ is a bounded subset of the dual of $C_0(\mathcal{M}_{\mathbb{N}_0}(\mathcal{S})\setminus\{0\})$ (the space of continuous functions $g$ on $\mathcal{M}_{\mathbb{N}_0}(\mathcal{S})$ with $\lim_{k(\mathcal{S})\to\infty} g(k) = 0$), so by applying the Banach-Alaoglu Theorem we get that $B_1$ is compact w.r.t. the vague topology. Since $\mathcal{L} \subset B_1$ it remains to argue that $\mathcal{L}$ is closed.

Let $(\lambda_n)_{n\in\mathbb{N}}$ be a sequence in $\mathcal{L}$ with vague limit $\lambda$. Then for each $n \in \mathbb{N}$ there exists $\tilde{\mu}_n \in M$ such that $c_{\lambda_n} \leq \tilde{\mu}_n$. Since $M$ is compact w.r.t. the weak topology, we can find $\tilde{\mu} \in M$ and a subsequence (which we will also denote by $(\tilde{\mu}_n)_{n\in\mathbb{N}}$) such that $\tilde{\mu}_n \to \tilde{\mu}$ weakly, as $n \to \infty$. We now argue that this implies that $c_\lambda \leq \tilde{\mu}$. Let $f \colon \mathcal{S} \to [0, \infty)$ be a continuous and bounded function. For any $R \in \mathbb{N}$, let $\chi_R \colon [0, \infty) \to [0, \infty)$ be a smooth function satisfying $\mathbb{1}_{[0,R]} \leq \chi_R \leq \mathbb{1}_{[0,R+1]}$, such that $R \mapsto \chi_R$ is increasing pointwise. Define

$$\Phi_R^f \colon \mathcal{M}_{\mathbb{N}_0}(\mathcal{S}) \to \mathbb{R}, \quad \Phi_R^f(k) = \left(\int_\mathcal{S} f(s)\, k(\mathrm{d}s)\right) \chi_R(k(\mathcal{S})) \qquad (5.3)$$

and note that $\Phi_R^f(k) \nearrow \int f(s)\, k(\mathrm{d}s)$, as $R \to \infty$, pointwise for any $k$. By Lemma 5.1 we have that $\Phi_R^f$ is compactly supported and it can be easily seen that continuity of $\chi_R$ and $f$ imply that $\Phi_R^f$ is continuous. Therefore

$$\int_\mathcal{S} f(s)\, c_\lambda(\mathrm{d}s) = \sup_{R\in\mathbb{N}} \int \Phi_R^f(k)\, \lambda(\mathrm{d}k) = \sup_{R\in\mathbb{N}} \lim_{n\to\infty} \int \Phi_R^f(k)\, \lambda_n(\mathrm{d}k)$$
$$\leq \sup_{R\in\mathbb{N}} \lim_{n\to\infty} \int_\mathcal{S} f(s)\, c_{\lambda_n}(\mathrm{d}s) \leq \lim_{n\to\infty} \int_\mathcal{S} f(s)\, \tilde{\mu}_n(\mathrm{d}s) = \int_\mathcal{S} f(s)\, \tilde{\mu}(\mathrm{d}s).$$

Since this holds for any $f \geq 0$, we can conclude that $c_\lambda \leq \mu$ and thus $\lambda \in \mathcal{L}$.

(2) We only sketch the construction of a vague limit point. Fix a sequence $(\varepsilon_i)_{i\in\mathbb{N}}$, with $\varepsilon_i \searrow 0$ as $i \to \infty$, and define $N_{\varepsilon_i} := \{y \in \mathcal{M}_{\leq 1}(\mathcal{S})\setminus\{0\} \colon |y| \geq \varepsilon_i\}$. Note that $\varepsilon_i \alpha_n(N_{\varepsilon_i}) \leq c_{\alpha_n}(\mathcal{S}) \leq 1$ implies that $\alpha_n(N_{\varepsilon_i}) \leq 1/\varepsilon_i$. Thus, for fixed $i \in \mathbb{N}$, the restricted measures $(\alpha_n|_{N_{\varepsilon_i}})_{n\in\mathbb{N}}$ are bounded and thus have a vaguely converging sub-sequence. By diagonalization we can construct for each $i \in \mathbb{N}$ a subsequence of $(\alpha_n)_{n\in\mathbb{N}}$ such that the restrictions to $N_{\varepsilon_i}$ converge vaguely to some $\alpha^{(i)} \in \mathcal{M}(N_{\varepsilon_i})$ and in such a way that $\alpha^{(i+1)}|_{N_{\varepsilon_i}} = \alpha^{(i)}$ holds for all $i \in \mathbb{N}$. Thus the monotone limit $\lim_{i\to\infty} \alpha^{(i)} =: \alpha$ on $\mathcal{M}_{\leq 1}(\mathcal{S})\setminus\{0\}$ is a countably additive extension of the measures $\alpha^{(i)}$, $i \in \mathbb{N}$. Since every compactly supported test function has its support contained in $N_{\varepsilon_i}$ for some $i$ by Lemma 5.1 one sees that in the vague topology $\alpha_n \to \alpha$.

To see that $\alpha \in \mathcal{A}$ (i.e., that $\mathcal{A}$ is closed) we can proceed as in the proof for $\mathcal{L}$. Note that for $\varepsilon > 0$ we can find a smooth function $\chi_\varepsilon \colon [0, \infty) \to [0, \infty)$ satisfying $\mathbb{1}_{[2\varepsilon,\infty)} \leq \chi_\varepsilon \leq \mathbb{1}_{[\varepsilon,\infty)}$, which allows us to define $\Phi_\varepsilon^f$ as above but by truncating with $\chi_\varepsilon$. $\square$

### 5.1 Discretization and approximation

Let $(P_m)_{m\in\mathbb{N}}$ be a sequence of finite partitions of $\mathcal{S}$ into non-empty sets. For $m \in \mathbb{N}$ we denote $P_m = \{A_{m,i} \colon i = 1, \ldots, |P_m|\}$. We say that $(P_m)_{m\in\mathbb{N}}$ is nested if for $m \leq n$ and any $A_{m,i} \in P_m$ there is $J \subset \{1, \ldots, |P_n|\}$, such that $A_{m,i} = \bigcup_{j\in J} A_{n,j}$. For any subset $A$ we write $\mathrm{diam}(A) := \sup\{d(x,y) \colon x, y \in A\}$. For any measure $\nu$ on $\mathcal{S}$ a set $A \subset \mathcal{S}$ is called a continuity set of $\nu$ if $\nu(\partial A) = 0$. The following lemma is a modification of Lemma 7.1 from [7].

**Lemma 5.3** *There exists a sequence of finite partitions $(P_m)_{m\in\mathbb{N}}$ of $\mathcal{S}$ with the following properties:*

1. *For any $m \in \mathbb{N}$ and any $i = 1, \ldots, |P_m|$ we have that $A_{m,i}$ is measurable and a continuity set of $\mu$ and $\mu_N$ for all $N \in \mathbb{N}$.*
2. *$(P_m)_{m\in\mathbb{N}}$ is nested.*
3. *It holds that*

$$\lim_{m\to\infty} \max_{i=1,\ldots,|P_m|} diam(A_{m,i}) = 0. \tag{5.4}$$

*Proof* The proof can be done as in the proof of Lemma 7.1 in [7] with two small modifications. Let $A \subset \mathcal{S}$ be the set of points that are atoms of $\mu$ or $\mu_N$ for some $N \in \mathbb{N}$. Then $A$ is still a countable set, so we may pick the balls $B_{mi}$ in a way such that they are continuity sets of $\mu$ and $\mu_N$ for all $N \in \mathbb{N}$. Since in our case the set $\mathcal{S}$ is compact, we can cover it with finitely many of these balls, hence we get the stronger property formulated in (5.4). □

In the following we always assume that $(P_m)_{m\in\mathbb{N}}$ has all the properties given in Lemma 5.3. For $m \in \mathbb{N}$ and any $A_{m,i} \in P_m$ we pick exactly one point $x_{m,i}$ from the set $A_m$ which we call the representative of $A_{m,i}$. We define

$$\mathcal{S}_m := \{x_{m,i} : i = 1, \ldots, |P_m|\} \tag{5.5}$$

and define the projection

$$\pi_m : \mathcal{S} \to \mathcal{S}_m, \quad x \mapsto x_{m,i} \tag{5.6}$$

where $i$ is such that $A_{m,i}$ is the unique set containing $x$.

Further we can lift the projection to the space $\mathcal{M}(\mathcal{S})$, the space of finite measures on $\mathcal{S}$. For any $m \in \mathbb{N}$ we introduce (by abuse of notation)

$$\pi_m : \mathcal{M}(\mathcal{S}) \to \mathcal{M}(\mathcal{S}_m), \quad \pi_m(\nu) := \nu \circ \pi_m^{-1} \quad \text{for } \nu \in \mathcal{M}(\mathcal{S}). \tag{5.7}$$

Going one level further, we also define

$$\pi_m : \mathcal{M}(\mathcal{M}(\mathcal{S})) \to \mathcal{M}(\mathcal{M}(\mathcal{S}_m)), \quad \pi_m(\lambda) := \lambda \circ \pi_m^{-1} \quad \text{for } \lambda \in \mathcal{M}(\mathcal{M}(\mathcal{S})). \tag{5.8}$$

We can now apply the projections to all the levels of our graph setting. On the first level of discretization, i.e., from $\mathcal{S}$ to $\mathcal{S}_m$, we approximate the type of a vertex by some type from the discrete set $\mathcal{S}_m$. On the second level, we approximate the type configuration of a vertex set which, depending on the context, may or may not be a cluster. On the third level, we approximate measures that count the multiplicities of type configurations and are therefore suited to register the number of clusters described by the different type configurations.

Fix $m \in \mathbb{N}$. We write $\pi_m(\mathbf{x}) = (\pi_m(x_1), \ldots, \pi_m(x_N)) \in (\mathcal{S}_m)^N$ for the discretized type sequence and denote by $\mu_N^{(m)} = \frac{1}{N}\sum_{i=1}^N \delta_{\pi_m(x_i)}$ its empirical measure. It is easy to check that $\mu_N^{(m)} = \pi_m(\mu_N)$. Since our partition $P_m$ has carefully been chosen such that the sets $A_{m,i} \in P_m$ are continuity sets of $\mu$, it holds that $\mu_N^{(m)}$ converges weakly to $\mu^{(m)} := \pi_m(\mu)$, as $N \to \infty$. We also define $\eta_{\pi_m(\mathbf{x})} \colon \mathcal{P}([N]) \to \mathcal{M}_{\mathbb{N}_0}(\mathcal{S}_m)$ as the discrete analog of the type registering measure $\eta_\mathbf{x}$ from (1.1), i.e., for any $A \subset [N]$ we have that $\eta_{\pi_m(\mathbf{x})} = \sum_{i \in A} \delta_{\pi_m(x_i)} \in \mathcal{M}_{\mathbb{N}_0}(\mathcal{S}_m)$. Abbreviating $\eta_m = \eta_{\pi_m(\mathbf{x})}$ we can now define

$$\mathrm{Mi}_N^{(m)} = \frac{1}{N}\sum_i \delta_{\eta_m(\mathcal{C}_i)} \quad \text{and} \quad \mathrm{Ma}_N^{(m)} = \sum_i \delta_{\frac{1}{N}\eta_m(\mathcal{C}_i)}. \tag{5.9}$$

It is straightforward to show that $\mathrm{Mi}_N^{(m)} = \pi_m(\mathrm{Mi}_N)$ and $\mathrm{Ma}_N^{(m)} = \pi_m(\mathrm{Ma}_N)$. Both empirical measures $\mathrm{Mi}_N^{(m)}$ and $\mathrm{Ma}_N^{(m)}$ evaluate type information given by the vertices of the clusters only roughly by approximating the type of each vertex by an element from $\mathcal{S}_m$. Note that $c_{\mathrm{Mi}_N^{(m)}} = \pi_m(c_{\mathrm{Mi}_N}) \leq \mu_N^{(m)}$ and $c_{\mathrm{Ma}_N^{(m)}} = \pi_m(c_{\mathrm{Ma}_N}) \leq \mu_N^{(m)}$. Hence, the natural state spaces for the discretized empirical measures are given by

$$\mathcal{L}_m = \{\lambda \in \mathcal{M}(\mathcal{M}_{\mathbb{N}_0}(\mathcal{S}_m)) \colon c_\lambda \leq \mu^{(m)} \text{ or } c_\lambda \leq \mu_N^{(m)} \text{ for some } N \in \mathbb{N}, \lambda(\{0\}) = 0\}, \tag{5.10}$$

$$\mathcal{A}_m = \{\alpha \in \mathcal{M}_{\mathbb{N}_0}(\mathcal{M}(\mathcal{S}_m)\setminus\{0\}) \colon c_\alpha \leq \mu^{(m)} \text{ or } c_\alpha \leq \mu_N^{(m)} \text{ for some } N \in \mathbb{N}\}, \tag{5.11}$$

and we endow them with vague topologies.

At this point it is important to note that we cannot apply the discrete LDP from Theorem 3.1 directly to get an LDP for the measure $\mathbb{P}_N((\mathrm{Mi}_N^{(m)}, \mathrm{Ma}_N^{(m)}) \in \cdot)$, since the edges of the random graph $\mathcal{G}_N$ are drawn according to the non-discretized types of the vertices. Before applying Theorem 3.1 one has to approximate the underlying graph itself by a discrete version.

Let $\kappa_N^{(m)} \colon \mathcal{S}_m \times \mathcal{S}_m \to [0, \infty)$, $N \in \mathbb{N}$, be a sequence of kernels on $\mathcal{S}_m$. We will specify later, in which sense they should be an approximation for $\kappa_N$. Consider the inhomogeneous random graph $\mathcal{G}_N^{(m)} = \mathcal{G}([N], \pi_m(\mathbf{x}), \frac{1}{N}\kappa_N^{(m)})$ and denote by $(\mathcal{C}_i^{(m)})_i$ the collection of the vertex sets of its connected components.

Instead of choosing just one approximating kernel, we will consider a lower and an upper approximation for $\kappa_N$, which will allow us to find upper and lower bounds for the distribution $\mathbb{P}_N((\mathrm{Mi}_N^{(m)}, \mathrm{Ma}_N^{(m)}) \in \cdot)$. For fixed $m \in \mathbb{N}$ let $\kappa_N^{(m,-)}$ and $\kappa_N^{(m,+)}$, $N \in \mathbb{N}$, be two sequences of kernels on $\mathcal{S}_m$ satisfying

$$\kappa_N^{(m,-)}(\pi_m(x), \pi_m(\hat{x})) \leq \kappa_N(x, \hat{x}) \leq \kappa_N^{(m,+)}(\pi_m(x), \pi_m(\hat{x})), \quad \forall x, \hat{x} \in \mathcal{S}, \forall N \in \mathbb{N}. \tag{5.12}$$

An obvious choice is given by

$$\kappa_N^{(m,-)}(r,s) := \inf\{\kappa_N(x,\hat{x}) \colon x, \hat{x} \in \mathcal{S},\, \pi_m(x) = r,\, \pi_m(\hat{x}) = s\}, \quad \text{for } r,s \in \mathcal{S}_m \tag{5.13}$$

$$\kappa_N^{(m,+)}(r,s) := \sup\{\kappa_N(x,\hat{x}) \colon x, \hat{x} \in \mathcal{S},\, \pi_m(x) = r,\, \pi_m(\hat{x}) = s\}, \quad \text{for } r,s \in \mathcal{S}_m. \tag{5.14}$$

Defining $\kappa^{(m,*)} = \lim_{N\to\infty} \kappa_N^{(m,*)}$ for both $* \in \{+,-\}$ it is obvious that

$$\kappa^{(m,-)}(\pi_m(x), \pi_m(\hat{x})) \leq \kappa(x,\hat{x}) \leq \kappa^{(m,+)}(\pi_m(x), \pi_m(\hat{x})), \quad \text{for all } x, \hat{x} \in \mathcal{S}, \tag{5.15}$$

and we see that

$$\lim_{m\to\infty} \kappa^{(m,*)}(\pi_m(x), \pi_m(\hat{x})) = \kappa(x,\hat{x}), \quad \text{uniformly in } x, \hat{x} \in \mathcal{S},\, * \in \{+,-\}. \tag{5.16}$$

For both $* \in \{+,-\}$ we will denote by $\mathbb{P}_N^{(m,*)}$ the probability measure corresponding to the graph $\mathcal{G}_N^{(m,*)} = \mathcal{G}([N], \pi_m(\mathbf{x}), \frac{1}{N}\kappa_N^{(m,*)})$.

The following comparison lemma shows how we can estimate the distribution $\mathbb{P}_N((\mathrm{Mi}_N^{(m)}, \mathrm{Ma}_N^{(m)}) \in \cdot)$ from below and above.

**Lemma 5.4** *Fix $m \in \mathbb{N}$ and let $\kappa_N^{(m,-)}$ and $\kappa_N^{(m,+)}$, $N \in \mathbb{N}$, be two sequences of kernels on $\mathcal{S}_m$ satisfying* (5.12). *Assume that $N \in \mathbb{N}$ is large enough such that $\frac{1}{N}\kappa_N^{(m,+)} \leq 1$. Then, for any $\ell \in \mathcal{M}_{\mathbb{N}_0}(\mathcal{M}_{\mathbb{N}_0}(\mathcal{S}_m))$,*

$$(\Delta_N^{(m)})^{-1} \mathbb{P}_N^{(m,-)}(N\mathrm{Mi}_N = \ell) \leq \mathbb{P}_N(N\mathrm{Mi}_N^{(m)} = \ell) \leq \mathbb{P}_N^{(m,+)}(N\mathrm{Mi}_N = \ell) \Delta_N^{(m)}, \tag{5.17}$$

*where*

$$\Delta_N^{(m)} = \prod_{r,s \in \mathcal{S}_m} \left( \frac{1 - \frac{1}{N}\kappa_N^{(m,-)}(r,s)}{1 - \frac{1}{N}\kappa_N^{(m,+)}(r,s)} \right)^{\frac{1}{2}N^2 \mu_N^{(m)}(r) \mu_N^{(m)}(s)}. \tag{5.18}$$

**Proof** As in the discrete setting of Sect. 3.1 we can identify elements from $\mathcal{M}_{\mathbb{N}_0}(\mathcal{S}_m)$ with elements from $\mathbb{N}_0^{\mathcal{S}_m}$ and identify $\ell$ with $(\ell_k)_{k \in \mathbb{N}_0^{\mathcal{S}_m}}$. For $* \in \{+,-\}$ and any $k \in \mathbb{N}_0^{\mathcal{S}_m}$ we write

$$p_N^{(m,*)}(k) = \mathbb{P}_N^{(m,*)}\big(\mathcal{G}(|k|, \mathbf{x}_k, \tfrac{1}{N}\kappa_N^{(m,*)}) \text{ is connected}\big), \tag{5.19}$$

with $\mathbf{x}_k \in \mathcal{S}_m^{|k|}$ any $|k|$-dimensional vector compatible with $k$. The main idea is that we identify in the exact formula (3.15) in Lemma 3.3 those terms that are increasing in $\kappa$ and those that are decreasing in $\kappa$. Indeed, the connection probabilities are increasing in $\kappa$ since the event of being connected is increasing in the edge parameter. Indeed, for any $k \in \mathbb{N}_0^{\mathcal{S}_m}$ and $\mathbf{x} \in \mathcal{S}^{|k|}$ such that $\pi_m(\sum_{i=1}^{|k|} \delta_{x_i}) = k$,

$$p_N^{(m,-)}(k) \leq \mathbb{P}_N(\mathcal{G}(|k|, \mathbf{x}, \tfrac{1}{N}\kappa_N) \text{ is connected }) \leq p_N^{(m,+)}(k).$$

Furthermore, the powers of $(1 - \frac{1}{N}\kappa(\cdot))$ that describe the probabilities of not being connected are increasing for negative powers, and decreasing for positive powers. With those observations and combinatorial factors we can estimate

$$\mathbb{P}_N(N\mathrm{Mi}_N^{(m)} = \ell) \leq \left(\prod_{r \in \mathcal{S}_m} (N\mu_N^{(m)}(r))!\right) \left(\prod_{k \in \mathbb{N}_0^{\mathcal{S}_m}} \frac{p_N^{(m,+)}(k)^{\ell_k}}{\ell_k! \prod_{r \in \mathcal{S}_m}(k_r!)^{\ell_k}}\right)$$
$$\times \left(\prod_{r,s \in \mathcal{S}_m} \left(1 - \tfrac{1}{N}\kappa_N^{(m,+)}(r,s)\right)^{\frac{1}{2}\sum_k k_r(N\mu_N^{(m)}(s) - k_s)\ell_k}\right) \Delta_N^{(m)}$$
(5.20)

The lower bound and its proof are analogous. $\square$

### 5.2 Projective system

Using our discretization scheme from Sect. 5.1 we now introduce a projective system that fits into the framework of [19, Section 4.6]. Recall the notions introduced at the beginning of Sect. 5.1, in particular the family of finite partitions $(P_m)_{m \in \mathbb{N}}$ that has the properties formulated in Lemma 5.3 and the definition (5.5) of the discretized type spaces $\mathcal{S}_m$. For any pair $m, n \in \mathbb{N}$ with $m \leq n$ we define

$$\pi_{m,n} \colon \mathcal{S}_n \to \mathcal{S}_m, \quad x_{n,j} \mapsto x_{m,i} \tag{5.21}$$

where $x_{m,i}$ is the representative of $A_{m,i}$ and $A_{m,i}$ is the unique set in $P_m$ containing $x_{n,j}$. As before, we lift this definition to the measure spaces by defining (with abuse of notation)

$$\pi_{m,n} \colon \mathcal{M}(\mathcal{S}_n) \to \mathcal{M}(\mathcal{S}_m), \quad \pi_{m,n}(\nu) = \nu \circ \pi_{m,n}^{-1}, \tag{5.22}$$

as well as

$$\pi_{m,n} \colon \mathcal{M}(\mathcal{M}(\mathcal{S}_n)) \to \mathcal{M}(\mathcal{M}(\mathcal{S}_m)), \quad \pi_{m,n}(\lambda) = \lambda \circ \pi_{m,n}^{-1}. \tag{5.23}$$

It is easy to check that restricting the latter mapping to $\mathcal{L}_n$ and $\mathcal{A}_n$, respectively, gives us two well-defined mappings $\pi_{m,n} \colon \mathcal{L}_n \to \mathcal{L}_m$ and $\pi_{m,n} \colon \mathcal{A}_n \to \mathcal{A}_m$, respectively. Now, $(\mathcal{L}_m, \pi_{m,n})_{m \leq n}$ (or $(\mathcal{A}_m, \pi_{m,n})_{m \leq n}$) is called a projective system, if

– for any $m \in \mathbb{N}$ the space $\mathcal{L}_m$ is a Hausdorff topological space,

- for any $m, n \in \mathbb{N}$ with $m \leq n$ the mapping $\pi_{m,n} \colon \mathcal{L}_n \to \mathcal{L}_m$ is continuous,
- for any $m, n, p \in \mathbb{N}$ with $m \leq n \leq p$ we have $\pi_{m,p} = \pi_{m,n} \circ \pi_{n,p}$.

The projective limit of the projective system $(\mathcal{L}_m, \pi_{m,n})_{m \leq n}$ is denoted by $\varprojlim \mathcal{L}_m$ and is defined as the subset of the product space $\prod_{m \in \mathbb{N}} \mathcal{L}_m$ that contains all elements $(\lambda_m)_{m \in \mathbb{N}}$ that satisfy $\pi_{m,n}(\lambda_n) = \lambda_m$ for any pair $m, n \in \mathbb{N}$ with $m \leq n$. The projective limit $\varprojlim \mathcal{L}_m$ is equipped with the topology that is induced by the product topology on $\prod_{m \in \mathbb{N}} \mathcal{L}_m$. We call this the projective limit topology. In particular, a sequence $\lambda^{(n)} \in \varprojlim \mathcal{L}_m$ converges to $\lambda \in \varprojlim \mathcal{L}_m$ as $n \to \infty$ if it holds for any $m \in \mathbb{N}$ that $\lambda_m^{(n)}$ converges to $\lambda_m$ as $n \to \infty$.

**Lemma 5.5** $(\mathcal{L}_m, \pi_{m,n})_{m \leq n}$ and $(\mathcal{A}_m, \pi_{m,n})_{m \leq n}$ are projective systems.

*Proof* The fact that $\mathcal{L}_m$ and $\mathcal{A}_m$ are Hausdorff topological spaces is a consequence of the equivalence of their topologies with the discrete topologies that we described in Sect. 3.1. For $m \leq n$, the continuity of $\pi_{m,n}$ is easily verified for the lowest level, i.e., on $\mathcal{S}_n \to \mathcal{S}_m$, and then lifted to the higher levels. Indeed, for each open set $\mathcal{O} \in \mathcal{M}(\mathcal{S}_m)$ (which can be identified as an open set in $\mathbb{R}_{\geq 0}^{\mathcal{S}_m} \setminus \{0\}$), we see that $\pi_{m,n}^{-1}(\mathcal{O})$ is an open set in $\mathcal{M}(\mathcal{S}_m)$. The same is true at the higher level with open sets in $\mathcal{L}_m$ and $\mathcal{A}_m$, which are images of open sets in $\mathcal{L}_n$ and $\mathcal{A}_n$. □

Throughout this section we will denote $\mathcal{L}_\infty = \varprojlim \mathcal{L}_m$ and $\mathcal{A}_\infty = \varprojlim \mathcal{A}_m$. The aim of this section is to prove the following

**Proposition 5.6** *The following assertions hold.*

1. *The set $\mathcal{L} \times \mathcal{A}$ can be identified with $\mathcal{L}_\infty \times \mathcal{A}_\infty$.*
2. *The projective limit topology on $\mathcal{L}_\infty \times \mathcal{A}_\infty$ is equivalent to the vague topology on $\mathcal{L} \times \mathcal{A}$.*

The proof will be a consequence of the following lemmas. In Lemma 5.8 we explain how to project from $\mathcal{L} \times \mathcal{A}$ to $\mathcal{L}_\infty \times \mathcal{A}_\infty$. As a direct consequence of Lemma 5.9 we get that this operation is continuous. Afterwards, we deal with the inverse operation. In Lemma 5.10 we show the existence of the inverse and in Lemma 5.12 we argue that it is continous.

Here is a basic property of the mappings $\pi_m$ and $\pi_{m,n}$, which we need to prepare for Lemma 5.8.

**Lemma 5.7** *Let $m \leq n$. Then the equality $\pi_m = \pi_{m,n} \circ \pi_n$ holds on all levels, i.e., for all mappings that were defined in* (5.6)–(5.8) *and* (5.21)–(5.23).

*Proof* On the lowest level, i.e., on $\mathcal{S} \to \mathcal{S}_m$, the equality is a direct consequence of the fact that $(P_m)_{m \in \mathbb{N}}$ is nested. On the higher levels, the equality follows by the definition of the image measure. □

**Lemma 5.8** *The following holds.*

1. *Let $\lambda \in \mathcal{L}$. Then for any $m \in \mathbb{N}$ we have that $\pi_m(\lambda) \in \mathcal{L}_m$. Further, we have that the sequence $(\pi_m(\lambda))_{m \in \mathbb{N}}$ is an element of the projective limit $\mathcal{L}_\infty$.*

2. Let $\alpha \in \mathcal{A}$. Then for any $m \in \mathbb{N}$ we have that $\pi_m(\alpha) \in \mathcal{A}_m$. Further, we have that the sequence $(\pi_m(\alpha))_{m \in \mathbb{N}}$ is an element of the projective limit $\mathcal{A}_\infty$.

Consequently, the mapping $\Pi: \mathcal{L} \times \mathcal{A} \to \mathcal{L}_\infty \times \mathcal{A}_\infty$, $(\lambda, \alpha) \mapsto \big((\pi_m(\lambda))_{m \in \mathbb{N}}, (\pi_m(\alpha)_{m \in \mathbb{N}})\big)$ is well-defined.

**Proof** (1) For any $\lambda \in \mathcal{L}$ and $m \in \mathbb{N}$ we have that $c_{\pi_m(\lambda)} = \pi_m(c_\lambda) \leq \pi_m(\tilde{\mu})$, where $\tilde{\mu} \in \{\mu\} \cup \{\mu_N : N \in \mathbb{N}\}$, hence $\pi_m(\lambda) \in \mathcal{L}_m$. As a direct consequence of Lemma 5.7 the sequence $(\pi_m(\lambda))_{m \in \mathbb{N}}$ satisfies the consistency condition and is therefore an element of $\mathcal{L}_\infty$.

(2) For $\alpha \in \mathcal{A}$ the proof is analogous. □

**Lemma 5.9** *Fix any $m \in \mathbb{N}$. The mappings $\pi_m: \mathcal{L} \to \mathcal{L}_m$ and $\pi_m: \mathcal{A} \to \mathcal{A}_m$ are continuous w.r.t. the vague topologies. Consequently, the mapping $\Pi$ defined as in Lemma 5.8 is continuous.*

**Proof** On the lowest level, i.e., for $\pi_m: \mathcal{S} \to \mathcal{S}_m$, it is obvious that $\pi_m$ is continuous on the set $\mathcal{S} \setminus \bigcup_{i=1}^{|P_m|} \partial A_{m,i}$. Consider the second level, i.e., $\pi_m: \mathcal{M}(\mathcal{S}) \to \mathcal{M}(\mathcal{S}_m)$. We claim that $\pi_m$ is continuous on

$$\mathcal{N}_0 := \{\nu \in \mathcal{M}(\mathcal{S}) : \nu(\partial A_{m,i}) = 0 \text{ for all } i = 1, \ldots, |P_m|\} \quad (5.24)$$

w.r.t. weak convergence: Given some $\nu \in \mathcal{N}_0$ and a sequence $(\nu_n)_{n \in \mathbb{N}}$ in $\mathcal{M}(\mathcal{S})$ such that $\nu_n \to \nu$ weakly as $n \to \infty$ and some continuous bounded function $f: \mathcal{S}_m \to \mathbb{R}$ we clearly have that

$$\int_{\mathcal{S}_m} f \, d\pi_m(\nu_n) = \int_{\mathcal{S}} f \circ \pi_m \, d\nu_n \to \int_{\mathcal{S}} f \circ \pi_m \, d\nu = \int_{\mathcal{S}_m} f \, d\pi_m(\nu),$$

as $n \to \infty$, since $f \circ \pi_m$ is continuous $\nu$-almost everywhere. Hence, $\pi_m(\nu_n) \to \pi_m(\nu)$ weakly as $n \to \infty$.

Now, consider $\pi_m: \mathcal{L} \to \mathcal{L}_m$. For $\lambda \in \mathcal{L}$ we have that $c_\lambda \leq \tilde{\mu}$, where $\tilde{\mu} \in \{\mu\} \cup \{\mu_N : N \in \mathbb{N}\}$. This implies that for any $i = 1, \ldots, |P_m|$ we have that

$$\lambda\big(\{k: k(\partial A_{m,i}) > 0\}\big) \leq \sum_{n=1}^{\infty} n\lambda\big(\{k: k(\partial A_{m,i}) = n\}\big) = c_\lambda(\partial A_{m,i}) \leq \tilde{\mu}(\partial A_{m,i}) = 0 \quad (5.25)$$

by Lemma (5.3). With other words, $\lambda$ is concentrated on a subset of the set $\mathcal{N}_0$ given in (5.24). Now, let $(\lambda_n)_{n \in \mathbb{N}}$ be a sequence in $\mathcal{L}$ that converges vaguely to $\lambda$ as $n \to \infty$. Then for any function $g: \mathcal{M}_{\mathbb{N}_0}(\mathcal{S}_m) \to \mathbb{R}$ that is continuous and compactly supported, we have that

$$\int_{\mathcal{M}_{\mathbb{N}_0}(\mathcal{S}_m)} g \, d\pi_m(\lambda_n) = \int_{\mathcal{M}_{\mathbb{N}_0}(\mathcal{S})} g \circ \pi_m \, d\lambda_n \to \int_{\mathcal{M}_{\mathbb{N}_0}(\mathcal{S})} g \circ \pi_m \, d\lambda$$
$$= \int_{\mathcal{M}_{\mathbb{N}_0}(\mathcal{S}_m)} g \, d\pi_m(\lambda),$$

as $n \to \infty$, since $g \circ \pi_m$ is compactly supported and continuous $\lambda$-almost everywhere. Hence, $\pi_m(\lambda_n) \to \pi_m(\lambda)$ vaguely as $n \to \infty$.

The proof for $\pi_m \colon \mathcal{A} \to \mathcal{A}_m$ is analogous. □

In the next lemmas we will deal with the construction of the inverse of the projection mapping $\Pi$ and verify its continuity. This requires, for any $m \in \mathbb{N}$, to identify measures $\lambda_m \in \mathcal{L}_m$ with measures $\bar{\lambda}_m \in \mathcal{L}$. To prepare for this identification we now define for any $m \in \mathbb{N}$

$$\mathcal{M}_m(\mathcal{S}) := \{\nu \in \mathcal{M}(\mathcal{S}) \colon \nu \text{ is concentrated on } \mathcal{S}_m\}, \qquad (5.26)$$

where for any measure $\nu$ on some measure space $\mathcal{X}$ and any measurable $\mathcal{U} \subset \mathcal{X}$ we say that $\nu$ is concentrated on $\mathcal{U}$ if $\nu(\mathcal{X}\setminus\mathcal{U}) = 0$. It is clear that $\mathcal{M}_m(\mathcal{S})$ can be identified with $\mathcal{M}(\mathcal{S}_m)$. Observe that this is possible because $\mathcal{S}_m \subset \mathcal{S}$, as we have defined it via the representatives of each partition. For any $\nu_m \in \mathcal{M}(\mathcal{S}_m)$ we will denote the corresponding element by $\bar{\nu}_m \in \mathcal{M}_m(\mathcal{S})$ and we will write $\bar{\pi}_m \colon \mathcal{M}(\mathcal{S}) \to \mathcal{M}_m(\mathcal{S})$ for the mapping that we obtain by concatenating $\pi_m \colon \mathcal{M}(\mathcal{S}) \to \mathcal{M}(\mathcal{S}_m)$ as defined in (5.7) with the operation $\nu_m \mapsto \bar{\nu}_m$. Then we can identify the space $\mathcal{L}_m$ that was defined in (5.10) with

$$\bar{\mathcal{L}}_m := \{\lambda \in \mathcal{L} \colon \lambda \text{ is concentrated on } \mathcal{M}_m(\mathcal{S}) \cap \mathcal{M}_{\mathbb{N}_0}(\mathcal{S})\}. \qquad (5.27)$$

For any $\lambda_m \in \mathcal{L}_m$ we will denote the corresponding element by $\bar{\lambda}_m \in \bar{\mathcal{L}}_m$ and we will write $\bar{\pi}_m \colon \mathcal{L} \to \bar{\mathcal{L}}_m$ for the mapping that we obtain by concatenating $\pi_m \colon \mathcal{L} \to \mathcal{L}_m$ with the operation $\lambda_m \mapsto \bar{\lambda}_m$. In the same way, we identify the space $\mathcal{A}_m$ that was defined in (5.11) with

$$\bar{\mathcal{A}}_m := \{\alpha \in \mathcal{A} \colon \alpha \text{ is concentrated on } \mathcal{M}_m(\mathcal{S}) \cap \mathcal{M}(\mathcal{S})\setminus\{0\}\}. \qquad (5.28)$$

Now we construct the inverse of the projection mapping $\Pi$ that was defined in Lemma 5.8.

**Lemma 5.10** *The following assertions hold.*

1. *Let $(\lambda_m)_{m \in \mathbb{N}} \in \mathcal{L}_\infty$. Then there exists a unique $\lambda \in \mathcal{L}$ such that $\lambda_m = \pi_m(\lambda)$ holds for all $m \in \mathbb{N}$.*
2. *Let $(\alpha_m)_{m \in \mathbb{N}} \in \mathcal{A}_\infty$. Then there exists a unique $\alpha \in \mathcal{A}$ such that $\alpha_m = \pi_m(\alpha)$ holds for all $m \in \mathbb{N}$.*

*Consequently, the mapping $\Pi$ defined in Lemma 5.8 is bijective with inverse $\Pi^{-1}$.*

**Proof** Fix $(\lambda_m)_{m \in \mathbb{N}} \in \mathcal{L}_\infty$. The idea is to identify for any $m \in \mathbb{N}$ the measure $\lambda_m \in \mathcal{L}_m$ uniquely with the element $\bar{\lambda}_m \in \bar{\mathcal{L}}_m$ and to prove that the sequence $(\bar{\lambda}_m)_{m \in \mathbb{N}}$ has a vague limit point in $\mathcal{L}$, which we will denote by $\lambda$. It then remains to show that $\pi_m(\lambda) = \lambda_m$ holds for any $m \in \mathbb{N}$.

Next, we will argue for the existence of a vague limit point of $(\bar{\lambda}_m)_{m \in \mathbb{N}}$. As an element of $\mathcal{L}_\infty$ the sequence $(\lambda_m)_{m \in \mathbb{N}}$ satisfies the consistency condition $\lambda_m = \pi_{m,n}(\lambda_n)$ for any $m \leq n$. Abbreviating $\mathcal{M}_{\mathbb{N}_0} := \mathcal{M}_{\mathbb{N}_0}(\mathcal{S}) \setminus \{0\}$ and $\mathcal{M}_{\mathbb{N}_0, m} := \mathcal{M}_{\mathbb{N}_0}(\mathcal{S}_m) \setminus \{0\}$, we get that

$$\bar{\lambda}_m(\mathcal{M}_{\mathbb{N}_0}) = \lambda_m(\mathcal{M}_{\mathbb{N}_0,m}) = \lambda_n(\pi_{m,n}^{-1}(\mathcal{M}_{\mathbb{N}_0,m})) = \lambda_n(\mathcal{M}_{\mathbb{N}_0,n}) = \bar{\lambda}_n(\mathcal{M}_{\mathbb{N}_0})$$

and consequently, $\bar{\lambda}_m$, $m \in \mathbb{N}$, is of constant total variation. Note that the measures $\bar{\lambda}_m$, $m \in \mathbb{N}$, are in the dual of $C_0(\mathcal{M}_{\mathbb{N}_0})$ and that, due to the Banach–Alaoglu theorem, they are compact w.r.t. the weak*-topology, which implies compactness w.r.t. the vague topology on $\mathcal{L}$. Hence, there exists a vague limit point $\lambda \in \mathcal{M}(\mathcal{M}_{\mathbb{N}_0})$ and a subsequence $(\bar{\lambda}_{m_i})_{i \in \mathbb{N}}$ in $\mathcal{L}$ converging vaguely to $\lambda$. Since $\mathcal{L}$ is compact by Lemma 5.2, we also have that $\lambda \in \mathcal{L}$.

Next, we fix $m \in \mathbb{N}$ and our goal is to show that $\lambda_m = \pi_m(\lambda)$. Observe that as a consequence of our identification between $\lambda_n$ and $\bar{\lambda}_n$ we have that $\pi_n(\bar{\lambda}_n) = \lambda_n$ for any $n \in \mathbb{N}$. Together with the consistency and Lemma 5.7 we get for $n \geq m$ that

$$\lambda_m = \pi_{m,n}(\lambda_n) = \pi_{m,n}(\pi_n(\bar{\lambda}_n)) = \pi_m(\bar{\lambda}_n).$$

Choosing a subsequence $(\bar{\lambda}_{n_i})_{i \in \mathbb{N}}$ that converges vaguely to $\lambda$, we get that $\lambda_m = \lim_{i \to \infty} \pi_m(\bar{\lambda}_{n_i}) = \pi_m(\lambda)$ where we used the continuity of the mapping $\pi_m$ that we showed in Lemma 5.9.

For a given $(\alpha_m)_{m \in \mathbb{N}} \in \mathcal{A}_\infty$ the proof is analogous. □

To prepare for the proof of the continuity of $\Pi^{-1}$ we need the following lemma.

**Lemma 5.11** *On bounded subsets of $\mathcal{M}(\mathcal{S})$, the mapping $\mathcal{M}(\mathcal{S}) \to \mathcal{M}(\mathcal{S})$, $\nu \mapsto \bar{\pi}_m(\nu)$, converges uniformly to the identity as $m \to \infty$.*

**Proof** Recall that we equip $\mathcal{M}(\mathcal{S})$ with the weak topology, which is generated by all the test integrals against continuous bounded functions $\mathcal{S} \to \mathbb{R}$. The weak topology on $\mathcal{M}(\mathcal{S})$ admits a number of metrisations, especially since $\mathcal{S}$ is compact. We introduce the dual bounded-Lipschitz distance given by

$$d_{\mathrm{BL}}(\nu, \hat{\nu}) := \sup_{\phi: \, \|\phi\|_{\mathrm{BL}} \leq 1} \left| \int_\mathcal{S} \phi \, \mathrm{d}\nu - \int_\mathcal{S} \phi \, \mathrm{d}\hat{\nu} \right|, \quad (5.29)$$

where $\|\phi\|_{\mathrm{BL}} = \|\phi\|_\infty + \mathrm{Lip}(\phi)$ and $\mathrm{Lip}(\phi)$ is the infimum of all Lipschitz constants of $\phi: \mathcal{S} \to \mathbb{R}$. Let $\phi: \mathcal{S} \to \mathbb{R}$ satisfy $\|\phi\|_{\mathrm{BL}} \leq 1$ then

$$\left| \int_\mathcal{S} \phi \, \mathrm{d}\pi_m(\nu) - \int_\mathcal{S} \phi \, \mathrm{d}\nu \right| \leq \sum_{i=1}^{|P_m|} \int_{A_{m,i}} \left| \phi(x_{m,i}) - \phi(x) \right| \nu(\mathrm{d}x)$$
$$\leq \max_{i=1,\ldots,|P_m|} \mathrm{diam}(A_{m,i}) \nu(\mathcal{S}), \quad (5.30)$$

where $x_{m,i}$ is the representative of $A_{m,i}$ as defined right before definition (5.5). Now, if $\mathcal{N} \subset \mathcal{M}(\mathcal{S})$ is bounded, i.e., $\sup_{\nu \in \mathcal{N}} \nu(\mathcal{S}) < \infty$, then $d_{\mathrm{BL}}(\nu, \pi_m(\nu))$ vanishes as $m \to \infty$ uniformly on $\mathcal{N}$ by assumption (5.4). □

In the next lemma we verify that the mapping $\Pi^{-1}$ constructed in Lemma 5.10 is continuous.

**Lemma 5.12** *The following assertions hold.*

1. *Let $\lambda^{(n)}, n \in \mathbb{N}$, be a sequence in $\mathcal{L}$. Assume for all $m \in \mathbb{N}$ that $\pi_m(\lambda^{(n)})$ converges to $\pi_m(\lambda)$ vaguely in $\mathcal{L}_m$ as $n \to \infty$. Then $\lambda^{(n)}$ converges vaguely to $\lambda$.*
2. *Let $\alpha^{(n)}, n \in \mathbb{N}$, be a sequence in $\mathcal{A}$. Assume for all $m \in \mathbb{N}$ that $\pi_m(\alpha^{(n)})$ converges to $\pi_m(\alpha)$ vaguely in $\mathcal{A}_m$ as $n \to \infty$. Then $\alpha^{(n)}$ converges vaguely to $\alpha$.*

*Consequently, the mapping $\Pi^{-1}$ constructed in Lemma 5.10 is continuous.*

**Proof** (1) Let $\lambda^{(n)}, n \in \mathbb{N}$, be a sequence in $\mathcal{L}$, such that for all $m \in \mathbb{N}$ it holds that $\pi_m(\lambda^{(n)})$ converges to $\pi_m(\lambda)$ vaguely in $\mathcal{L}_m$ as $n \to \infty$. Abbreviate $\mathcal{M}_{\mathbb{N}_0} := \mathcal{M}_{\mathbb{N}_0}(\mathcal{S}) \setminus \{0\}$ and fix a continuous and compactly supported function $g \colon \mathcal{M}_{\mathbb{N}_0} \to \mathbb{R}$. Recall the identification of $\mathcal{L}_m$ and $\bar{\mathcal{L}}_m$ introduced before Lemma 5.10. Clearly, we have for any $m \in \mathbb{N}$ that

$$\left| \int g \, d\lambda^{(n)} - \int g \, d\lambda \right|$$
$$\leq \left| \int g \, d\lambda^{(n)} - \int g \, d\bar{\pi}_m(\lambda^{(n)}) \right| + \left| \int g \, d\bar{\pi}_m(\lambda^{(n)}) - \int g \, d\bar{\pi}_m(\lambda) \right|$$
$$+ \left| \int g \, d\bar{\pi}_m(\lambda) - \int g \, d\lambda \right|$$
$$\leq \left( |\lambda^{(n)}| + |\lambda| \right) \|g - g \circ \bar{\pi}_m\|_\infty + \left| \int g \, d\bar{\pi}_m(\lambda^{(n)}) - \int g \, d\bar{\pi}_m(\lambda) \right|. \quad (5.31)$$

Note that $\lambda^{(n)} \in \mathcal{L}$ implies that $|\lambda^{(n)}| \leq c_{\lambda^{(n)}}(\mathcal{S}) \leq \sup_{N \in \mathbb{N}} \mu_N(\mathcal{S}) \vee \mu(\mathcal{S}) = 1$ for all $n \in \mathbb{N}$, and in the same way $|\lambda| \leq 1$. Further, the support of the function $g$ is compact, and thus bounded. So by Lemma 5.11 we have that the mappings $\bar{\pi}_m$ restricted to the support of $g$ converge uniformly to the identity, as $m \to \infty$. Hence, we can first choose $m \in \mathbb{N}$ sufficently large such that $(|\lambda^{(n)}| + |\lambda|) \|g - g \circ \bar{\pi}_m\|_\infty$ is arbitrarily small, uniformly in $n$. Then we can use that, $\bar{\pi}_m(\lambda^{(n)}) \to \bar{\pi}_m(\lambda)$, as $n \to \infty$, holds by assumption, so the second summand on the right-hand side of (5.31) vanishes as $n \to \infty$.

(2) Let $\alpha^{(n)}, n \in \mathbb{N}$, be a sequence in $\mathcal{A}$, such that for all $m \in \mathbb{N}$ it holds that $\pi_m(\alpha^{(n)})$ converges to $\pi_m(\alpha)$ vaguely in $\mathcal{A}_m$ as $n \to \infty$. Recall that for all $n \in \mathbb{N}$ the measures $\alpha^{(n)}$ are concentrated on $\mathcal{M}_{\leq 1}(\mathcal{S})$, i.e., on subprobability measures on $\mathcal{S}$, since for any $y \in \mathrm{supp}(\alpha^{(n)})$, we have that $y(\mathcal{S}) \leq c_{\alpha^{(n)}}(\mathcal{S}) \leq \sup_{N \in \mathbb{N}} \mu_N(\mathcal{S}) \vee \mu(\mathcal{S}) = 1$. The same holds for $\alpha$. Hence, without loss of generality we can show vague convergence by considering any continuous compactly supported test function $g \colon \mathcal{M}_{\leq 1}(\mathcal{S}) \setminus \{0\} \to \mathbb{R}$. By Lemma 5.1 there exists $\varepsilon > 0$ such that the support of $g$ is contained in $N_\varepsilon := \{\nu \in \mathcal{M}_{\leq 1}(\mathcal{S}) \setminus \{0\} \colon \nu(\mathcal{S}) \geq \varepsilon\}$. Analogously to (5.31) we get

$$\left| \int g \, d\alpha^{(n)} - \int g \, d\alpha \right| \leq \left( \alpha^{(n)}(N_\varepsilon) + \alpha(N_\varepsilon) \right) \sup_{y \in N_\varepsilon} |g(y) - g \circ \bar{\pi}_m(y)|$$

$$+ \left| \int g \, \mathrm{d}\bar{\pi}_m(\alpha^{(n)}) - \int g \, \mathrm{d}\bar{\pi}_m(\alpha) \right|. \tag{5.32}$$

To bound the first summand observe the following: Since $\alpha^{(n)} \in \mathcal{A}$ one has that $\varepsilon \alpha^{(n)}(N_\varepsilon) \leq c_{\alpha^{(n)}}(\mathcal{S}) \leq 1$ and hence $\alpha^{(n)}(N_\varepsilon) \leq 1/\varepsilon$. The same holds for $\alpha$ and hence $\alpha^{(n)}(N_\varepsilon) + \alpha(N_\varepsilon) \leq 2/\varepsilon$. Also, by Lemma 5.11 the supremum in (5.32) vanishes, as $m \to \infty$. So, by first choosing $m$ large enough and then using that $\bar{\pi}_m(\alpha^{(n)}) \to \bar{\pi}_m(\alpha)$, as $n \to \infty$, the right-hand side of (5.32) vanishes. □

### 5.3 Dawson–Gärtner and identification of the rate function

As an easy consequence of Lemma 5.4 and the LDP of Theorem 3.1, we obtain upper and lower LDP bounds for $(\mathrm{Mi}_N^{(m)}, \mathrm{Ma}_N^{(m)})$ with different rate functions. In order to formulate this (in particular, to identify the rate functions) we need to introduce additional notation.

For dealing with the microscopic clusters, we need the discretized version of the connection parameter $\tau$ defined in (1.9), which now has to be understood with respect to the discretized kernels. For $* \in \{+, -\}$ and $k \in \mathcal{M}(\mathcal{S}_m) \setminus \{0\}$ we write

$$\tau_m^{(*)}(k) = \sum_{T \in \mathcal{T}(k)} \prod_{\{i,j\} \in E(T)} \kappa^{(m,*)}(r_i, r_j), \tag{5.33}$$

where $(r_i)_{i=1,\ldots,|k|} \in \mathcal{S}_m^{|k|}$ is such that $k = \sum_{i=1}^{|k|} \delta_{r_i} \in \mathcal{M}_{\mathbb{N}_0}(\mathcal{S}_m)$ and $\mathcal{T}(k)$ is the set of spanning trees on $[|k|]$ and we recall that $|k| = k(\mathcal{S}_m)$. Further, we define for $* \in \{+, -\}$

$$C_m^{(*)} = \frac{1}{2} \langle \mu^{(m)}, \kappa^{(m,*)} \mu^{(m)} \rangle, \tag{5.34}$$

$$\hat{I}_{\mathrm{Mi}}^{(m,*)}(\lambda) = \mathbb{H}(\lambda | \mathbb{Q}_{\mu^{(m)}}) - 1 - \langle \lambda, \log \tau_m^{(*)} \rangle + |c_\lambda| - |\lambda|, \quad \lambda \in \mathcal{L}_m, \tag{5.35}$$

$$\hat{I}_{\mathrm{Me}}^{(m,*)}(\nu) = \left\langle \nu, \log \frac{\mathrm{d}\nu}{(\kappa^{(m,*)}\nu) \, \mathrm{d}\mu^{(m)}} \right\rangle, \quad \nu \in \mathcal{M}(\mathcal{S}_m). \tag{5.36}$$

where we recall that $\mathbb{Q}_{\mu^{(m)}}$ is the distribution of a Poisson point process on $\mathcal{S}_m$ with intensity measure $\mu^{(m)}$ and that $\mathbb{H}$ denotes relative entropy between non-normalized measures, defined as in (1.8). Again, we adopt the convention that $\hat{I}_{\mathrm{Me}}^{(m,*)}(\nu) = \infty$ if $\frac{\mathrm{d}\nu}{(\kappa^{(m,*)}\nu) \, \mathrm{d}\mu^{(m)}}$ does not exist. We have to be more careful in the definition of the macroscopic rate. For $\alpha \in \mathcal{A}_m$ we define

$$\hat{I}_{\mathrm{Ma}}^{(m,*)}(\alpha)$$
$$= \begin{cases} \int_{\mathcal{M}(\mathcal{S}_m) \setminus \{0\}} \alpha(\mathrm{d}y) \left[ \left\langle y, \log \frac{\mathrm{d}y}{(1 - e^{-\kappa^{(m,*)}y}) \, \mathrm{d}\mu^{(m)}} \right\rangle - \frac{1}{2} \langle y, \kappa^{(m,*)} y \rangle \right] & \text{if } \alpha \text{ is connectable,} \\ \infty & \text{otherwise,} \end{cases}$$
$$\tag{5.37}$$

where again we define $I_{\text{Ma}}(\alpha) = \infty$, if it is not true that $\alpha$-almost everywhere the density in the log-term exists. The definition of $\hat{I}_{\text{Ma}}$ also takes into account the possibility that the discretized kernel might not be irreducible. In particular if $\kappa$ is irreducible, it is not always true that $\kappa^{(m,-)}$ is also irreducible. In that case we have to additionally assume that $\alpha$ is connectable with respect to $\kappa^{(m,*)}$ to get a finite rate, as it is formulated in the generalized version of Theorem 3.1, i.e., Theorem 3.2. We will comment on this in detail below.

When estimating the distribution of the pair $(\text{Mi}^{(m)}, \text{Ma}^{(m)})$ under the measure $\mathbb{P}_N$ by Lemma 5.4 we will get an additional error term. To deal with that we define for $(\lambda, \alpha) \in \mathcal{L}_m \times \mathcal{A}_m$

$$I^{(m,+)}(\lambda, \alpha) = \begin{cases} \hat{I}_{\text{Mi}}^{(m,+)}(\lambda) + \hat{I}_{\text{Ma}}^{(m,+)}(\alpha) + \hat{I}_{\text{Me}}^{(m,+)}(\mu^{(m)} - c_\lambda - c_\alpha) + C_m^{(-)} & \text{if } c_\lambda + c_\alpha \leq \mu^{(m)}, \\ \infty & \text{otherwise,} \end{cases} \tag{5.38}$$

$$I^{(m,-)}(\lambda, \alpha) = \begin{cases} \hat{I}_{\text{Mi}}^{(m,-)}(\lambda) + \hat{I}_{\text{Ma}}^{(m,-)}(\alpha) + \hat{I}_{\text{Me}}^{(m,-)}(\mu^{(m)} - c_\lambda - c_\alpha) + C_m^{(+)} & \text{if } c_\lambda + c_\alpha \leq \mu^{(m)}, \\ \infty & \text{otherwise.} \end{cases} \tag{5.39}$$

**Corollary 5.13** (LDP bounds for $(\text{Mi}_N^{(m)}, \text{Ma}_N^{(m)})$ under $\mathcal{G}_N$) *Assume that $\mu_N$ converges to $\mu$ as $N \to \infty$. Let $\kappa_N$ converge to a continuous kernel $\kappa$ that is irreducible w.r.t. $\mu$. Fix $m \in \mathbb{N}$ and let $\kappa_N^{(m,-)}$ and $\kappa_N^{(m,+)}$, $N \in \mathbb{N}$, be two sequences of kernels on $\mathcal{S}_m$ satisfying (5.12). Then the distribution of $(\text{Mi}_N^{(m)}, \text{Ma}_N^{(m)})$ under $\mathbb{P}_N$ satisfies, as $N \to \infty$, the upper large-deviations bound with rate function $I^{(m,+)}$ and the lower large-deviations bound with rate function $I^{(m,-)}$.*

*Proof* It is easy to verify that

$$\lim_{N \to \infty} \frac{1}{N} \log \Delta_N^{(m)} = \tfrac{1}{2} \langle \mu^{(m)}, \kappa^{(m,+)} \mu^{(m)} \rangle - \tfrac{1}{2} \langle \mu^{(m)}, \kappa^{(m,-)} \rangle = C_m^{(+)} - C_m^{(-)},$$

where $\Delta_N^{(m)}$ is defined as in (5.18). Note that irreducibility of $\kappa$ implies irreducibility of $\kappa^{(m,+)}$. Hence, for the upper large-deviations bound we can apply Theorem 3.1 after using the upper bound from Lemma 5.4. For the lower bound we use the lower bound from Lemma 5.4 and the lower bound from Theorem 3.2 that also applies if $\kappa^{(m,-)}$ is reducible. □

Of course, the basic idea is to use Corollary 5.13 for very large $m$ and the hope is that the discretized rate functions approximate the rate function $I$ given in Theorem 1.1. However, the problem is that the lower bound can be arbitrarily bad due to the following issue. Observe that although $\kappa$ is assumed to be irreducible with respect to $\mu$, as a consequence of the definition of the lower approximation $\kappa^{(m,-)}$ given in (5.13) we have to deal with the possibility that $\kappa^{(m,-)}$ is not irreducible with respect to $\mu \circ \pi_m^{-1}$. This might even be the case for all $m \in \mathbb{N}$. To illustrate this, we give a brief example.

**Example 5.14** Choose $\mathcal{S} = [0, 1]$, let $\mu$ be the Lebesgue measure and $\kappa_N(x, x') = \kappa(x, x') = xx'$ for $x, x' \in \mathcal{S}$ and all $N \in \mathbb{N}$. Let $\{P_m\}_{m \in \mathbb{N}}$ be any nested partition for $\mathcal{S}$ and let $\mathcal{S}_m, m \in \mathbb{N}$, be any choice for the sets of representative points. Then for any $m \in \mathbb{N}$ there exists some set $A_m \in P_m$ such that $0 \in \partial A_m$ and $\mu(A_m) > 0$. Let $x_m$ be the representative of $A_m$. It can be directly verified that for any $x'_m \in \mathcal{S}_m$ we have that $\kappa^{(m,-)}(x_m, x'_m) = 0$, although $\mu^{(m)}(\{x_m\}) = \mu(A_m) > 0$. On the other hand, assumption (5.4) implies that $\mu(A_m) \to 0$ as $m \to \infty$ so there exists $m_0$ such that for $m \geq m_0$ we also have that $\mu^{(m)}(\mathcal{S}_m \setminus \{x_m\}) > 0$ and hence $\kappa^{(m,-)}$ is reducible with respect to $\mu^{(m)}$ for all $m \geq m_0$. ◊

Regarding the approximation of the rate function $I$ by its discrete version the problem of (missing) irreducibility enters our analysis only via the function $\hat{I}_{\text{Ma}}^{(m,-)}$. Indeed, for $\hat{I}_{\text{Ma}}^{(m,+)}$ the additional case distinction given in (5.37) is not necessary, since $\kappa^{(m,+)}$ is always irreducible and hence every $\alpha$ is connectable with respect to $\kappa^{(m,+)}$. In order to make the lower approximation work, we formulate additional assumptions on $\alpha$. Since they are not satisfied for any $\alpha$ we have to find a way how to deal with the other cases.

**Lemma 5.15** (Identification of the rate function) *Let $\kappa : \mathcal{S} \times \mathcal{S} \to [0, \infty)$ be continuous and irreducible with respect to $\mu$. Assume that $\kappa^{(m,+)}, \kappa^{(m,-)} : \mathcal{S}_m \times \mathcal{S}_m \to [0, \infty)$ are given such that $\kappa^{(m,-)} \leq \kappa^{(m,+)}$ and*

$$\lim_{m \to \infty} \kappa^{(m,*)}(\pi_m(x), \pi_m(y)) = \kappa(x, y), \quad \text{uniformly in } x, y \in \mathcal{S}, * \in \{+, -\},$$

*monotonously decreasing for $* = +$ and monotonously increasing for $* = -$. We abbreviate $\Pi_m(\lambda, \alpha) := (\pi_m(\lambda), \pi_m(\alpha))$ for $m \in \mathbb{N}$, $(\lambda, \alpha) \in \mathcal{L} \times \mathcal{A}$. Then the rate function $I$ introduced in Theorem 1.1 satisfies the following.*

1. *For any $(\lambda, \alpha) \in \mathcal{L} \times \mathcal{A}$ it holds that $I^{(m,+)}(\Pi_m(\lambda, \alpha)) \nearrow I(\lambda, \alpha)$ as $m \to \infty$.*
2. *Let $(\lambda, \alpha) \in \mathcal{L} \times \mathcal{A}$ and assume that $\alpha$ satisfies the following: there exists some $m_0 \in \mathbb{N}$ such that for all $m \geq m_0$ the measures $\pi_m(\alpha)$ are connectable with respect to $\kappa^{(m,-)}$. Then $I^{(m,-)}(\Pi_m(\lambda, \alpha)) \to I(\lambda, \alpha)$ as $m \to \infty$.*

**Proof** Fix $\lambda \in \mathcal{L}$ and $\alpha \in \mathcal{A}$. Also, we will denote $\lambda_m = \pi_m(\lambda)$ and $\alpha_m = \pi_m(\alpha)$.

We first assume that $c_\lambda + c_\alpha \leq \mu$. It is straightforward to show that this implies that $c_{\lambda_m} + c_{\alpha_m} \leq \mu^{(m)}$ for all $m \in \mathbb{N}$. We denote $\nu_m = \mu^{(m)} - c_{\lambda_m} - c_{\alpha_m}$.

We will see in the proof that in each of the terms that we handle, the main part is an entropy between two image measures under $\pi_m$ plus a perturbation and other terms that will turn out to converge monotonically as $m \to \infty$. Then we will use [22, Prop. 15.6], which says that the named entropy converges, as $m \to \infty$, to its supremum over $m$.

For this, we will handle each of the four terms in (5.34)–(5.36) separately.

*Step 1: term $C_m^{(*)}$.* It is easy to deduce that $\lim_{m \to \infty} C_m^{(+)} = \frac{1}{2} \langle \mu, \kappa \mu \rangle = \lim_{m \to \infty} C_m^{(-)}$ and that $C_m^{(*)}$ is decreasing in $m$ for $* = +$ and increasing for $* = -$. Note that $I^{(m,+)}$ is defined with $C_m^{(-)}$, so it has the right direction of monotonicity.

*Step 2: term* $\hat{I}_{\mathrm{Mi}}^{(m,*)}(\lambda_m)$. Let us turn to the part $\hat{I}_{\mathrm{Mi}}^{(m,*)}(\lambda_m)$. Recall that $\lambda_m = \lambda \circ \pi_m^{-1}$.

Since $\lambda_m = \lambda \circ \pi_m^{-1}$ and $\mathbb{Q}_{\mu^{(m)}} = \mathbb{Q}_\mu \circ \pi_m^{-1}$ by the mapping theorem for Poisson point processes, [22, Prop. 15.6] implies that $\mathbb{H}(\lambda_m|\mathbb{Q}_{\mu^{(m)}})$ converges towards $\mathbb{H}(\lambda|\mathbb{Q}_\mu)$, and $\mathbb{H}(\lambda_m|\mathbb{Q}_{\mu^{(m)}})$ is increasing in $m$. According to our assumption in (5.16), using the monotone convergence theorem, we see that the term $-\langle \lambda_m, \log \tau_m^{(*)} \rangle$ converges towards $-\langle \lambda, \log \tau \rangle$ (also if $-\langle \lambda, \log \tau \rangle = \infty$). For $* = +$ the convergence is from below, as desired. It is easy to verify that $|c_{\lambda_m}| - |\lambda_m| = |c_\lambda| - |\lambda|$ holds for all $m$. Hence, we have shown that $\hat{I}_{\mathrm{Mi}}^{(m,+)}(\pi_m(\lambda)) \nearrow \hat{I}_{\mathrm{Mi}}(\lambda)$ and $\hat{I}_{\mathrm{Mi}}^{(m,-)}(\pi_m(\lambda)) \to \hat{I}_{\mathrm{Mi}}(\lambda)$ as $m \to \infty$.

*Step 3: term* $\hat{I}_{\mathrm{Me}}^{(m,*)}(\nu_m)$. Now we turn to the mesoscopic term $\hat{I}_{\mathrm{Me}}^{(m,*)}(\nu_m)$. Note that for $\nu = \mu - c_\lambda - c_\alpha$ the definition of the image measure implies that $\nu_m = \nu \circ \pi_m^{-1}$. Further, we have that

$$\hat{I}_{\mathrm{Me}}^{(m,*)}(\nu_m) = \left\langle \nu_m, \log \frac{\nu_m}{\mu^{(m)}} \right\rangle - \langle \nu_m, \log(\kappa^{(m,*)} \nu_m) \rangle.$$

Hence, we may apply [22, Prop. 15.6] to the first term and see that $\langle \nu_m, \log \frac{\mathrm{d}\nu_m}{\mathrm{d}\mu^{(m)}} \rangle$ converges as $m \to \infty$ to its supremum on $m$, and for the second term we use (5.16) to see that $\lim_{m \to \infty} -\langle \nu_m, \log(\kappa^{(m,*)} \nu_m) \rangle = -\langle \nu, \log(\kappa \nu) \rangle$ holds by the monotone convergence theorem (also if $-\langle \nu, \log(\kappa \nu) \rangle = \infty$). For $* = +$ the sequence is increasing as desired.

*Step 4: term* $\hat{I}_{\mathrm{Ma}}^{(m,*)}(\alpha_m)$. Finally, we turn to the macroscopic term $\hat{I}_{\mathrm{Ma}}^{(m,*)}(\alpha_m)$. Note that irreducibility of $\kappa$ with respect to $\mu$ implies irreducibility of $\kappa^{(m,+)}$ with respect to $\mu^{(m)}$. Hence the measures $\pi_m(\alpha), m \in \mathbb{N}$, are connectable. In the setting of statement (2) this is true for $m \geq m_0$ by assumption. Therefore, without loss of generality, we assume $m \geq m_0$. Applying the definition of the image measure we have that

$$\hat{I}_{\mathrm{Ma}}^{(m,*)}(\alpha_m) = \int \alpha(\mathrm{d}y) \left\langle y_m, \log \frac{\mathrm{d}y_m}{\mathrm{d}\mu^{(m)}} \right\rangle - \int \alpha(\mathrm{d}y) \Big[ \langle y_m, \log(1 - e^{-\kappa^{(m,*)} y_m}) \rangle + \frac{1}{2} \langle y_m, \kappa^{(m,*)} y_m \rangle \Big],$$

where we wrote $y_m = y \circ \pi_m^{-1}$. Now the convergence (including the desired monotonicity) can be argued in the same way as for the microscopic and the mesoscopic terms above, using the monotonicity of the entropy and our assumption in (5.16).

This finishes the proof of (1) and (2) in the case $c_\lambda + c_\alpha \leq \mu$.

Finally, we are considering the case that $\lambda \in \mathcal{L}$ and $\alpha \in \mathcal{A}$ do not satisfy $c_\lambda + c_\alpha \leq \mu$. Then there exists a $\mu$-continuity set $A \subset \mathcal{S}$ such that $c_\lambda(A) + c_\alpha(A) > \mu(A)$. Then with $A_m = \pi_m(A)$ we have that $\pi_m^{-1}(A_m) \supseteq A$ and hence $c_{\lambda_m}(A_m) + c_{\alpha_m}(A_m) \geq c_\lambda(A) + c_\alpha(A)$. Now, let $0 < \varepsilon < c_\lambda(A) + c_\alpha(A) - \mu(A)$. By Lemma 5.11 we have that $\bar{\pi}_m(\mu) \to \mu$ weakly, as $m \to \infty$. Hence, we can choose $m_0$ large enough such that $\mu(A) \geq \bar{\pi}_m(\mu)(A) - \varepsilon/2 = \pi_m(\mu)(A_m) - \varepsilon/2$ holds for all $m \geq m_0$. Consequently, $c_{\lambda_m}(A_m) + c_{\alpha_m}(A_m) \geq \pi_m(\mu)(A_m) + \varepsilon/2$. Therefore $I^{(m,*)}(\Pi_m(\lambda, \alpha)) = \infty$ for any $m \geq m_0$ and $* \in \{+, -\}$. □

Now we derive an upper and a lower bound LDP for the distribution of $(\mathrm{Mi}_N, \mathrm{Ma}_N)$ under $\mathcal{G}_N$ by following the two parts of the proof of the Dawson–Gärtner theorem, [19, Theorem 4.6.1] and using the fact that due to Proposition 5.6 it is sufficient to work with open and closed sets from the projective limit topology.

**Lemma 5.16** (LDP—upper bound) *Suppose that all assumptions of Theorem 1.1 are satisfied. Then the distribution of $(\mathrm{Mi}_N, \mathrm{Ma}_N)$ under $\mathbb{P}_N$ satisfies the upper bound part of the LDP with rate function $I$ as defined in Theorem 1.1.*

*Proof* Fix a set $F \subset \mathcal{L} \times \mathcal{A}$ that is closed with respect to the vague topology. Then by Proposition 5.6 the set $F$ is also closed with respect to the projective limit topology. For any $m \in \mathbb{N}$ we use the notation $\Pi_m(\lambda, \alpha) := (\pi_m(\lambda), \pi_m(\alpha))$ and recall that $(\mathrm{Mi}_N^{(m)}, \mathrm{Ma}_N^{(m)}) = \Pi_m(\mathrm{Mi}_N, \mathrm{Ma}_N)$. Therefore,

$$\limsup_{N \to \infty} \frac{1}{N} \log \mathbb{P}_N\big((\mathrm{Mi}_N, \mathrm{Ma}_N) \in F\big) \leq \limsup_{N \to \infty} \frac{1}{N} \log \mathbb{P}_N\big((\mathrm{Mi}_N^{(m)}, \mathrm{Ma}_N^{(m)}) \in \Pi_m(F)\big)$$
$$\leq - \inf_{\Pi_m(F)} I^{(m,+)} = - \inf_F I^{(m,+)} \circ \Pi_m,$$

where we used Corollary 5.13 for the second inequality. Since the left-hand side does not depend on $m$, we can proceed with the supremum over $m \in \mathbb{N}$ on the right-hand side. By Lemma 5.15 we have that $\sup_m I^{(m,+)} \circ \Pi_m = I$, which implies the claim. □

It remains to prove the lower bound of the LDP formulated in Theorem 1.1. Since the approximation of $I$ from below via $I^{(m,-)}$ works only in the case where $\alpha$ satisfies the assumptions given in the second statement of Lemma 5.15, we have to do some additional work. The idea is the following: given some $\alpha$ that does not fulfill the assumptions, we will first approximate $\alpha$ by some suitable choice for which the assumptions hold. Then we can apply Lemma 5.15.

**Lemma 5.17** *Let $\alpha \in \mathcal{A}$ with $c_\alpha \leq \mu$ and $I_{\mathrm{Ma}}(\alpha) < \infty$. Then there exists a sequence $(\alpha^{(\delta,\varepsilon)})_{\delta>0, \varepsilon>0}$ in $\mathcal{A}$ such that the following properties hold:*

*(1) for fixed $\delta > 0$ and $\varepsilon > 0$ there exists $m_0 = m_0(\delta)$ such that for all $m \geq m_0$ the measures $\pi_m(\alpha^{(\delta,\varepsilon)})$ are connectable with respect to $\kappa^{(m,-)}$;*
*(2) $\alpha^{(\delta,\varepsilon)} \to \alpha$ as $\delta \to 0$ and $\varepsilon \to 0$ with respect to the vague topology;*
*(3) for any $\lambda \in \mathcal{L}$ we have that $I(\lambda, \alpha^{(\delta,\varepsilon)}) \to I(\lambda, \alpha)$ as $\delta \to 0$ and $\varepsilon \to 0$.*

*Proof* We always write $\alpha = \sum_{i \in J} \delta_{y_i}$, where $J$ is a countable set. The idea is to pick some type $x \in \mathcal{S}$ and to restrict the measures $y_i \in \mathcal{M}_{\leq 1}(\mathcal{S}) \setminus \{0\}$, $i \in J$, to a subset $\mathcal{S}_\delta$ of $\mathcal{S}$ that contains all types $x' \in \mathcal{S}$ that can be connected to $x$ by using a finite sequence of intermediate types $x_{h-1}, x_h$ for which we have $\kappa(x_{h-1}, x_h) \geq \delta$. Then for large enough $m$ connectivity is preserved with respect to $\kappa^{(m,-)}$, which will imply (1). The parameter $\varepsilon > 0$ is only introduced to deal with the fact that $J$ might be infinite and ensures the convergence claimed in (3).

Fix some $x \in \mathrm{supp}(\mu)$. For $\delta > 0$ define

$$\mathcal{S}_\delta := \mathcal{S}_\delta(x) := \{x' \in \mathcal{S} \colon \exists k \in \mathbb{N}, \exists x_0, x_1, \ldots, x_k \in \mathcal{S}, \, x_0 = x, \, x_k = x', \, \kappa(x_{h-1}, x_h) \geq \delta, \, \forall h \in [k]\}. \tag{5.40}$$

We first show that $\mu(\mathcal{S}\setminus S_\delta) \to 0$ as $\delta \to 0$, which we will need for the proof of (2) and (3). Observe that $S_\delta \subset S_{\delta'}$ if $\delta > \delta'$ and put $S_0 := \bigcup_{\delta>0} S_\delta$. Since $\kappa$ is irreducible with respect to $\mu$ and $x \in \mathrm{supp}(\mu)$ it is easy to see that for $\delta$ small enough we have that $\mu(S_\delta) > 0$ and hence $\mu(S_0) > 0$. We now argue that $\kappa = 0$ $\mu$-almost everywhere on $S_0 \times \mathcal{S}\setminus S_0$. Assume the contrary, i.e., $\int_{S_0} \int_{\mathcal{S}\setminus S_0} \kappa(x', x'')\mu(\mathrm{d}x')\mu(\mathrm{d}x'') > 0$. Then by continuity of $\kappa$ we find sets of positive measure $A \subset S_0$ and $B \subset \mathcal{S}\setminus S_0$ where $\delta' := \inf_{x' \in A, x'' \in B} \kappa(x', x'') > 0$, which is a contradiction to the fact that $B \subset \mathcal{S}\setminus S_0$. By the irreducibility of $\kappa$ with respect to $\mu$, the facts that $\mu(S_0) > 0$ and $\kappa = 0$ holds $\mu$-almost everywhere on $S_0 \times \mathcal{S}\setminus S_0$ imply that $\mu(\mathcal{S}\setminus S_0) = 0$, so by continuity of measures we have that $\mu(\mathcal{S}\setminus S_\delta) \to 0$ as $\delta \to 0$.

Now for any $\delta \geq 0$ and $\varepsilon \geq 0$ we define

$$\alpha^{(\delta,\varepsilon)} := \sum_{i \in J_\varepsilon} \delta_{y_i^{(\delta)}}, \quad \text{where} \quad y^{(\delta)} := y(\cdot \cap S_\delta) \text{ for any } y \in \mathcal{M}_{\leq 1}(\mathcal{S}) \quad (5.41)$$

and where $J_\varepsilon = \{i \in J : |y_i| > \varepsilon\}$. Note that $|c_\alpha| \leq 1$ implies that for $\varepsilon > 0$ the set $J_\varepsilon$ is finite.

Now we show that $(\alpha^{(\delta,\varepsilon)})_{\delta,\varepsilon>0}$ has the three properties.

(1) Fix $\varepsilon > 0$ and $\delta > 0$. Now take $m_0$ such that $\|\kappa - \kappa^{(m,-)} \circ \pi_m\|_\infty \leq \delta/2$ holds for all $m \geq m_0$. Let $m \geq m_0$. Then we have that $\kappa^{(m,-)} \circ \pi_m$ is irreducible on $S_\delta$, since $\kappa(x_{i-1}, x_i) \geq \delta$ implies $\kappa^{(m,-)} \circ \pi_m(x_{i-1}, x_i) \geq \delta/2$. With other words, $\kappa^{(m,-)}$ is irreducible on $\pi_m^{-1}(S_\delta)$. For any $i \in J$ we have that $\mathrm{supp}(y_i^{(\delta)}) \subset S_\delta$ by construction and hence $\mathrm{supp}(\pi_m(y_i^{(\delta)})) \subset \pi_m^{-1}(S_\delta)$. Therefore, we get for all $m \geq m_0$ that $\pi_m(\alpha^{(\delta,\varepsilon)})$ is connectable with respect to $\kappa^{(m,-)}$.

(2) It is straightforward to show that for all $i \in J$ and $\delta > 0$ we have $d_{\mathrm{BL}}(y_i, y_i^{(\delta)}) \leq y_i(\mathcal{S}\setminus S_\delta) \leq \mu(\mathcal{S}\setminus S_\delta)$, where $d_{\mathrm{BL}}$ is the metric defined in (5.29) that induces the weak topology on $\mathcal{M}(\mathcal{S})$. For any continuous compactly supported test function $f \colon \mathcal{M}_{\leq 1}(\mathcal{S})\setminus\{0\} \to \mathbb{R}$ there exists $\varepsilon_f > 0$ such that $f = 0$ on $\{y \in \mathcal{M}_{\leq 1}(\mathcal{S})\setminus\{0\} \colon |y| \leq \varepsilon_f\}$, so for $\varepsilon \leq \varepsilon_f$ (including the case $\varepsilon = 0$) we have

$$\left|\int \alpha^{(\delta,\varepsilon)}(\mathrm{d}y)f(y) - \int \alpha(\mathrm{d}y)f(y)\right| = \sum_{i \in J_{\varepsilon_f}} |f(y_i^{(\delta)}) - f(y_i)|$$

$$\leq |J_{\varepsilon_f}| \max_{i \in J_{\varepsilon_f}} |f(y_i^{(\delta)}) - f(y_i)|. \quad (5.42)$$

Observe that the right-hand side converges to 0 as $\delta \to 0$. This implies (2).

(3) By lower-semicontinuity of the rate function it is clear that $I(\lambda, \alpha) \leq \lim_{\delta,\varepsilon \to 0} I(\lambda, \alpha^{(\delta,\varepsilon)})$. So we only have to deal with the other estimate, i.e., we need to find an upper bound for

$$I(\lambda, \alpha^{(\delta,\varepsilon)}) - I(\lambda, \alpha) = \left(I_{\mathrm{Me}}(\mu - c_\lambda - c_\alpha) - I_{\mathrm{Me}}(\mu - c_\lambda - c_{\alpha^{(\delta,\varepsilon)}})\right)$$
$$+ \left(I_{\mathrm{Ma}}(\alpha^{(\delta,\varepsilon)}) - I_{\mathrm{Ma}}(\alpha^{(0,\varepsilon)})\right) + \left(I_{\mathrm{Ma}}(\alpha^{(0,\varepsilon)}) - I_{\mathrm{Ma}}(\alpha)\right).$$
(5.43)

Let $\gamma > 0$. It is straightforward to verify that $c_{\alpha^{(\delta,\varepsilon)}} \to c_\alpha$ weakly as $\delta, \varepsilon \to 0$. Since $I_{\text{Me}}$ is continuous we can choose $\delta_0$ and $\varepsilon_0$ such that $|I_{\text{Me}}(\mu - c_\lambda - c_\alpha) - I_{\text{Me}}(\mu - c_\lambda - c_{\alpha^{(\delta,\varepsilon)}})| \leq \gamma/3$ for all $\delta \leq \delta_0, \varepsilon \leq \varepsilon_0$.

Recall that

$$I_{\text{Ma}}(\alpha) = \int \alpha(\mathrm{d}y) f_{\text{Ma}}(y), \quad \text{where} \quad f_{\text{Ma}}(y) = \langle y, \log \tfrac{\mathrm{d}y}{(1-e^{-\kappa y})\mathrm{d}\mu} \rangle - \tfrac{1}{2}\langle y, \kappa y \rangle, \tag{5.44}$$

and where we interpret the log term as equal to $+\infty$ if the density does not exist.

We have that

$$I_{\text{Ma}}(\alpha) = I_{\text{Ma}}(\alpha^{(0,\varepsilon)}) + I_{\text{Ma}}(\textstyle\sum_{i \notin J_\varepsilon} \delta_{y_i}) \geq I_{\text{Ma}}(\alpha^{(0,\varepsilon)}) + I_{\text{Ma}}(\textstyle\sum_{i \notin J_\varepsilon} y_i),$$

where the inequality follows from (7.1), which is proved in Lemma 7.1. Note that $|\sum_{i \notin J_\varepsilon} y_i| \to 0$ as $\varepsilon \to 0$ and $f_{\text{Ma}}(y) \to 0$ as $|y| \to 0$. So we can choose $\varepsilon \leq \varepsilon_0$ such that $I_{\text{Ma}}(\alpha^{(0,\varepsilon)}) - I_{\text{Ma}}(\alpha) \leq \gamma/3$.

As a last step we want to choose $\delta = \delta(\varepsilon) \leq \delta_0$ such that $|I_{\text{Ma}}(\alpha^{(\delta,\varepsilon)}) - I_{\text{Ma}}(\alpha^{(0,\varepsilon)})| \leq \gamma/3$, so it remains to show that

$$\lim_{\delta \to 0} I_{\text{Ma}}(\alpha^{(\delta,\varepsilon)}) = I_{\text{Ma}}(\alpha^{(0,\varepsilon)}). \tag{5.45}$$

Notice that the definition of $y^{(\delta)}$ given in (5.41) implies that $|y| - \mu(\mathcal{S} \setminus \mathcal{S}_\delta) \leq |y^{(\delta)}| \leq |y|$. Therefore, we can choose $\delta$ small enough such that $\mu(\mathcal{S} \setminus \mathcal{S}_\delta) \leq \tfrac{1}{2} \min_{i \in J_\varepsilon}(|y_i| - \varepsilon)$ to ensure that for all $i \in J_\varepsilon$ we have that $|y_i^{(\delta)}| > \varepsilon$. Now, we choose a function $\chi_\varepsilon \colon \mathcal{M}_{\leq 1}(\mathcal{S}) \setminus \{0\} \to \mathbb{R}$ that is equal to one on $\{y \in \mathcal{M}_{\leq 1}(\mathcal{S}) \setminus \{0\} \colon |y| > \varepsilon\}$, equal to zero on $\{y \in \mathcal{M}_{\leq 1}(\mathcal{S}) \setminus \{0\} \colon |y| \leq \varepsilon/2\}$ and continuous. In particular, the function $\chi_\varepsilon f_{\text{Ma}}$ is then compactly supported and we have that for all $\delta' \leq \delta$

$$I_{\text{Ma}}(\alpha^{(\delta',\varepsilon)}) = \int \alpha^{(\delta',0)}(\mathrm{d}y) \chi_\varepsilon(y) f_{\text{Ma}}(y).$$

Technically, we still have the problem that $f_{\text{Ma}}$ can take values in $\mathbb{R} \cup \{+\infty\}$. But our assumption $I_{\text{Ma}}(\alpha) < \infty$ implies that $f_{\text{Ma}}(y_i) < \infty$ for all $i \in J_\varepsilon$, so by continuity of $f_{\text{Ma}}$ and the finiteness of $J_\varepsilon$ we can tune $\delta' \leq \delta$ such that uniformly for all $i \in J_\varepsilon$ we have $f_{\text{Ma}}(y_i^{(\delta)}) \leq C$ for some constant $C$. We already showed in (2) that $\alpha^{(\delta',0)} \to \alpha$, so having established continuity and compactly supportedness of the function in the integral we get that

$$\lim_{\delta' \to 0} I_{\text{Ma}}(\alpha^{(\delta',\varepsilon)}) = \lim_{\delta' \to 0} \int \alpha^{(\delta',0)}(\mathrm{d}y) \chi_\varepsilon(y)(f_{\text{Ma}}(y) \wedge C)$$
$$= \int \alpha(\mathrm{d}y) \chi_\varepsilon(y)(f_{\text{Ma}}(y) \wedge C) = I_{\text{Ma}}(\alpha^{(0,\varepsilon)}).$$

Altogether the right-hand side of (5.43) can be bounded by $\gamma$, which proves the claim. □

**Lemma 5.18** (LDP—lower bound) *Suppose that all the assumptions of Theorem 1.1 are satisfied. Then the distribution of $(Mi_N, Ma_N)$ under $\mathbb{P}_N$ satisfies the lower bound part of the LDP with rate function $I$.*

*Proof* Fix a set $G \in \mathcal{L} \times \mathcal{A}$ that is open with respect to the vague topology and a point $(\lambda, \alpha) \in G$. We will show that for any $\gamma > 0$ we have that

$$\liminf_{N \to \infty} \frac{1}{N} \log \mathbb{P}_N\big((Mi_N, Ma_N) \in G\big) \geq -I(\lambda, \alpha) - \gamma. \tag{5.46}$$

Note that in the case where $I(\lambda, \alpha) = \infty$, the claimed estimate always holds, so we assume that $I(\lambda, \alpha) < \infty$. Let $\gamma > 0$. We will use the approximating sequence $(\alpha^{(\delta,\varepsilon)})_{\delta>0,\varepsilon>0}$ constructed in Lemma 5.17. Let $\varepsilon > 0$ and $\delta > 0$ be small enough such that $\alpha^{(\delta,\varepsilon)} \in G$ and $|I(\lambda, \alpha) - I(\lambda, \alpha^{(\delta,\varepsilon)})| \leq \gamma/2$. Again, we write $\Pi_m(\tilde{\lambda}, \tilde{\alpha}) = (\pi_m(\tilde{\lambda}), \pi_m(\tilde{\alpha}))$ and also $\Pi_{m,n}(\tilde{\lambda}, \tilde{\alpha}) = (\pi_{m,n}(\tilde{\lambda}), \pi_{m,n}(\tilde{\alpha}))$ for any $n \geq m$. y Proposition 5.6 the set $G$ is also open with respect to the projective limit topology. It is a general fact that the set

$$\mathcal{B}_1 := \{\Pi_m^{-1}(U_m) : m \in \mathbb{N}, \ U_m \subset \mathcal{L}_m \times \mathcal{A}_m \text{ open}\} \tag{5.47}$$

is a basis of the projective limit topology. From now on, we fix some large $m_0 \in \mathbb{N}$ which we will specify later. We claim that also the set

$$\mathcal{B}_{m_0} := \{\Pi_m^{-1}(U_m) : m \geq m_0, \ U_m \subset \mathcal{L}_m \times \mathcal{A}_m \text{ open}\} \tag{5.48}$$

is a basis of the projective limit topology and argue this as follows: note that for any $m < m_0$ we can take any $n \geq m_0$ and use that $\pi_m = \pi_{m,n} \circ \pi_n$ holds by Lemma 5.7 to derive that $\Pi_m^{-1}(U_m) = \Pi_n^{-1}(\Pi_{m,n}^{-1}(U_m))$ for any open set $U_m$. Since the set $\Pi_{m,n}^{-1}(U_m)$ is again open due to the continuity of $\pi_{m,n}$, we get that $\Pi_m^{-1}(U_m) \in \mathcal{B}_{m_0}$. Altogether we have that $\mathcal{B}_1 \subset \mathcal{B}_{m_0}$.

Having established that $\mathcal{B}_{m_0}$ is a basis for the projective limit topology, we may pick $m \geq m_0$ and an open set $U_m \subset \mathcal{L}_m \times \mathcal{A}_m$ such that $(\lambda, \alpha^{(\delta,\varepsilon)}) \in \Pi_m^{-1}(U_m) \subset G$. Therefore, we see that

$$\liminf_{N \to \infty} \frac{1}{N} \log \mathbb{P}_N\big((Mi_N, Ma_N) \in G\big) \geq \liminf_{N \to \infty} \frac{1}{N} \log \mathbb{P}_N\big((Mi_N^{(m)}, Ma_N^{(m)}) \in U_m\big)$$
$$\geq -\inf_{U_m} I^{(m,-)} \geq -I^{(m,-)}(\Pi_m(\lambda, \alpha^{(\delta,\varepsilon)})), \tag{5.49}$$

where we used Corollary 5.13 for the second inequality. Using Lemmas 5.15 and 5.17 we can pick $m_0$ large enough such that for all $m \geq m_0$ the measures $\pi_m(\alpha^{(\delta,\varepsilon)})$ are connectable and $|I^{(m,-)}(\Pi_m(\lambda, \alpha^{(\delta,\varepsilon)})) - I(\lambda, \alpha^{(\delta,\varepsilon)})| \leq \gamma/2$. Altogether, this gives (5.46). □

## 6 The minimizers of $I_{\text{Mi}}$

In this section, we derive an explicit description of the minimizer(s) $\lambda$ of $I_{\text{Mi}}$ under the constraint $c_\lambda = c$ for any $c \in \mathcal{M}(\mathcal{S})$ satisfying $c \leq \mu$. This will allow us to solve the optimization problem in (2.6), i.e., to identify the minimizer(s) of the rate function for the LDP for $\text{Mi}_N$ in Theorem 2.3. It will also be used as an important intermediate step in deriving the full optimization of the rate function $I$ of the LDP in Theorem 2.1, our main result. Recall the notation that we introduced in Sect. 2.1, in particular the definition of $\Sigma(\kappa, c)$ from (2.3). Here is the main result of this section.

**Proposition 6.1** (Minimizers of $I_{Mi}$) *Fix a probability measure $\mu$ on $\mathcal{S}$ and a kernel $\kappa$ on $\mathcal{S} \times \mathcal{S}$ that is nonnegative and continuous.*

*Let $c \in \mathcal{M}(\mathcal{S})$ be a measure such that $c \leq \mu$.*

*(i) Assume that $\Sigma(\kappa, c) \leq 1$. Then*

$$\inf_{\lambda \in \mathcal{L}: \, c_\lambda = c} I_{Mi}(\lambda) = \left\langle c, \log \frac{\mathrm{d}c}{\mathrm{d}\mu} \right\rangle + \frac{1}{2} \langle c, \kappa(\mu - c) \rangle, \tag{6.1}$$

*and the infimum is attained in the unique minimizer $\lambda_c$ defined in (2.1).*

*(ii) Assume that $\Sigma(\kappa, c) > 1$. Then*

$$\inf_{\lambda \in \mathcal{L}: \, c_\lambda = c} I_{Mi}(\lambda) \geq \inf_{\lambda \in \mathcal{L}: \, c_\lambda = b^*} I_{Mi}(\lambda) + I_{Me}(c - b^*) \tag{6.2}$$

*where $b^* = b^*(c) \in \mathcal{M}(\mathcal{S})$ is the minimal, non-trivial (i.e., not equal to c) solution to (2.9) and satisfies $\Sigma(\kappa, b^*) = 1$.*

It is interesting to notice that one can see the phase transition already from the sole consideration of $I_{\text{Mi}}$. We will refer to (i) and (ii) as to the sub- and supercritical cases, respectively.

In the case where $\mathcal{S}$ is finite we can actually prove an equality in (6.2). In the general case we also expect this to be true, but did not attempt a proof, since the inequality will be enough to prove our main results, Theorems 2.3 and 2.1.

The proof is naturally divided into Sects. 6.1–6.4 according to the distinctions between finite $\mathcal{S}$ (the discrete case) or general compact $\mathcal{S}$ and between the sub- and supercritical cases. In Sect. 6.1 we construct minimizers for subcritical measures $c$ for finite $\mathcal{S}$; we analyze if the only candidate $\lambda$ (coming from the Euler–Lagrange equations) satisfies the constraint $c_\lambda = c$, which requires the result about the multivariate power series from Sect. 4.1. In Sect. 6.2 we generalize the results to a general compact type space via an approximation argument. In order to deal with supercritical measures $c$ for $\mathcal{S}$ a finite set, we also rely on combinatorial results in Sect. 6.3. Afterwards, we handle the general supercritical case in Sect. 6.4.

### 6.1 The discrete, subcritical case

In this section, we formulate and prove the main assertions about the minimizers of $I_{\text{Mi}}$ in the discrete case, i.e., the case of a finite type space $\mathcal{S}$. Here we will be using the notation of linear algebra, i.e., measures $\lambda \in \mathcal{L}$ on $\mathcal{M}_{\mathbb{N}_0}(\mathcal{S})$ will be written as sequences $(\lambda_k)_{k \in \mathbb{N}_0^{\mathcal{S}}}$.

Recall the definition of the rate function $I_{\mathrm{Mi}}$ from Theorem 3.1 as well as the notation for the integrated type-configuration $c_r(\lambda) = \sum_{k \in \mathbb{N}_0^{\mathcal{S}}} \lambda_k k_r$ for $r \in \mathcal{S}$ introduced in (3.5). We write $[0, \mu]$ for the set of all $c \in [0, \infty)^{\mathcal{S}}$ satisfying $0 \leq c_r \leq \mu_r$ for any $r \in \mathcal{S}$.

For $c \in [0, \infty)^{\mathcal{S}}$, define $\lambda(c) = (\lambda_k(c))_{k \in \mathbb{N}_0^{\mathcal{S}}}$ by

$$\lambda_k(c) = \tau(k) \prod_{s \in \mathcal{S}} \frac{(c_s e^{-(\kappa c)_s})^{k_s}}{k_s!}, \qquad k \in \mathbb{N}_0^{\mathcal{S}} \tag{6.3}$$

and note that this definition is the discrete analog of the general form of the minimizer in (2.1).

The aim of the present Sect. 6.1 is to verify the subcritical case of Proposition 6.1 in the discrete setting, which we restate here quickly.

**Proposition 6.2** *Let $c = (c_s)_{s \in \mathcal{S}}$ be in $[0, \mu]$. Assume that $\Sigma(\kappa, c) \leq 1$. Then*

$$\inf_{\lambda \in \mathcal{L}: \, c(\lambda) = c} I_{\mathrm{Mi}}(\lambda) = \left\langle c, \log \frac{c}{\mu} \right\rangle + \frac{1}{2} \langle c, \kappa(\mu - c) \rangle, \tag{6.4}$$

*and the infimum is attained in the unique minimizer $\lambda(c)$ defined in* (6.3).

To derive the form of the minimizer given in (6.3), we will start by giving a short heuristic. First note that $I_{\mathrm{Mi}}$ is a strictly convex function and that $\{\lambda : c(\lambda) = c\}$ is a convex set, which implies that there is at most one minimizer. Assume that a minimizer $\lambda^*$ exists in the interior of $\{\lambda : c(\lambda) = c\}$. Then by formally writing down the Euler–Lagrange equations, one can see that

$$\lambda_k^* = \tau(k) \prod_{s \in \mathcal{S}} \frac{\theta_s^{k_s}}{k_s!}, \qquad k \in \mathbb{N}_0^{\mathcal{S}}, \tag{6.5}$$

where $\theta = (\theta_s)_{s \in \mathcal{S}}$ is some non-negative real-valued vector. Note that $\theta$ has to be chosen in such a way that $c(\lambda^*) = c$, i.e., for every $r \in \mathcal{S}$ the multivariate power series

$$c_r(\lambda^*) = \sum_{k \in \mathbb{N}_0^{\mathcal{S}}} \tau(k) k_r \prod_s \frac{\theta_s^{k_s}}{k_s!} \tag{6.6}$$

converges with limit $c_r$. We already encountered in Sect. 4.1 that for $\theta = c e^{-\kappa c}$ the power series on the right-hand side of (6.6) has the right value. The following is just a reformulation of the results from Lemma 4.1 and Proposition 4.2 using the notation of the present section.

**Corollary 6.3** *Let $c = (c_s)_{s \in \mathcal{S}}$ in $[0, \mu]$ and assume that $\Sigma(\kappa, c) \leq 1$. Then for $\lambda^* = \lambda(c)$ we have that $c(\lambda^*) = c$.*

A rigorous argument showing that this choice uniquely minimizes $I_{\text{Mi}}$ can be found at the end of this section. The identification of the optimal rate $I_{\text{Mi}}(\lambda(c))$ needs an additional property of the minimizer, namely a formula of its total mass, which we derive in Lemma 6.5. For this we use the following recursive formula.

**Lemma 6.4** *Let $k \in \mathbb{N}_0^{\mathcal{S}}$. Then we have the recursion*

$$\sum_{r,s \in \mathcal{S}} \kappa(r,s) \sum_{m, \widetilde{m} \in \mathbb{N}_0^{\mathcal{S}}: m + \widetilde{m} = k} \Big( \prod_{u \in \mathcal{S}} \frac{k_u!}{m_u! \widetilde{m}_u!} \Big) \tau(m) m_r \tau(\widetilde{m}) \widetilde{m}_s = 2(|k| - 1) \tau(k). \tag{6.7}$$

*Proof* Let $(x_i)_{i \in [|k|]} \in \mathcal{S}^{|k|}$ be a vector compatible to $k$, i.e., $k = \sum_i \delta_{x_i}$ and recall the definition of $\tau(k)$ from (3.6). For $i, j \in [|k|]$ with $i \neq j$ define

$$W_{i,j} := \sum_{T \in \mathcal{T}(k): \{i,j\} \in E(T)} \prod_{\{v,w\} \in E(T)} \kappa(x_v, x_w), \tag{6.8}$$

i.e., $W_{i,j}$ is the total weight of trees containing the edge $\{i, j\}$. Observe, that each tree $T$ on $[|k|]$ contains exactly $|k| - 1$ edges and, for each edge $\{i, j\} \in E(T)$ the weight of $T$ appears once in $W_{i,j}$ and once in $W_{j,i}$. Thus, the weight of $T$ is counted $2(|k| - 1)$ times in the sum $\sum_{i \neq j} W_{i,j}$, which implies that

$$2(|k| - 1) \tau(k) = \sum_{i \neq j} W_{i,j}. \tag{6.9}$$

Now, for a fixed pair of types $r, s \in \mathcal{S}$ consider the weights of trees containing an edge connecting some type $r$ with some type $s$ vertex, i.e., consider $\sum_{i \neq j: x_i = r, x_j = s} W_{i,j}$. Notice that each tree contributing to this weight can be decomposed into an edge of weight $\kappa(r, s)$ and two trees $T_r$ and $T_s$ with roots of type $r$ and $s$ respectively. (The term 'root' is here only used to mark a certain vertex, not to give some directed structure.) This implies the formula

$$\sum_{i \neq j: x_i = r, x_j = s} W_{i,j} = \kappa(r, s) \sum_{m + \widetilde{m} = k} \Big( \prod_{u \in \mathcal{S}} \frac{k_u!}{m_u! \widetilde{m}_u!} \Big) \tau(m) m_r \tau(\widetilde{m}) \widetilde{m}_s. \tag{6.10}$$

Here we used formula (4.5) to collect the weight coming from the possible choices of $T_r$ and $T_s$, which is $\tau(m) m_r$ and $\tau(\widetilde{m}) \widetilde{m}_s$ respectively. Formula (6.7) now follows by summing over all possible pairs $r, s \in \mathcal{S}$. □

**Lemma 6.5** *Let $c = (c_s)_{s \in \mathcal{S}}$ be non-negative with $\Sigma(\kappa, c) \leq 1$. Then for $\lambda(c)$ defined as in* (6.3) *we have that*

$$|\lambda(c)| = \sum_{k \in \mathbb{N}_0^{\mathcal{S}}} \lambda_k(c) = \sum_{r \in \mathcal{S}} c_r - \frac{1}{2} \langle c, \kappa c \rangle. \tag{6.11}$$

*Proof* Writing $\lambda^* = \lambda(c)$ and using that $c(\lambda^*) = c$ we show the equivalent equation

$$\sum_{k \in \mathbb{N}_0^{\mathcal{S}}} \lambda_k^* (|k| - 1) = \frac{1}{2} \langle c, \kappa c \rangle. \tag{6.12}$$

For fixed $k \in \mathbb{N}_0^{\mathcal{S}}$ the recursive equation (6.7) for $\tau(k)$ easily implies

$$\frac{1}{2} \sum_{m, \widetilde{m} \in \mathbb{N}_0^{\mathcal{S}} : m + \widetilde{m} = k} \sum_{r, s \in \mathcal{S}} \lambda_m^* m_r \kappa(r, s) \lambda_{\widetilde{m}}^* \widetilde{m}_s = \lambda_k^* (|k| - 1).$$

With the assumption $\Sigma(\kappa, c) \leq 1$ all series in the next equations converge (absolutely) by Corollary 6.3, so by rearranging terms we get that

$$\sum_{k \in \mathbb{N}_0^{\mathcal{S}}} \lambda_k^* (|k| - 1) = \frac{1}{2} \sum_{m \in \mathbb{N}_0^{\mathcal{S}}} \sum_{\widetilde{m} \in \mathbb{N}_0^{\mathcal{S}}} \sum_{r, s \in \mathcal{S}} \lambda_m^* m_r \kappa(r, s) \lambda_{\widetilde{m}}^* \widetilde{m}_s = \frac{1}{2} \langle c, \kappa c \rangle.$$

□

Combining the results from Corollary 6.3 and Lemma 6.5 we can now give the proof of Proposition 6.2:

*Proof of Proposition 6.2* Assume that $\Sigma(\kappa, c) \leq 1$. Define $\lambda^* = \lambda(c)$ as in (6.3). By Corollary 6.3 we have that $c_r(\lambda^*) = \sum_k \lambda_k^* k_s = c_s$ for all $s \in \mathcal{S}$. Now, take any $\lambda$ satisfying $c_s(\lambda) = \sum_k \lambda_k k_s = c_s$ for all $s \in \mathcal{S}$. Then from (3.7), using the formula from Lemma 6.5, we get

$$\begin{aligned} I_{\mathrm{Mi}}(\lambda) &= \sum_k \lambda_k \log \frac{\lambda_k}{\lambda_k^*} + \sum_k \lambda_k \log \left( \prod_{s \in \mathcal{S}} \left( \frac{c_s e^{-(\kappa c)_s}}{\mu_s} \right)^{k_s} \right) \\ &\quad + \sum_k \lambda_k (|k| - 1) + \frac{1}{2} \langle c, \kappa \mu \rangle \\ &= \mathbb{H}(\lambda | \lambda^*) - |\lambda^*| + \left\langle c, \log \frac{c}{\mu} \right\rangle + |c| + \frac{1}{2} \langle c, \kappa \mu \rangle - \langle c, \kappa c \rangle \\ &\geq \left\langle c, \log \frac{c}{\mu} \right\rangle + \frac{1}{2} \langle c, \kappa (\mu - c) \rangle, \end{aligned} \tag{6.13}$$

where we wrote $\mathbb{H}(\lambda | \lambda^*) = \langle \lambda, \log \frac{\lambda}{\lambda^*} \rangle + |\lambda^*| - |\lambda|$ for the entropy and used that $\mathbb{H}(\lambda | \lambda^*) \geq 0$ with equality if and only if $\lambda = \lambda^*$. □

### 6.2 The general subcritical case

In this section we derive Proposition 6.1(i). The proof is similar to the proof of the discrete variant in Sect. 6.1. Again, there is an explicit candidate for the minimizer, but one has to prove that it is admissible, and we need to identify its total mass. This is done in the analogs of Corollary 6.3 and Lemma 6.5, see Lemmas 6.7 and 6.8, whose proofs proceed via a discrete approximation based on the material of Sect. 5.1.

In the current case of a general compact metric type space $\mathcal{S}$, the candidate for a minimizer is given in terms of a Poisson point process, see (6.14). Recall that we write $\mathbb{Q}_\theta$ for the distribution of a Poisson point process $\mathbb{X} = (X_i)_{i \in I}$ in $\mathcal{S}$ with intensity measure $\theta \in \mathcal{M}(\mathcal{S})$. We write $k = \sum_i \delta_{X_i} \in \mathcal{M}_{\mathbb{N}_0}(\mathcal{S})$ for the measure induced by the random point cloud. Note that the points $X_i$ do not have to be distinct with positive probability, if $\theta$ has no Lebesgue density. We start by noting a simple fact about the densities between absolute continuous Poisson point processes.

**Lemma 6.6** *Let $\theta, \hat{\theta} \in \mathcal{M}(\mathcal{S})$ with $\hat{\theta} \ll \theta$. Then $\mathbb{Q}_{\hat{\theta}} \ll \mathbb{Q}_\theta$ and*

$$\frac{\mathrm{d}\mathbb{Q}_{\hat{\theta}}}{\mathrm{d}\mathbb{Q}_\theta}(k) = \mathrm{e}^{\langle k, \log \frac{\mathrm{d}\hat{\theta}}{\mathrm{d}\theta}\rangle + \theta(\mathcal{S}) - \hat{\theta}(\mathcal{S})}, \quad k \in \mathcal{M}_{\mathbb{N}_0}(\mathcal{S}).$$

Recall the definition of $\tau(k)$ introduced in (1.9). Also recall that for a fixed $c \in \mathcal{M}$ according to definition (2.1) the candidate for the minimizer of $I_{\mathrm{Mi}}$ under the constraint $c_\lambda = c$ has the form

$$\lambda_c(\mathrm{d}k) = \mathrm{e}^{\theta_c(\mathcal{S})}\tau(k)\mathbb{Q}_{\theta_c}(\mathrm{d}k), \quad \text{where } \theta_c(\mathrm{d}r) = \mathrm{e}^{-\kappa c(r)}c(\mathrm{d}r). \tag{6.14}$$

We first provide a generalized version of Corollary 6.3 and Lemma 6.5. We first impose the stricter condition $\Sigma(\kappa, c) < 1$.

**Lemma 6.7** *Let $c \in \mathcal{M}(\mathcal{S})$ with $c \leq \mu$. Assume that $\Sigma(\kappa, c) < 1$. Then the following holds.*

1. *For any continuous test function $f : \mathcal{S} \to [0, \infty)$ we have that*

$$\int_{\mathcal{M}_{\mathbb{N}_0}(\mathcal{S})} \lambda_c(\mathrm{d}k)\langle k, f\rangle = \langle c, f\rangle. \tag{6.15}$$

2. *The total mass of $\lambda_c$ is given by*

$$|\lambda_c| = \int_{\mathcal{M}_{\mathbb{N}_0}(\mathcal{S})} \lambda_c(\mathrm{d}k) = c(\mathcal{S}) - \frac{1}{2}\langle c, \kappa c\rangle. \tag{6.16}$$

*Proof* We focus on showing Eq. (6.15); the proof of (6.16) is similar (see the end of the proof). Abbreviating $\theta := \theta_c$ and inserting the definition of $\lambda_c$ we have to prove that

$$\mathbb{Q}_\theta\Big[\tau(k)\langle f, k\rangle \mathrm{e}^{\theta(\mathcal{S})}\Big] = \langle f, c\rangle, \tag{6.17}$$

where we conceive $k$ as an $\mathcal{M}_{\mathbb{N}_0}(\mathcal{S})$-valued random variable on the left-hand side. The idea is to deduce the equality from the one that we have in the finite-type case by using the discretization scheme from Sect. 5.1. Recall the notation from Sect. 5.1, where we discretized the compact metric space $\mathcal{S}$ into finite spaces $\mathcal{S}_m$, $m \in \mathbb{N}$, and defined the projections $\pi_m$, $m \in \mathbb{N}$, on different spaces in equations (5.6)–(5.8). For $k \in \mathcal{M}_{\mathbb{N}_0}(\mathcal{S})$ we will again identify the discretized measure $\pi_m(k)$ with an element of $\mathbb{N}_0^{\mathcal{S}_m} \setminus \{0\}$. Via $\mathcal{S}_m \subset \mathcal{S}$ the function $f$ can be restricted to $\mathcal{S}_m$ and write $f_m = f|_{\mathcal{S}_m}$. Also, we write $c_m := \pi_m(c)$ and identify it with a vector $(c_m(r))_{r \in \mathcal{S}_m}$. Recall the definitions of the discretized kernels. Let $\kappa_m \in \{\kappa^{(m,+)}, \kappa^{(m,-)}\}$, where $\kappa^{(m,\star)}$ for $\star = \pm$ is defined as in (5.13) and (5.14). Denote $\theta_m(r) := e^{-(\kappa_m c_m)(r)} c_m(r)$, $r \in \mathcal{S}_m$. For $k \in \mathbb{N}_0^{\mathcal{S}_m}$ let $\tau_m(k)$ be defined as in (3.6), but with respect to $\kappa_m$. Fix a continuous function $f : \mathcal{S} \to \mathbb{R}$. Our aim is to show that

$$\mathbb{Q}_\theta\left[\tau(k)\langle f, k\rangle e^{\theta(\mathcal{S})}\right] = \lim_{m \to \infty} \mathbb{Q}_{\theta_m}\left[\tau_m(k)\langle f_m, k\rangle e^{\theta_m(\mathcal{S}_m)}\right] = \langle f, c\rangle, \quad (6.18)$$

which finishes the proof of (6.15).

We start with proving the second equality of (6.18). It is straightforward to show that

$$\left|\Sigma(\kappa_m, c_m) - \Sigma(\kappa, c)\right| \leq \|\kappa_m \circ \pi_m - \kappa\|_\infty. \quad (6.19)$$

and hence $\Sigma(\kappa_m, c_m) \to \Sigma(\kappa, c)$, as $m \to \infty$ and so we will have for large $m \in \mathbb{N}$ that $\Sigma(\kappa_m, c_m) < 1$. Then we get

$$\mathbb{Q}_{\theta_m}\left[\tau_m(k)\langle f_m, k\rangle e^{\theta_m(\mathcal{S}_m)}\right] = \sum_{r \in \mathcal{S}_m} f_m(r) \sum_{k \in \mathbb{N}_0^{\mathcal{S}_m}} \tau_m(k) k_r \prod_{s \in \mathcal{S}_m} \frac{(\theta_m(s))^{k_s}}{k_s!}$$

$$= \sum_{r \in \mathcal{S}_m} f_m(r) c_m(r) = \int_\mathcal{S} f(\pi_m(x)) c(dx)$$

$$\to \int_\mathcal{S} f(x) c(dx) = \langle f, c\rangle, \quad m \to \infty.$$

where the second equation only holds if $m$ is large enough, and thus $\Sigma(\kappa_m, c_m) \leq 1$ due to Lemma 4.1(i) and Proposition 4.2.

Now we show the first equation of (6.18). Note that $\mathbb{Q}_{\theta_m}$ is a point process on $\mathcal{S}_m$, whereas $\mathbb{Q}_\theta$ is a point process on $\mathcal{S}$. However, by defining an intensity measure on $\mathcal{S}$ by $\bar{\theta}_m(dx) := e^{-(\kappa_m c_m)(\pi_m(x))} c(dx)$ we have that $\theta_m = \bar{\theta}_m \circ \pi_m^{-1}$ and hence $\mathbb{Q}_{\theta_m} = \mathbb{Q}_{\bar{\theta}_m} \circ \pi_m^{-1}$ holds by the mapping theorem for Poisson point processes. Therefore, according to Lemma 6.6,

$$\mathbb{Q}_{\theta_m}\left[\tau_m(k)\langle f_m, k\rangle e^{\theta_m(\mathcal{S}_m)}\right] = \mathbb{Q}_{\bar{\theta}_m}\left[\tau_m(\pi_m(k))\langle f, \pi_m(k)\rangle e^{\theta_m(\mathcal{S}_m)}\right] = \mathbb{Q}_\theta\left[\Psi_m^f\right], \quad (6.20)$$

with

$$\Psi_m^f(k) := \tau_m(\pi_m(k))\langle f, \pi_m(k)\rangle e^{\langle k, \kappa c - (\kappa_m c_m)\circ \pi_m\rangle} e^{\theta(\mathcal{S})}, \quad k \in \mathcal{M}_{\mathbb{N}_0}(\mathcal{S}). \quad (6.21)$$

Note that $\Psi_m^f$ converges pointwise to $\Psi^f$, where $\Psi^f(k) = \tau(k)\langle f, k\rangle e^{\theta(\mathcal{S})}$, $k \in \mathcal{M}_{\mathbb{N}_0}(\mathcal{S})$. Hence, the first equation of (6.18) immediately follows as soon as we have given an argument for interchanging the limit as $m \to \infty$ and the integration with respect to $\mathbb{Q}_\theta$. We will be using Lebesgue's theorem about dominated convergence for that. Let us introduce a majorant. Recalling the definition (5.33) of $\tau_m^{(+)}$ we define $\widehat{\Psi}_m$ by

$$\widehat{\Psi}_m(k) := \tau_m^{(+)}(\pi_m(k))|\pi_m(k)|e^{\langle k, \kappa c - (\kappa^{(m,-)} c_m)\circ \pi_m\rangle} e^{\theta(\mathcal{S})}, \quad k \in \mathcal{M}_{\mathbb{N}_0}(\mathcal{S}). \quad (6.22)$$

Then, since $\kappa^{(m,-)} \leq \kappa_m \leq \kappa^{(m,+)}$ we clearly have for any $k \in \mathcal{M}_{\mathbb{N}_0}(\mathcal{S})$ that $\Psi_m^f(k) \leq \|f\|_\infty \widehat{\Psi}_m(k)$ and $\widehat{\Psi}_m(k) \leq \widehat{\Psi}_{m_0}(k)$ if $m \geq m_0$. Hence, $\widehat{\Psi}_{m_0}$ is a majorant. It remains to show that there exists $m_0$ such that $\mathbb{Q}_\theta(\widehat{\Psi}_{m_0}) < \infty$, then the majorant $\widehat{\Psi}_{m_0}$ is integrable. Arguing as in (6.20) we have that

$$\mathbb{Q}_\theta(\widehat{\Psi}_m) = \sum_{k \in \mathbb{N}_0^{\mathcal{S}_m}} \tau_m^{(+)}(k)|k| \prod_{s \in \mathcal{S}_m} \frac{\theta_m^{(-)}(s)^{k_s}}{k_s!}, \quad (6.23)$$

where $\theta_m^{(-)}(s) = e^{-(\kappa^{(m,-)} c_m)(s)} c_m(s)$, for $s \in \mathcal{S}_m$. Let $\chi_m := \chi(\kappa^{(m,+)}, \theta_m^{(-)})$ be defined as in (4.11). We can argue as in the proof of Lemma 4.5 to get that, for any $n \in \mathbb{N}$,

$$\sum_{k \in \mathbb{N}_0^{\mathcal{S}_m}\,:\,|k|=n} \tau_m^{(+)}(k)|k| \prod_{s \in \mathcal{S}_m} \frac{(\theta_m^{(-)}(s))^{k_s}}{k_s!} \leq e^{o(n)} e^{-n[\chi_m - 1]}, \quad (6.24)$$

(where the $e^{o(n)}$-term is actually given by $|\mathcal{M}_1^{(n)}(\mathcal{S}_m)| \sum_{r \in \mathcal{S}_m} \Delta_r(n\nu)$). Abbreviating $\delta_m = \|\kappa^{(m,+)} - \kappa^{(m,-)}\|_\infty$, we further have that

$$\chi_m = \inf_{\nu \in \mathcal{M}_1(\mathcal{S}_m)} \left\{\sum_{r \in \mathcal{S}_m} \nu_r \log \frac{\nu_r}{(\kappa^{(m,+)}\nu)_r (\theta_m^{(+)})_r} - \langle \kappa^{(m,+)} c_m - \kappa^{(m,-)} c_m, \nu\rangle\right\}$$
$$\geq \chi(\kappa^{(m,+)}, \theta_m^{(+)}) - \delta_m = \Sigma_m - \log \Sigma_m - \delta_m,$$

where $\Sigma_m := \Sigma(\kappa^{(m,+)}, c_m)$. Now, choose $\varepsilon > 0$ small enough such that $\Sigma(\kappa, c) + \varepsilon < 1$ and $\Sigma(\kappa, c) - \log(\Sigma(\kappa, c)) > 1 + \varepsilon$. It holds that $\Sigma_m \to \Sigma(\kappa, c)$, as $m \to \infty$, and clearly we have that $\delta_m \to 0$, as $m \to \infty$. Additionally, we use that the function $x \mapsto \phi(x) := x - \log x$ is continuous and decreasing on $[0, 1]$. Hence, we find $m_0$ such that $\delta_{m_0} \leq \varepsilon/2$, as well as $\Sigma_{m_0} \leq \Sigma(\kappa, c) + \varepsilon < 1$ and $\phi(\Sigma_{m_0}) \geq \phi(\Sigma(\kappa, c)) - \varepsilon/2$. Consequently,

$$\chi_{m_0} \geq \Sigma_{m_0} - \log \Sigma_{m_0} - \frac{\varepsilon}{2} \geq \Sigma(\kappa, c) - \log \Sigma(\kappa, c) - \varepsilon > 1,$$

which altogether implies that

$$\mathbb{Q}_\theta(\widehat{\Psi}_{m_0}) = \sum_{k \in \mathbb{N}_0^{\mathcal{S}_m}} \tau_m^{(+)}(k) |k| \prod_{s \in \mathcal{S}_m} \frac{(\theta_m^{(-)}(s))^{k_s}}{k_s!} \leq \sum_{n \in \mathbb{N}} e^{o(n)} e^{-n[\chi_{m_0} - 1]} < \infty. \tag{6.25}$$

Thus, Lebesgue's theorem of dominated convergence is applicable and (6.18) follows.

Equation (6.16) can be shown in the same way and relies on the discrete version of the equation, derived in Lemma 6.5. □

**Lemma 6.8** *The statement of Lemma 6.7 is also true under the assumption $\Sigma(\kappa, c) = 1$.*

**Proof** The idea is to construct a sequence $c^{(n)} \in \mathcal{M}(\mathcal{S})$ with $\Sigma(\kappa, c^{(n)}) < 1$ such that $\theta^{(n)} := \theta_{c^{(n)}} \nearrow \theta_c =: \theta$ monotonically as $n \to \infty$.

Recall that $\mathcal{S}$ is compact and $\kappa$ is continuous, hence a standard argument (see e.g. [7, Lemma 5.15]) shows that the operator $T_{\kappa,c}$ is a positive Hilbert–Schmidt operator and therefore has a non-negative eigenfunction corresponding to the eigenvalue $\Sigma(\kappa, c)$. By the assumption $\Sigma(\kappa, c) = 1$ we can find a function $g \colon \mathcal{S} \to (0, \infty)$ such that $T_{\kappa,c} g = g$. For any $n \in \mathbb{N}$ define $c^{(n)} \in \mathcal{M}(\mathcal{S})$ via $\frac{dc^{(n)}}{dc} := 1 - \frac{1}{n} g$. Then for $n$ large enough $c^{(n)}(A) < c(A)$ for any measurable $A \subset \mathcal{S}$ with $\int_A g \, dc > 0$. In particular $\Sigma(\kappa, c^{(n)}) < 1$. (An ad-hoc argument in the case that $\kappa$ is irreducible is as follows: Pick an $L^2(c^{(n)})$-normalized positive eigenfunction $g_n$ of $T_{\kappa,c^{(n)}}$ corresponding to the eigenvalue $\Sigma(\kappa, c^{(n)})$ and observe that $\widetilde{g}_n(x) = g_n(x)(1 - \frac{1}{n} g(x))^{1/2}$ is $L^2(c)$ normalized and that $\Sigma(\kappa, c^{(n)}) = \|T_{\kappa,c^{(n)}} g_n\|_{L^2(c^{(n)})} < \|T_{\kappa,c} \widetilde{g}_n\|_{L^2(c)} \leq \|T_{\kappa,c}\| = \Sigma(\kappa, c)$. If $\kappa$ is reducible, then apply this argument to the irreducible components.) Now, observe that for any $n \in \mathbb{N}$ we have

$$\frac{d\theta^{(n)}}{d\theta} = e^{\kappa(c - c^{(n)})} \frac{dc^{(n)}}{dc} = e^{T_{\kappa,c}(\frac{1}{n} g)} \left(1 - \frac{1}{n} g\right) = e^{\frac{1}{n} g} \left(1 - \frac{1}{n} g\right),$$

and the right-hand side converges pointwise monotonically to 1, as $n \to \infty$.

Now, fix any continuous test function $f \colon \mathcal{S} \to [0, \infty)$. Then, by monotone convergence and the fact that we can apply Lemma 6.7 to $c^{(n)}$ for all $n \in \mathbb{N}$, we get that

$$\begin{aligned}
\mathbb{Q}_\theta \big[ \tau(k) \langle k, f \rangle e^{\theta(\mathcal{S})} \big] &= \lim_{n \to \infty} \mathbb{Q}_\theta \big[ \tau(k) \langle k, f \rangle e^{\theta(\mathcal{S})} e^{\langle k, \log \frac{d\theta^{(n)}}{d\theta} \rangle} \big] \\
&= \lim_{n \to \infty} \mathbb{Q}_{\theta^{(n)}} \big[ \tau(k) \langle k, f \rangle e^{\theta^{(n)}(\mathcal{S})} \big] \\
&= \lim_{n \to \infty} \langle c^{(n)}, f \rangle = \langle c, f \rangle.
\end{aligned}$$

The same argument shows that $\mathbb{Q}_\theta\left[\tau(k)e^{\theta(\mathcal{S})}\right] = c(\mathcal{S}) - \frac{1}{2}\langle c, \kappa c\rangle$. □

Now we can identify the minimizers of $I_{\text{Mi}}$. The following is a variant of Proposition 6.2 in the general setting. Once having established Lemmas 6.7 and 6.8, the proof in the general setting follows the ones in the discrete setting. Recall that $\mu$ is the reference probability measure on $\mathcal{S}$.

**Lemma 6.9** (Minimizers of $I_{\text{Mi}}$) *Assume that $c \in \mathcal{M}(\mathcal{S})$ with $c \leq \mu$ satisfying $\Sigma(\kappa, c) \leq 1$. Then the unique minimizer $\lambda$ of $I_{\text{Mi}}$ under the assumption $c_\lambda = c$ is equal to $\lambda_c$ defined in (2.1) and*

$$\min_{\lambda \in \mathcal{L}:\, c_\lambda = c} I_{\text{Mi}}(\lambda) = I_{\text{Mi}}(\lambda_c) = -\frac{1}{2}\langle c, \kappa c\rangle + \langle c, \log \tfrac{dc}{d\mu}\rangle.$$

**Proof** Note that $\lambda_c$ is admissible, according to Lemmas 6.7 and 6.8, since $c_{\lambda^*} = c$.

Using Lemma 6.6 and the fact that $\mu$ is a probability measure we can rewrite the measure $\mathbb{Q}_\mu$ as

$$\mathbb{Q}_\mu(\mathrm{d}k) = e^{\theta(c)(\mathcal{S})}\mathbb{Q}_{ce^{-\kappa*c}}(\mathrm{d}k)e^{\langle k, \kappa*c\rangle}e^{-\langle k, \log \tfrac{dc}{d\mu}\rangle}e^{-1}$$

Now, writing $\mathcal{M} = \mathcal{M}_{\mathbb{N}_0}(\mathcal{S})$ we get for any $\lambda \in \mathcal{L}$ satisfying $c_\lambda = c$ that

$$\begin{aligned}
I_{\text{Mi}}(\lambda) &= \left\langle \lambda, \log \frac{\mathrm{d}\lambda}{\tau \mathrm{d}\mathbb{Q}_\mu}\right\rangle + c(\mathcal{S}) - 2\lambda(\mathcal{M}) + \frac{1}{2}\langle c, \kappa\mu\rangle \\
&= \left\langle \lambda, \log \frac{\mathrm{d}\lambda}{\mathrm{d}\lambda_c}\right\rangle + \left\langle c, \log \frac{\mathrm{d}c}{\mathrm{d}\mu}\right\rangle - \langle c, \kappa c\rangle + c(\mathcal{S}) - \lambda(\mathcal{M}) + \frac{1}{2}\langle c, \kappa\mu\rangle \\
&= \mathbb{H}(\lambda|\lambda_c) - \lambda_c(\mathcal{M}) + \left\langle c, \log \frac{\mathrm{d}c}{\mathrm{d}\mu}\right\rangle - \langle c, \kappa c\rangle + c(\mathcal{S}) + \frac{1}{2}\langle c, \kappa\mu\rangle \\
&= \mathbb{H}(\lambda|\lambda_c) + \left\langle c, \log \frac{\mathrm{d}c}{\mathrm{d}\mu}\right\rangle + \frac{1}{2}\langle c, \kappa(\mu - c)\rangle.
\end{aligned}$$

We used the fact that $\lambda_c(\mathcal{M}) = c(\mathcal{S}) - \frac{1}{2}\langle c, \kappa c\rangle$, which was derived in the last statement of Lemma 6.7. Since $\mathbb{H}(\lambda|\lambda_c) \geq 0$ and $\mathbb{H}(\lambda|\lambda_c) = 0$ if and only if $\lambda = \lambda_c$, the claim follows. □

### 6.3 The discrete, supercritical case

In this section, we assume again that $\mathcal{S}$ is a finite space and investigate the case where the measure $c$ (the one that formulates the constraint) is supercritical, meaning $\Sigma(\kappa, c) > 1$. The aim of this section is to verify the following result.

**Proposition 6.10** (Discrete, supercritical case) *Let $c \in [0, \mu]$ with $\Sigma(\kappa, c) > 1$. Then*

$$\inf_{\lambda:\, c(\lambda)=c} I_{\text{Mi}}(\lambda) = \inf_{\lambda:\, c(\lambda)=b^*} I_{\text{Mi}}(\lambda) + I_{\text{Me}}(c - b^*) \tag{6.26}$$

where $b^* = b^*(c)$ is the minimal non-trivial (i.e., not equal to $c$) solution to

$$\kappa(c - b^*)b^* = c - b^*, \qquad b^* \leq c, \tag{6.27}$$

and satisfies $\Sigma(\kappa, b^*) = 1$.

Indeed, one possible realization of the rate (6.26) is given by constructing a minimizer as in formula (6.3) with respect to the (sub-)critical parameter $b^*$ and realizing the remaining part $c - b^*$ by means of a diverging sequence $k^{(n)}$, such that the mesoscopic rate term appears.

The proof will be a consequence of the next lemmas. In Lemma 6.11 we derive an upper and a lower bound for $\inf_{\lambda:\, c_\lambda = c} I_{\text{Mi}}(\lambda)$ and in Lemma 6.12 we show that they coincide, if there are solutions to the fixed point equation (6.27). We will postpone the proof about existence of solutions to Sect. 6.4, Lemma 6.14.

**Lemma 6.11** *Let $c \in [0, \mu]$ with $\Sigma(\kappa, c) > 1$. For $b \in [0, c]$ with $\Sigma(\kappa, b) \leq 1$ we put*

$$F_c(b) = \left\langle c, \log \frac{b}{\mu} \right\rangle + |c - b| + \frac{1}{2} \langle b, \kappa b \rangle - \langle c, \kappa b \rangle + \frac{1}{2} \langle c, \kappa \mu \rangle, \tag{6.28}$$

$$G_c(b) = \left\langle b, \log \frac{b}{\mu} \right\rangle - \frac{1}{2} \langle b, \kappa b \rangle + \left\langle c - b, \log \frac{c - b}{(\kappa(c - b))\mu} \right\rangle + \frac{1}{2} \langle c, \kappa \mu \rangle. \tag{6.29}$$

*Then*

$$\sup_{b \in [0,c]:\, \Sigma(\kappa,b) \leq 1} F_c(b) \leq \inf_{\lambda:\, c(\lambda) = c} I_{\text{Mi}}(\lambda) \leq \inf_{b \in [0,c]:\, \Sigma(\kappa,b) \leq 1} G_c(b). \tag{6.30}$$

*Proof* We first show the first inequality in (6.30). Fix $b \in [0, c]$ with $\Sigma(\kappa, b) \leq 1$. Let $\lambda^* := \lambda(b)$ be given as in (6.3). We proceed in the same way as in the proof of Proposition 6.2, but use that this time $|\lambda^*| = |b| - \frac{1}{2} \langle b, \kappa b \rangle$. Then for any $\lambda$ with $c(\lambda) = c$ we have

$$I_{\text{Mi}}(\lambda) = \mathbb{H}(\lambda|\lambda^*) - |\lambda^*| + \left\langle c, \log \frac{b}{\mu} \right\rangle - \langle c, \kappa b \rangle + |c| + \frac{1}{2} \langle c, \kappa \mu \rangle \geq F_c(b).$$

We now prove the upper bound in (6.30). Fix $b \in [0, c]$ with $\Sigma(\kappa, b) \leq 1$.

*Case 1*: First, we assume that $\kappa$ is irreducible with respect to $c - b$. Let $\lambda^* := \lambda(b)$ be defined as in (6.3). For $n \in \mathbb{N}$ define $k^{(n)} := \lfloor n(c - b) \rfloor$ and write $R_n := |k^{(n)}|$ and $b^{(R_n)} := \sum_{|k| \leq R_n} \lambda_k^* k$. We define

$$\lambda_k^{(n)} := \begin{cases} \lambda_{\mathbf{e}_s}^* + \varepsilon_s & \text{if } k = \mathbf{e}_s \text{ for some } s \in \mathcal{S} \\ \lambda_k^* & \text{if } 1 < |k| < R_n \\ \frac{1}{n} & \text{if } k = k^{(n)} \\ 0 & \text{else} \end{cases} \tag{6.31}$$

with $\varepsilon := c - b^{(R_n)} - \frac{1}{n}k^{(n)}$, which ensures that $c(\lambda^{(n)}) = c$ holds for all $n \in \mathbb{N}$, but is negligible in the limit, i.e., $\lim_{n \to \infty} \lambda^{(n)}_{\mathbf{e}_s} = \lambda^*_{\mathbf{e}_s}$ for all $s \in \mathcal{S}$. Note that due to the irreducibility assumption we have that $\tau(k^{(n)}) > 0$ if $n$ is large enough, which ensures that $I_{\text{Mi}}(\lambda^{(n)})$ is finite.

Using the notation $I^{(R)}_{\text{Mi}}(\lambda)$ introduced in (3.22) we get

$$I_{\text{Mi}}(\lambda^{(n)}) = I^{(R_n-1)}_{\text{Mi}}(\lambda^{(n)}) + \frac{1}{2}\langle c - b^{(R_n)}, \kappa\mu\rangle + \frac{1}{n}\log \frac{\frac{1}{n}\prod_s k_s^{(n)}!}{\tau(k^{(n)})\mathrm{e}^{1-|k^{(n)}|}\prod_s \mu_s^{k_s^{(n)}}}$$

Denote the last summand as $A_n$. By using the formula (4.6) from Lemma 4.3 for some $r \in \text{supp}(c)$ as well as Stirling's formula for the factorial terms, we have that, as $n \to \infty$

$$A_n = o(1) + \sum_{s \in \mathcal{S}}(c_s - b_s)\log \frac{c_s - b_s}{(\kappa(c-b))_s)\mu_s} + \frac{1}{n}\log \frac{n(c_r - b_r)\prod_s \kappa(c-b)_s}{\mathrm{e}\Delta_r(c-b)}, \tag{6.32}$$

where $\Delta_r$ is defined in (4.7), which can be easily extended to arguments in $[0, \infty)^{\mathcal{S}}$. Note that by construction $\Delta_r(c-b) > 0$. Clearly, the last summand in (6.32) is of order $o(1)$. By the construction of $\lambda^{(n)}$ it is immediate that $\lim_{n \to \infty} I^{(R_n-1)}_{\text{Mi}}(\lambda^{(n)}) = \langle b, \log(b/\mu)\rangle + \frac{1}{2}\langle b, \kappa(\mu-b)\rangle$, so altogether we get that

$$\lim_{n\to\infty} I_{\text{Mi}}(\lambda^{(n)}) = \left\langle b, \log\frac{b}{\mu}\right\rangle + \frac{1}{2}\langle b, \kappa(\mu-b)\rangle + \langle c-b, \log\frac{c-b}{(\kappa(c-b))\mu}\rangle$$
$$+ \frac{1}{2}\langle c-b, \kappa\mu\rangle = G_c(b).$$

*Case 2*: If $\kappa$ is reducible with respect to $c - b$, we can find a decomposition of $\text{supp}(c-b)$ into disjoint sets $\mathcal{S}_j$ such that $\kappa$ restricted to $\mathcal{S}_j \times \mathcal{S}_j$ is irreducible and $\kappa|_{\mathcal{S}_i \times \mathcal{S}_j} = 0$ for $i \neq j$. Then we have to modify the construction of $\lambda^{(n)}$ given above by putting mass $\frac{1}{n}$ on each of the meso-particles $k^{(j,n)} := \lfloor n(c-b)\rfloor \mathbb{1}_{\mathcal{S}_j}$. We omit the details. □

The following lemma completes the proof of Proposition 6.10. Its assumptions are verified later and in more generality in Lemma 6.14.

**Lemma 6.12** *Let $c \in [0, \mu]$ with $\Sigma(\kappa, c) > 1$. If there exists a non-trivial solution $b^* \in [0, c]$ to (6.27) and $\Sigma(\kappa, b^*) = 1$, then $F(b^*) = G(b^*)$. Consequently, Eq. (6.26) holds.*

**Proof** Using the fixed point equation (6.27) we can substitute $|c-b^*| = \langle b^*, \kappa(c-b^*)\rangle$ to rewrite $F_c(b^*)$ and $\langle c-b^*, \log[(c-b^*)/\kappa(c-b^*)]\rangle = \langle c-b^*, \log b^*\rangle$ to rewrite $G_c(b^*)$. Then $F_c(b^*) = G_c(b^*)$. □

### 6.4 The general supercritical case

Building on the results of the previous subsection we derive a slightly weaker result than Proposition 6.10 for the general case, which will still be enough to derive the optimal rates for the contraction principle as well as to fully optimize the rate function $I$.

**Lemma 6.13** *Fix $c \in \mathcal{M}(\mathcal{S})$ with $c \leq \mu$ and $\Sigma(\kappa, c) > 1$. Then*

$$\inf_{\lambda \in \mathcal{L}:\, c_\lambda = c} I_{\mathrm{Mi}}(\lambda) \geq \inf_{\lambda \in \mathcal{L}:\, c_\lambda = b^*} I_{\mathrm{Mi}}(\lambda) + I_{\mathrm{Me}}(c - b^*), \tag{6.33}$$

*where $b^* = b^*(c) \in \mathcal{M}(\mathcal{S})$ is the minimal, non-trivial (i.e., not equal to c) solution to* (2.9) *and satisfies $\Sigma(\kappa, b^*) = 1$.*

***Sketch of proof*** We can generalize the proof of the lower bound of Lemma 6.11 and the definition of $F_c$ to obtain

$$\inf_{\lambda:\, c_\lambda = c} I_{\mathrm{Mi}}(\lambda) \geq \left\langle c, \log \frac{\mathrm{d}b}{\mathrm{d}\mu} \right\rangle + c(\mathcal{S}) - b(\mathcal{S}) + \frac{1}{2}\langle b, \kappa b\rangle - \langle c, \kappa b\rangle =: F_c(b)$$

for any $b \in \mathcal{M}(\mathcal{S})$ with $b \leq c$ and $\Sigma(\kappa, b) \leq 1$. This relies on the admissibility of the (auxiliary) minimizers $\lambda_b$ proved in Lemmas 6.7 and 6.8. The lower bound is obtained by writing everything with entropies as in the proof of Lemma 6.9.

Observe that if $b^*$ is a solution of (2.9), then one can argue as in the proof of Lemma 6.12 to see that

$$F_c(b^*) = \left\langle b^*, \log \frac{\mathrm{d}b^*}{\mathrm{d}\mu} \right\rangle - \frac{1}{2}\langle b^*, \kappa b^*\rangle + I_{\mathrm{Me}}(c - b^*) + \frac{1}{2}\langle c, \kappa \mu\rangle =: G_c(b^*)$$

□

**Lemma 6.14** (Solutions to (2.9)) *Fix $c \in \mathcal{M}(\mathcal{S})$.*

*(i) Assume that $\kappa$ is irreducible w.r.t. $c$ and $\Sigma(\kappa, c) > 1$. Then there exists exactly one solution $b^*$ to* (2.9) *that satisfies $b^* \neq c$. Further, it holds that $\Sigma(\kappa, b^*) = 1$.*
*(ii) If $\Sigma(\kappa, c) \leq 1$, then the only solution $b^*$ to* (2.9) *is given by the trivial solution $b^* = c$.*
*(iii) Assume that $\kappa$ is reducible w.r.t. $c$ and $\Sigma(\kappa, c) > 1$, then there exists at least one solution $b^*$ to* (2.9) *with $b^* \neq c$. Moreover, there exists a unique minimal solution $b_*$ to* (2.9) *(which is minimal in the sense that $b_* \leq b^*$ holds c-almost everywhere for all solutions $b^*$ of* (2.9)*). Further, we have that $\Sigma(\kappa, b^*) > 1$ for all solutions $b^*$ with $b^* \neq b_*$ and $\Sigma(\kappa, b_*) = 1$.*

***Proof*** We will study the existence and uniqueness of non-trivial solutions $f^* \colon \mathcal{S} \to [0, 1)$ to

$$T_{\kappa,c} f^* = \frac{f^*}{1 - f^*}. \tag{6.34}$$

By substituting $b^* = (1 - f^*)c$ it is easily seen that solving (6.34) is equivalent to solving (2.9).

(i) *Existence*: We once more reformulate (6.34) by substituting $g^* = f^*/(1 - f^*)$ (which is equivalent to $f^* = g^*(1 + g^*)$). Then (6.34) is equivalent to

$$U(g^*) := T_{\kappa,c}\Big(\frac{g^*}{1+g^*}\Big) = g^*, \qquad (6.35)$$

i.e., we are searching for a fixed point of $U$. Note that $g^*(s)/(1 + g^*(s)) \leq 1$ for all $s \in \mathcal{S}$. Together with the fact that $\kappa$ is non-negative this implies that any solution $g^*$ of (6.35) satisfies $g^* \leq T_{\kappa,c}\mathbb{1}$. Hence, it suffices to study the operator $U$ on the domain $D = \{g \colon \mathcal{S} \to \mathbb{R} \colon 0 \leq g \leq T_{\kappa,c}\mathbb{1}\}$. We construct a solution iteratively by defining $g_0 := T_{\kappa,c}\mathbb{1}$ and $g_m := U(g_{m-1})$ for $m \in \mathbb{N}$. Since the function $x \mapsto x/(1 + x)$ is strictly increasing on $[0, \infty)$ and $\kappa$ is non-negative, we have that $g \leq g'$ implies $U(g) \leq U(g')$. Since $g_1 \leq g_0$ we can iterate this argument to show that $g_m \leq g_{m-1}$ holds for any $m \in \mathbb{N}$. Therefore the limit $g^* := \lim_{m \to \infty} g_m \in D$ exists and by the continuity of $U$ it satisfies (6.35). We claim that our assumptions on $\kappa$ and $c$ imply that $g^* > 0$: By the assumptions that $\kappa$ is irreducible w.r.t. $c$, $\mathcal{S}$ is compact and $\kappa$ is continuous, $T_{\kappa,c}$ is a positive, irreducible and compact operator. Therefore there exists a strictly positive eigenfunction $v$ of $T_{\kappa,c}$ with eigenvalue $\Sigma(\kappa, c) > 1$. Note that the function $T_{\kappa,c}v$ is continuous (by compactness of $\mathcal{S}$ and continuity of $\kappa$), hence $v$ is continuous and by compactness of $\mathcal{S}$ it is also bounded. So without loss of generality we can pick $v$ such that $v(s) \in (0, 1]$ for any $s \in \mathcal{S}$. Observe that by the irreducibility assumption $g_0 = T_{\kappa,c}\mathbb{1} > 0$. Now, pick $\delta > 0$ such that $\Sigma(\kappa, c) \geq 1 + \delta$ and $g_0 \geq \delta v$. Observe that for any $g$ with $g \geq \delta v$ we have that

$$U(g) = T_{\kappa,c}\Big(\frac{g}{1+g}\Big) \geq T_{\kappa,c}\Big(\frac{\delta v}{1+\delta v}\Big) \geq \frac{\delta}{1+\delta}T_{\kappa,c}v \geq \delta v.$$

Hence $g_m \geq \delta v$ holds for all $m \in \mathbb{N}$. Consequently, $g^* \geq \delta v > 0$.

Additionally, we claim that $g^*$ is the maximal solution to (6.35). Let $\tilde{g}^*$ be any other solution to (6.35), then we necessarily have that $\tilde{g}^* \leq T_{\kappa,c}\mathbb{1} = g_0$. It follows by the monotonicity of $U$ that $\tilde{g}^* = U(\tilde{g}^*) \leq g_1$ and, iteratively, $\tilde{g}^* \leq g_m$ for any $m \in \mathbb{N}$. Hence $\tilde{g}^* \leq g^*$. By the equivalence of (6.34) and (6.35) and monotonicity of $x \mapsto x/(1 + x)$ we have that $f^* = g^*/(1 + g^*)$ is the maximal solution to (6.34).

*Uniqueness*: Assume towards a contradiction that $f^*$ and $\tilde{f}^*$ are non-trivial solutions of (6.34) and $f^* \neq \tilde{f}^*$ on a set $A \subset \mathcal{S}$ with $c(A) > 0$. Without loss of generality we can assume that $f^*$ is the maximal solution (as constructed in the existence part), i.e., $\tilde{f}^* \leq f^*$ on $\mathcal{S}$ and $\tilde{f}^* < f^*$ everywhere on $A$. For any $h \colon \mathcal{S} \to [0, 1]$ put $\Psi(h) := (1 - h)T_{\kappa,c}h$. Then $\Psi(f^*) = f^*$ and $\Psi(\tilde{f}^*) = \tilde{f}^*$ and we have that

$$\Psi\Big(\frac{f^*}{2} + \frac{\tilde{f}^*}{2}\Big) = \frac{1}{2}\Psi(f^*) + \frac{1}{2}\Psi(\tilde{f}^*) + \frac{1}{4}(f^* - \tilde{f}^*)T_{\kappa,c}(f^* - \tilde{f}^*)$$
$$\geq \frac{1}{2}\Psi(f^*) + \frac{1}{2}\Psi(\tilde{f}^*) = \frac{1}{2}f^* + \frac{1}{2}\tilde{f}^*,$$

where the inequality relies on the fact that $\kappa$ is non-negative and $\tilde{f}^* \leq f^*$. Note that we even have a strict inequality on the set $A$. Now, define $h := (f^* - \tilde{f}^*)/2$, then $\tilde{f}^* + h = (\tilde{f}^* + f^*)/2$ and we already argued that $\Psi(\tilde{f}^* + h) \geq \tilde{f}^* + h$ where the inequality is strict on $A$. On the other hand, as $h \geq 0$, we get

$$\Psi(\tilde{f}^* + h) = (1 - \tilde{f}^* - h)T_{\kappa,c}(\tilde{f}^* + h) \leq (1 - \tilde{f}^*)T_{\kappa,c}(\tilde{f}^* + h),$$

so altogether $(1 - \tilde{f}^*)T_{\kappa,c}(\tilde{f}^* + h) \geq \tilde{f}^* + h$ with strict inequality on $A$. Using (6.34) for $\tilde{f}^*$ and the symmetry of $\kappa$ we get

$$\langle c, \tilde{f}^* T_{\kappa,c}(\tilde{f}^* + h) \rangle > \left\langle c, \frac{\tilde{f}^*}{1 - \tilde{f}^*}(\tilde{f}^* + h) \right\rangle = \langle c, (T_{\kappa,c}\tilde{f}^*)(\tilde{f}^* + h) \rangle$$
$$= \langle c, \tilde{f}^* T_{\kappa,c}(\tilde{f}^* + h) \rangle,$$

which is a contradiction. Hence the solution to (6.34) is unique up to sets that have measure zero w.r.t. $c$, which implies uniqueness of $b^*$.

We now argue that $\Sigma(\kappa, b^*) = 1$. Our procedure is very similar to the one used in the proof of Lemma 6.6 in [7]. Let $w \colon \mathcal{S} \to \mathbb{R}$ be an eigenfunction of $T_{\kappa,b^*}$ with eigenvalue $a$. Then using (6.34) and the symmetry of $\kappa$, we get that

$$\langle c, w f^* \rangle = \langle c, w(1 - f^*)T_{\kappa,c} f^* \rangle = \langle c, f^* T_{\kappa,c}(w(1 - f^*)) \rangle$$
$$= \langle c, f^* T_{\kappa,b^*} w \rangle = a \langle c, f^* w \rangle. \qquad (6.36)$$

Hence, we either have that $\langle c, w f^* \rangle = 0$ or $a = 1$. By the Krein–Rutman Theorem (the extension of the Perron–Frobenius Theorem to positive compact operators) the eigenfunction $w$ corresponding to the largest eigenvalue of $T_{\kappa,b^*}$ is non-negative and non-trivial, so $\langle c, w f^* \rangle > 0$ and hence $\Sigma(\kappa, b^*) = 1$. (Interestingly, we have constructed $b^*$ in such a way that all other eigenfunctions $\widetilde{w}$ of $T_{\kappa,b^*}$ satisfy $\langle c, \widetilde{w} f^* \rangle = 0$.)

(ii) Let $\Sigma(\kappa, c) \leq 1$. Assume towards a contradiction that $f^*$ is a solution to (6.34) and $\int_A f^* \, dc > 0$ for some open set $A \subset \mathcal{S}$ such that $c(A) > 0$. Using the substitution $f^* = g/(1+g)$ we have that Eq. (6.34) is equivalent to $T_{\kappa,c}(g/(1+g)) = g$ and since $\kappa$ is continuous and $\mathcal{S}$ is compact the left-hand side $T_{\kappa,c}(g/(1+g))$ is a continuous function, which implies that both $g$ and $f^*$ are continuous functions. So we can find an $\varepsilon > 0$ and a set $A_\varepsilon \subset A$ such that $f^* \geq \varepsilon$ on $A_\varepsilon$ and $c(A_\varepsilon) > 0$. Therefore,

$$T_{\kappa,c} f^* = \frac{f^*}{1 - f^*} \geq \frac{f^*}{1 - \varepsilon}$$

holds on $A_\varepsilon$ and $T_{\kappa,c} f^* \geq f^*$ holds on $\mathcal{S}$. This implies that $\|T_{\kappa,c} f^*\|_{L^2(c)} > \|f^*\|_{L^2(c)}$ and hence $\Sigma(\kappa, c) > 1$ in contradiction to our assumption.

(iii) Since $\kappa$ is reducible w.r.t. $c$, we find a decomposition of $\mathrm{supp}(c)$ into (countable many) disjoint sets $S_j$, $j \in J$, such that $\kappa^{(j)} = \kappa|_{S_j \times S_j}$ is irreducible with respect to $c^{(j)}$, the restriction of $c$ to $S_j$ for any $j \in J$, and $\kappa|_{S_i \times S_j} = 0$ holds $c$-almost everywhere, if $i \neq j$. Let $J' := \{j \in J : \Sigma(\kappa^{(j)}, c^{(j)}) > 1\}$ and note that $J' \neq \emptyset$.

By (i) we get that for any $j \in J'$ there exists a function $f^{(j)}: S_j \to [0, 1)$ that solves (6.34) on $S_j$ and $f^{(j)} > 0$. By (ii) we get that for any $j \in J \setminus J'$ the only function $f: S_j \to [0, 1)$ that solves (6.34) on $S_j$ is equal to $0$ $c^{(j)}$-almost everywhere. Now for $\sigma = (\sigma_j)_{j \in J'} \in \{0, 1\}^{J'}$ define $f^{(\sigma)}: \mathcal{S} \to [0, 1)$ by $f^{(\sigma)}(s) = \sigma_j f^{(j)}(s)$, if $s \in S_j$ for some $j \in J'$ and $f^{(\sigma)}(s) = 0$ for $s \in \mathcal{S} \setminus \bigcup_{j \in J'} S_j$. It can now be easily checked that all solutions to (6.34) are given by

$$\mathcal{F} := \left\{ f^{(\sigma)}: \sigma = (\sigma_j)_{j \in J'} \in \{0, 1\}^{J'} \right\} \tag{6.37}$$

and that $\mathcal{F}$ contains at least one non-trivial solution. Write $b^{(\sigma)} = (1 - f^{(\sigma)})c$ and note that all possible solutions of (2.9) are of this form. Clearly, the minimal solution $b^*$ of (2.9) is given via the maximal solution in $\mathcal{F}$, i.e., by choosing $\sigma \equiv 1$.

We will now investigate the quantities $\Sigma(\kappa, b^{(\sigma)})$ for any choice of $\sigma$. First, let $\sigma$ be such that there exists some $j \in J'$ with $\sigma_j = 0$. Then for $r \in S_j$ and any function $h: \mathcal{S} \to \mathbb{R}$

$$T_{\kappa, b^{(\sigma)}} h(r) = \int_{S_j} \kappa(r, s) h(s) (1 - \sigma_j f^{(j)}(s)) \, c(\mathrm{d}s) = T_{\kappa^{(j)}, c^{(j)}} h(r).$$

Therefore, given an eigenfunction $g^{(j)}$ of $T_{\kappa^{(j)}, c^{(j)}}$ that corresponds to the eigenvalue $\Sigma(\kappa^{(j)}, c^{(j)})$, we can construct an eigenfunction $g$ for $T_{\kappa, b^{(\sigma)}}$ with the same eigenvalue by choosing $g = g^{(j)}$ on $S_j$ and $g = 0$ on $\mathcal{S} \setminus S_j$. Hence, $\Sigma(\kappa, b^{(\sigma)}) = \Sigma(\kappa^{(j)}, c^{(j)}) > 1$. Now, consider $\sigma \equiv 1$. Then we can argue as in (6.36) to show that $\Sigma(\kappa, b^{(\sigma)}) = 1$. □

## 7 Analysis of minimizers of the rate function in Theorem 1.1

In this section we provide the final steps needed for the optimization of the rate function and prove Theorems 2.3 and 2.1. Since the arguments for the remaining steps do not rely on discrete combinatorics (as it was the case in Sect. 6), we will immediately work in the general setting. In Sect. 7.1 we study a constrained optimization problem for the functions $I_{\text{Me}}$ and $I_{\text{Ma}}$. In Sect. 7.2 we prove the explicit form of the rate functions that is derived by applying the contraction principle and formulated in Theorem 2.3. Section 7.3 presents the last step for a full optimization of the rate function $I$ that gives Theorem 2.1. In Sect. 7.4 we derive the Flory equation that we formulated in Proposition 2.7.

### 7.1 The minimizers of $I_{\text{Ma}}$ and $I_{\text{Me}}$

Complementary to what was done in Sect. 6 for the function $I_{\text{Mi}}$ we will solve the analogous optimization problems for the functions $I_{\text{Me}}$ and $I_{\text{Ma}}$ defined in (1.12) and (1.13). To optimize the function $I_{\text{Me}}$ it is beneficial to combine it with $I_{\text{Ma}}$. We will again fix a measure $c \in \mathcal{M}(\mathcal{S})$ to formulate the constraint. In contrast to the result of Proposition 6.1 it will turn out that we do not have to distinguish between the cases $\Sigma(\kappa, c) \leq 1$ and $\Sigma(\kappa, c) > 1$.

**Lemma 7.1** (Minimizers of $I_{Me}$ and $I_{Ma}$) *Let $c \in \mathcal{M}(\mathcal{S})$ be such that $c \leq \mu$. Then*

$$\inf\{I_{Ma}(\alpha): \alpha \in \mathcal{A}, c(\alpha) = c\} = I_{Ma}(\delta_c), \tag{7.1}$$

$$\inf\{I_{Ma}(\alpha) + I_{Me}(\nu): \alpha \in \mathcal{A}, \nu \in \mathcal{M}(\mathcal{S}), c(\alpha) + \nu = c\} = I_{Ma}(\delta_c) + I_{Me}(0). \tag{7.2}$$

*If $\kappa$ is irreducible with respect to $c$, then the minimizers are unique.*

In order to prove the result above, we need the following lemma.

**Lemma 7.2** *Let $c \in \mathcal{M}(\mathcal{S})$ and let $\kappa$ be irreducible with respect to $c$. Let $\alpha \in \mathcal{A}$ be such that $c(\alpha) = c$ and assume that $\alpha = \sum_{i \in I} \delta_{y^{(i)}}$ with $|I| \geq 2$ and $I_{Ma}(\alpha) < \infty$. Then for any fixed $i \in I$ there is a measurable set $A \subset \mathcal{S}$ such that $y^{(i)}(A) > 0$ and $\kappa(c - y^{(i)})(x) > 0$ for all $x \in A$.*

*Proof* Denote $y := y^{(i)}$, $S_1 = \mathrm{supp}(y)$ and $S_2 = \mathrm{supp}(c(\alpha) - y)$. We first study the case where the sets $S_1$ and $S_2$ are disjoint. Assume towards a contradiction that for $y$-almost every $r$ we have that $\kappa(c - y)(r) = 0$. Then $\kappa|_{S_1 \times S_2} = 0$ which is equivalent to saying that $\kappa|_{S_1 \times \mathcal{S} \setminus S_1} = 0$ holds $c$-almost everywhere. Since $\kappa$ is irreducible with respect to $c$ we get that either $c(S_1) = 0$ or $c(\mathcal{S} \setminus S_1) = 0$. Consequently, either $y(S_1) = 0$ or $(c - y)(\mathcal{S} \setminus S_1) = 0$, i.e., either $y = 0$ or $c - y = 0$ in contradiction to our assumptions.

Now, assume that $S_1$ and $S_2$ are not disjoint. In that case we can pick $j \in I \setminus \{i\}$ in such a way that with $\hat{y} := y^{(j)}$ there exists $r_0 \in \mathrm{supp}(y) \cap \mathrm{supp}(\hat{y})$, which implies that for any open neighborhood $A_0 \subset \mathcal{S}$ of $r_0$ we have that $y(A_0) > 0$ and $\hat{y}(A_0) > 0$. By our assumption $I_{Ma}(\alpha) < \infty$, we have that $-\infty < \langle \hat{y}, \log(1 - e^{-\kappa \hat{y}}) \rangle$. This, together with the uniform continuity of $\kappa$ implies that we can find a neighborhood $A_0$ of $r_0$ such that $\kappa \hat{y}\big|_{A_0} > 0$. Now the claim follows, since $\kappa(c - y) \geq \kappa \hat{y}$ and $y(A_0) > 0$. □

Now we are ready to prove Lemma 7.1.

*Proof of Lemma 7.1* Writing $\widehat{I}_{Ma}(\alpha) = I_{Ma}(\alpha) - \frac{1}{2}\langle c(\alpha), \kappa\mu \rangle$ and $\widehat{I}_{Me}(\nu) = I_{Me}(\nu) - \frac{1}{2}\langle \nu, \kappa\mu \rangle$ it suffices to prove Eqs. (7.1) and (7.2) for $\widehat{I}_{Ma}$ and $\widehat{I}_{Me}$ since the difference does not depend on $c$. Observe that $\widehat{I}_{Ma}(\alpha) = A(\alpha) + B(\alpha)$, where

$$A(\alpha) = \int_{\mathcal{S}} \mu(\mathrm{d}r) \int \alpha(\mathrm{d}y) \frac{\mathrm{d}y}{\mathrm{d}\mu}(r) \log \frac{\frac{\mathrm{d}y}{\mathrm{d}\mu}(r)}{\kappa y(r)}, \tag{7.3}$$

$$B(\alpha) = \int_{\mathcal{S}} \mu(\mathrm{d}r) \int \alpha(\mathrm{d}y) \frac{\mathrm{d}y}{\mathrm{d}\mu}(r) \log \frac{\kappa y(r)}{e^{\frac{1}{2}\kappa y(r)} - e^{-\frac{1}{2}\kappa y(r)}}. \tag{7.4}$$

Let $\alpha \in \mathcal{A}$ with $c(\alpha) = c$. Note that $\int \alpha(\mathrm{d}y) \kappa y(r) = \kappa c(r)$ for $r \in \mathcal{S}$. With $\phi(u) = u \log u$ we use Jensen's inequality to get that

$$A(\alpha) = \int_{\mathcal{S}} \mu(\mathrm{d}r) \kappa c(r) \int \alpha(\mathrm{d}y) \frac{\kappa y(r)}{\kappa c(r)} \phi\left(\frac{\frac{\mathrm{d}y}{\mathrm{d}\mu}(r)}{\kappa y(r)}\right)$$

$$\geq \int_{\mathcal{S}} \mu(\mathrm{d}r) \kappa c(r) \phi \left( \int \alpha(\mathrm{d}y) \frac{\frac{\mathrm{d}y}{\mathrm{d}\mu}(r)}{\kappa c(r)} \right)$$

$$= \int_{\mathcal{S}} \mu(\mathrm{d}r) \frac{\mathrm{d}c}{\mathrm{d}\mu}(r) \log \frac{\frac{\mathrm{d}c}{\mathrm{d}\mu}(r)}{\kappa c(r)} = A(\delta_c), \tag{7.5}$$

where we used that $\int \alpha(\mathrm{d}y) \frac{\mathrm{d}y}{\mathrm{d}\mu}(r) = \frac{\mathrm{d}c}{\mathrm{d}\mu}(r)$. We now derive a corresponding lower bound for $B(\alpha)$. Note that the function $u \mapsto u/\sinh(u) =: \psi(u)$ is strictly decreasing on $[0, \infty)$ and that for any $y \in \mathrm{supp}(\alpha)$ we have that $\kappa y \leq \kappa c$. Therefore

$$B(\alpha) = \int_{\mathcal{S}} \mu(\mathrm{d}r) \int \alpha(\mathrm{d}y) \frac{\mathrm{d}y}{\mathrm{d}\mu}(r) \log \psi\left(\tfrac{1}{2}\kappa y(r)\right)$$

$$\geq \int_{\mathcal{S}} \mu(\mathrm{d}r) \int \alpha(\mathrm{d}y) \frac{\mathrm{d}y}{\mathrm{d}\mu}(r) \log \psi\left(\tfrac{1}{2}\kappa c(r)\right)$$

$$= \int_{\mathcal{S}} \mu(\mathrm{d}r) \frac{\mathrm{d}c}{\mathrm{d}\mu}(r) \log \psi\left(\tfrac{1}{2}\kappa c(r)\right) = B(\delta_c). \tag{7.6}$$

This implies Eq. (7.1). Now, let $\alpha \in \mathcal{A}$ and $\nu \in \mathcal{M}(\mathcal{S})$ be such that $c(\alpha) + \nu = c$. Note that

$$\widehat{I}_{\mathrm{Ma}}(\alpha) + \widehat{I}_{\mathrm{Me}}(\nu) = A(\alpha + \delta_\nu) + B(\alpha).$$

Since $c(\alpha + \delta_\nu) = c$ holds, we can use the estimate from before to get that $A(\alpha + \delta_\nu) \geq A(\delta_c)$. To get the estimate for $B(\alpha)$, observe that we still have that $\kappa y \leq \kappa c$ for any $y \in \mathrm{supp}(\alpha)$. So,

$$B(\alpha) \geq \int_{\mathcal{S}} \mu(\mathrm{d}r) \frac{\mathrm{d}c(\alpha)}{\mathrm{d}\mu}(r) \log \psi\left(\tfrac{1}{2}\kappa c(r)\right) \geq \int_{\mathcal{S}} \mu(\mathrm{d}r) \frac{\mathrm{d}c}{\mathrm{d}\mu}(r) \log \psi\left(\tfrac{1}{2}\kappa c(r)\right) = B(\delta_c),$$

where the second estimate is due to the the fact that $\psi(u) \in [0, 1]$ for $u \geq 0$ and thus $\log \psi(\tfrac{1}{2}\kappa c(r)) \leq 0$. Combining the estimates gives Eq. (7.2).

For both uniqueness claims we rely on Lemma 7.2 above.

To show uniqueness of the minimizer in Eq. (7.1) assume that $\alpha \in \mathcal{A}$ with $c(\alpha) = c$ and $\alpha \neq \delta_c$. Without loss of generality we can assume that $I_{\mathrm{Ma}}(\alpha) < \infty$. Now, we only have to note that the inequality in the estimate (7.6) is a strict inequality, since by Lemma 7.2 we have that $\kappa y < \kappa c$ holds on some set $A$, for which $y(A) > 0$ and the function $\psi$ is strictly decreasing.

To show uniqueness of the minimizer in Eq. (7.2) assume that $\alpha \in \mathcal{A}$ and $\nu \in \mathcal{M}(\mathcal{S})$ with $\nu \neq 0$ and $c(\alpha) + \nu = c$. Then by the same arguments as before

$$B(\alpha) \geq B(\delta_{c(\alpha)}) = \int c(\alpha)(\mathrm{d}r) \log \psi\left(\kappa c(\alpha)(r)\right) > \int c(\alpha)(\mathrm{d}r) \log \psi\left(\kappa c(r)\right) \geq B(\delta_c)$$

holds, if $\kappa c > \kappa c(\alpha)$ on some measurable set $A \subset \mathcal{S}$ for which $c(\alpha)(A) > 0$. To see that the latter condition is satisfied, we apply Lemma 7.2 to the measure $\delta_{c(\alpha)} + \delta_\nu$. This proves the claim. □

### 7.2 Minimization for the contraction principles

Here, we will exploit the work of Sects. 6 and 7.1 to prove Theorem 2.3, which is an application of the contraction principle but also provides an explicit solution for the optimization problem. When studying the optimization problem in the large deviation principle for $\mathrm{Ma}_N$, we encounter a functional that combines rates coming from the microscopic and the mesoscopic part. Its optimization is derived in the following lemma.

**Lemma 7.3** *Fix $c \in \mathcal{M}(\mathcal{S})$ with $c \leq \mu$. For $b \in \mathcal{M}(\mathcal{S})$ with $b \leq c$ and $\Sigma(\kappa, b) \leq 1$ let $G_c$ be as in* (6.29)

$$G_c(b) := \left\langle b, \log \frac{\mathrm{d}b}{\mathrm{d}\mu} \right\rangle - \frac{1}{2}\langle b, \kappa b\rangle + \left\langle c - b, \log \frac{\mathrm{d}(c-b)}{\kappa(c-b)\mathrm{d}\mu} \right\rangle + \frac{1}{2}\langle c, \kappa\mu\rangle.$$

*Then the following holds.*

1. *If $\Sigma(\kappa, c) \leq 1$, then*

$$\min\{G_c(b) \colon b \in \mathcal{M}(\mathcal{S}), b \leq c, \Sigma(\kappa, b) \leq 1\} = G_c(c), \qquad (7.7)$$

  *and $c$ is the unique minimizer.*
2. *If $\Sigma(\kappa, c) > 1$, then*

$$\min\{G_c(b) \colon b \in \mathcal{M}(\mathcal{S}), b \leq c, \Sigma(\kappa, b) \leq 1\} = G_c(b^*(c)) \qquad (7.8)$$

  *and $b^*(c)$ is the unique minimizer, which is given as the minimal, non-trivial (i.e., not equal to c) solution of* (2.9), *and it satisfies $\Sigma(\kappa, b^*(c)) = 1$.*

**Proof** (1) We use that $\frac{1}{2}\langle c, \kappa c\rangle - \frac{1}{2}\langle b, \kappa b\rangle = \langle b, \kappa(c-b)\rangle + \frac{1}{2}\langle c-b, \kappa(c-b)\rangle$ holds by the symmetry of $\kappa$. Therefore,

$$G_c(b) - G_c(c) = \left\langle b, \log \frac{\mathrm{d}b}{\mathrm{d}c} \right\rangle - \frac{1}{2}\langle b, \kappa b\rangle + \left\langle c - b, \log \frac{1 - \frac{\mathrm{d}b}{\mathrm{d}c}}{\kappa(c-b)} \right\rangle + \frac{1}{2}\langle c, \kappa c\rangle$$

$$= \left\langle b, \log \frac{\frac{\mathrm{d}b}{\mathrm{d}c}}{\mathrm{e}^{-\kappa(c-b)}} \right\rangle + \left\langle c - b, \log \frac{1 - \frac{\mathrm{d}b}{\mathrm{d}c}}{\kappa(c-b)\mathrm{e}^{-\frac{1}{2}\kappa(c-b)}} \right\rangle.$$

We claim that $\kappa(c-b)\mathrm{e}^{-\frac{1}{2}\kappa(c-b)} \leq 1 - \mathrm{e}^{-\kappa(c-b)}$ holds pointwise. Indeed, the function $\psi(z) := 1 - \mathrm{e}^{-z} - z\mathrm{e}^{-\frac{z}{2}}, z \geq 0$, satisfies that $\psi(0) = 0$ and $\psi'(z) = \mathrm{e}^{-\frac{z}{2}}(\mathrm{e}^{-\frac{z}{2}} - (1 - \frac{z}{2})) \geq 0$ for any $z \geq 0$, implying that $\psi(z) \geq 0$ for any $z \geq 0$. So the claim holds, since $\kappa(c-b) \geq 0$ holds pointwise. Therefore, we can estimate

$$G_c(b) - G_c(c) \geq \left\langle b, \log \frac{\frac{\mathrm{d}b}{\mathrm{d}c}}{\mathrm{e}^{-\kappa(c-b)}} \right\rangle + \left\langle c - b, \log \frac{1 - \frac{\mathrm{d}b}{\mathrm{d}c}}{1 - \mathrm{e}^{-\kappa(c-b)}} \right\rangle$$

$$= \int_\mathcal{S} c(\mathrm{d}r) \Big[ \frac{\mathrm{d}b}{\mathrm{d}c}(r) \log \frac{\frac{\mathrm{d}b}{\mathrm{d}c}(r)}{\mathrm{e}^{-\kappa(c-b)(r)}} $$

$$+ \left(1 - \frac{db}{dc}(r)\right) \log \frac{1 - \frac{db}{dc}(r)}{1 - e^{-\kappa(c-b)(r)}} \right] \geq 0,$$

with equality if and only if $b = c$. The last inequality can be seen by applying Jensen's inequality to the function $x \mapsto x \log x$ or by noting the following: For any point $r \in \mathcal{S}$ the term in brackets in the last line is an entropy between the Bernoulli distribution with (success) parameter $\frac{db}{dc}(r)$ and the Bernoulli distribution with parameter $e^{-\kappa(c-b)(r)}$ and therefore non-negative.

(2) Define $F_c$ as the generalized analog of (6.28), i.e., for $b \leq c$

$$F_c(b) = \left\langle c, \log \frac{db}{d\mu} \right\rangle + (c-b)(\mathcal{S}) + \frac{1}{2} \langle b, \kappa b \rangle - \langle c, \kappa b \rangle + \frac{1}{2} \langle c, \kappa \mu \rangle.$$

Let $b, b' \leq c$ with $\Sigma(\kappa, b) \leq 1$ and $\Sigma(\kappa, b') \leq 1$. We want to show that $F_c(b') \leq G_c(b)$. By rearranging terms one can see that

$$G_c(b) - F_c(b') = \mathbb{H}(c-b|\kappa(c-b)b') + \mathbb{H}(b|b') - \frac{1}{2}\langle b-b', \kappa(b-b')\rangle. \quad (7.9)$$

Now, given the signed measure $b - b'$ we use the Hahn decomposition theorem to decompose $\mathcal{S}$ into two disjoint sets $S_+, S_-$ with $S_+ \cup S_- = \mathcal{S}$ such that $\delta_+(\cdot) := (b-b')(\cdot \cap S_+)$ and $-\delta_-(\cdot) := -(b-b')(\cdot \cap S_-)$ are non-negative measures and $b - b' = \delta_+ - \delta_-$. Observe that

$$\frac{1}{2}\langle b-b', \kappa(b-b')\rangle = \frac{1}{2}\int_{S_+} (b-b')(ds)\,\kappa(b-b')(s)$$
$$+ \frac{1}{2}\int_{S_-} (b'-b)(ds)\,\kappa(b'-b)(s).$$

We write $f_{b',b} := \frac{db'}{db}\mathbb{1}_{S_+}$ and denote by $\langle \cdot, \cdot \rangle_b$ the inner product on $L^2(b)$, i.e. $\langle f, g \rangle_b = \int f(s)g(s)\,b(ds)$. Note that by the symmetry of $\kappa$ we have $\langle f, T_{\kappa,b}g\rangle_b = \langle g, T_{\kappa,b}f\rangle_b$, so we have that $\Sigma(\kappa, b) = \sup_{f \neq 0} \frac{\langle f, T_{\kappa,b}f\rangle_b}{\langle f, f\rangle_b} \leq 1$. Then

$$\frac{1}{2}\int_{S_+} (b-b')(ds)\,\kappa(b-b')(s) = \frac{1}{2}\langle \mathbb{1}_{S_+} - f_{b',b}, T_{\kappa,b}(\mathbb{1}_{S_+} - f_{b',b})\rangle_b$$
$$\leq \frac{1}{2}\Sigma(\kappa, b)\langle \mathbb{1}_{S_+} - f_{b',b}, \mathbb{1}_{S_+} - f_{b',b}\rangle_b \leq \frac{1}{2}\int_{S_+} b(ds)\,(1 - f_{b',b}(s))^2. \quad (7.10)$$

An elementary analysis shows that $\frac{1}{2}(1-x)^2 \leq -\log x + x - 1$ for $x \in (0, 1]$ with equality if and only if $x = 1$, and since $b' \leq b$ on $S_+$ implies that $f_{b',b}(s) \in (0, 1]$ for $s \in S_+$, we get that

$$\frac{1}{2}\int_{S_+} b(ds)\,(1 - f_{b',b}(s))^2 \leq \int_{S_+} b(ds)\left(-\log f_{b',b}(s) + f_{b',b}(s) - 1\right)$$
(7.11)

Denote $f_{b,b'} := \frac{\mathrm{d}b}{\mathrm{d}b'} \mathbb{1}_{S_-}$. Interchanging the roles of $b$ and $b'$ and replacing $S_+$ by $S_-$ one can argue as in (7.10) to show that

$$\frac{1}{2}\int_{S_-}(b'-b)(\mathrm{d}s)\kappa(b'-b)(s) \leq \frac{1}{2}\int_{S_-} b'(\mathrm{d}s)(1-f_{b,b'}(s))^2.$$

An elementary analysis shows that $\frac{1}{2}(1-x)^2 \leq x\log x + 1 - x$ for $x \in (0,1]$ with equality if and only if $x = 1$, and since $b \leq b'$ on $S_-$ implies that $f_{b,b'}(s) \in [0,1]$ for $s \in S_-$, we get that

$$\frac{1}{2}\int_{S_-} b'(\mathrm{d}s)(1-f_{b,b'}(s))^2 \leq \int_{S_+} b'(\mathrm{d}s)\Big(f_{b,b'}(s)\log f_{b,b'}(s) + 1 - f_{b,b'}(s)\Big). \tag{7.12}$$

Note that the two expressions on the right-hand sides of (7.11) and (7.12) sum up to $\mathbb{H}(b|b')$, hence we have shown that

$$\mathbb{H}(b|b') - \frac{1}{2}\langle b-b', \kappa(b-b')\rangle \geq 0$$

and we have equality if and only if $b = b'$. Using this in Eq. (7.9) and the fact that the first entropy term in (7.9) is always non-negative, we get that

$$\sup_{b \leq c:\ \Sigma(\kappa,b) \leq 1} F_c(b) \leq \inf_{b \leq c:\ \Sigma(\kappa,b) \leq 1} G_c(b).$$

Note that $G_c(b) - F_c(b') = 0$ if and only if the following conditions are satisfied: (i) $b = b'$, (ii) $b$ is a solution of (2.9) and (iii) $\Sigma(\kappa,b) = 1$. By Lemma 6.14 the only choice is given by $b = b' = b_*$, where $b_*$ is the unique minimal solution of (2.9). Hence, the uniqueness claim holds. □

*Proof of Theorem 2.3* The projection $(\lambda,\alpha) \mapsto \lambda$ is continuous with respect to the vague topology, so the contraction principle gives that the LDP for $\mathrm{Mi}_N$ holds with rate function

$$\begin{aligned}\mathcal{I}_{\mathrm{Mi}}(\lambda) &= \inf_{\alpha:\ c(\alpha) \leq \mu - c(\lambda)} I(\lambda,\alpha) \\ &= I_{\mathrm{Mi}}(\lambda) + \inf_{\alpha \in \mathcal{A}, \nu \in \mathcal{M}(\mathcal{S}):\ c(\alpha)+\nu=\mu-c(\lambda)}\Big(I_{\mathrm{Ma}}(\alpha) + I_{\mathrm{Me}}(\nu)\Big)\end{aligned} \tag{7.13}$$

assuming $c(\lambda) \leq \mu$ (the other case is trivial). By Eq. (7.2) of Lemma 7.1 we immediately get the representation for $\mathcal{I}_{\mathrm{Mi}}$ claimed in Eq. (2.6).

The projection $(\lambda, \alpha) \mapsto \alpha$ is continuous with respect to the chosen topology, so the contraction principle gives that the LDP for $\mathrm{Ma}_N$ holds with rate function

$$\mathcal{I}_{\mathrm{Ma}}(\alpha) = \inf_{\lambda:\, c(\lambda) \leq \mu - c(\alpha)} I(\lambda, \alpha)$$
$$= I_{\mathrm{Ma}}(\alpha) + \inf_{c \in [0, \mu - c(\alpha)]} \left( \inf_{\lambda:\, c(\lambda) = c} I_{\mathrm{Mi}}(\lambda) + I_{\mathrm{Me}}(\mu - c(\alpha) - c) \right), \tag{7.14}$$

assuming $c(\alpha) \leq \mu$ (the other case is trivial). Define $\mu_\alpha = \mu - c(\alpha)$.

Now, assume that $\Sigma(\kappa, \mu_\alpha) \leq 1$. Then for any fixed $c \in \mathcal{M}(\mathcal{S})$ with $c \leq \mu_\alpha$ we have that $\Sigma(\kappa, c) \leq 1$, so according to Eq. (6.1) of Proposition 6.1 we have

$$\inf_{\lambda:\, c(\lambda) = c} I_{\mathrm{Mi}}(\lambda) + I_{\mathrm{Me}}(\mu_\alpha - c) = G_{\mu_\alpha}(c)$$

with $G_{\mu_\alpha}(c)$ defined as in (6.29). By Lemma 7.3 we have that $\min_{c \leq \mu_\alpha} G_{\mu_\alpha}(c) = G_{\mu_\alpha}(\mu_\alpha)$, which implies Eq. (2.7) under the assumption $\Sigma(\kappa, \mu_\alpha) \leq 1$.

Now, assume that $\Sigma(\kappa, \mu_\alpha) > 1$. Then by Proposition 6.1 we have

$$\inf_{\lambda:\, c(\lambda) = c} I_{\mathrm{Mi}}(\lambda) + I_{\mathrm{Me}}(\mu_\alpha - c)$$
$$= \begin{cases} G_{\mu_\alpha}(c) & \text{if } \Sigma(\kappa, c) \leq 1, \\ G_c(b^*(c)) + \langle \mu_\alpha - c, \log \frac{\mu_\alpha - c}{\kappa(\mu_\alpha - c)} \rangle & \text{if } \Sigma(\kappa, c) > 1. \end{cases}$$

In the case $\Sigma(\kappa, c) > 1$, one can use the same argument as in (7.5) to show that $G_c(b^*(c)) + \langle \mu_\alpha - c, \log \frac{\mu_\alpha - c}{\kappa(\mu_\alpha - c)} \rangle > G_{\mu_\alpha}(b^*(c))$, where $b^* = b^*(c)$ is given as in (2.9). In particular, $\Sigma(\kappa, b^*(c)) = 1$ holds, so any possible minimizer has to be in the set $\{c: \Sigma(\kappa, c) \leq 1\}$. Now, recall that due to Lemma 7.3

$$\inf_{c:\, \Sigma(\kappa, c) \leq 1} G_{\mu_\alpha}(c) = G_{\mu_\alpha}(b^*(\mu_\alpha))$$

This proves the claim of equation (2.7) under the assumption $\Sigma(\kappa, \mu_\alpha) > 1$. □

### 7.3 The minimizers of $I$

*Proof of Theorem 2.1* Note that $\inf I(\lambda, \alpha) = \inf_{c \in \mathcal{M}(\mathcal{S}):\, c \leq \mu} \mathcal{J}(c)$ where for fixed $c \in \mathcal{M}(\mathcal{S})$ with $c \leq \mu$ we define

$$\mathcal{J}(c) = \inf_{\lambda, \alpha:\, c(\lambda) = c} I(\lambda, \alpha) =: \begin{cases} \mathcal{J}_{\leq 1}(c) & \text{if } \Sigma(\kappa, c) \leq 1 \\ \mathcal{J}_{>1}(c) & \text{otherwise.} \end{cases} \tag{7.15}$$

By Proposition 6.1 and Lemma 7.1 we have that

$$\mathcal{J}_{\leq 1}(c) = \left\langle c, \log \frac{\mathrm{d}c}{\mathrm{d}\mu} \right\rangle + \frac{1}{2} \langle c, \kappa(\mu - c) \rangle + I_{\mathrm{Ma}}(\delta_{\mu-c})$$

$$\mathcal{J}_{>1}(c) \geq \left\langle b^*, \log \frac{\mathrm{d}b^*}{\mathrm{d}\mu} \right\rangle + \frac{1}{2} \langle b^*, \kappa(\mu - b^*) \rangle + I_{\mathrm{Me}}(c - b^*) + I_{\mathrm{Ma}}(\delta_{\mu-c})$$

with $b^* = b^*(c)$ characterized in (2.9).

We will start by minimizing the function $\mathcal{J}_{\leq 1}$ over all $c \in \mathcal{M}(\mathcal{S})$ with $c \leq \mu$. Rearranging terms, we get that

$$\begin{aligned}\mathcal{J}_{\leq 1}(c) &= \int_\mathcal{S} \mu(\mathrm{d}r) \left( \frac{\mathrm{d}c}{\mathrm{d}\mu}(r) \log \frac{\frac{\mathrm{d}c}{\mathrm{d}\mu}(r)}{\mathrm{e}^{-\kappa(\mu-c)(r)}} + \left(1 - \frac{\mathrm{d}c}{\mathrm{d}\mu}(r)\right) \log \frac{1 - \frac{\mathrm{d}c}{\mathrm{d}\mu}(r)}{1 - \mathrm{e}^{-\kappa(\mu-c)(r)}} \right) \\ &= \int_\mathcal{S} \mu(\mathrm{d}r) H(\beta^{(r)} | \gamma^{(r)}),\end{aligned}$$

where for $r \in \mathcal{S}$ we defined $\beta^{(r)}$ and $\gamma^{(r)}$ to be Bernoulli distributions with success rate $\frac{\mathrm{d}c}{\mathrm{d}\mu}(r)$ and $\mathrm{e}^{-\kappa(\mu-c)(r)}$, respectively (note that $c \leq \mu$ implies that $\frac{\mathrm{d}c}{\mathrm{d}\mu}(r) \leq 1$ for all $r \in \mathcal{S}$). By Jensen's inequality we have that $H(\beta^{(r)} | \gamma^{(r)}) \geq 0$ for all choices of $\beta^{(r)}$ and $\gamma^{(r)}$ and $H(\beta^{(r)} | \gamma^{(r)}) = 0$ if and only if $\beta^{(r)} = \gamma^{(r)}$, that is if and only if $\frac{\mathrm{d}c}{\mathrm{d}\mu}(r) = \mathrm{e}^{-\kappa(\mu-c)(r)}$. The minimizer of $\mathcal{J}_{\leq 1}$ is therefore characterized by (2.4).

In the case $\Sigma(\kappa, \mu) \leq 1$ (which implies $\Sigma(\kappa, c) \leq 1$ for all $c \in \mathcal{M}(\mathcal{S})$ with $c \leq \mu$), Lemma 4.1 states that $\mu$ is the only solution to (2.4). Note that even though we formulated Lemma 4.1 only for the case of a finite type space, it is valid also under our general assumptions due to Theorem 6.2 and Theorem 6.7 in [7].

Now assume that $\Sigma(\kappa, \mu) > 1$. Then for any $c \in \mathcal{M}(\mathcal{S})$ with $c \leq \mu$ and $\Sigma(\kappa, c) > 1$, Lemma 7.1 implies that $I_{\mathrm{Me}}(c - b^*) + I_{\mathrm{Ma}}(\delta_{\mu-c}) > I_{\mathrm{Ma}}(\delta_{\mu-b^*})$, where $b^* = b^*(c)$ is given as in (2.9) and satisfies $\Sigma(\kappa, b^*) = 1$. Therefore, $\mathcal{J}_{>1}(c) > \mathcal{J}_{\leq 1}(b^*(c))$, which implies that the minimizer of $\mathcal{J}$ lies in the set $\{c \colon \Sigma(\kappa, c) \leq 1\}$. (Note, that in this way, $\mu$ is ruled out as a minimizer, although it solves Eq. (2.4)). By the analysis of $\mathcal{J}_{\leq 1}$ above, the minimizer of $\mathcal{J}$ is given by a solution to (2.4) satisfying $\Sigma(\kappa, c) \leq 1$. Applying the second part of Lemma 4.1 (again under generalized assumptions) finishes the proof. □

**Remark 7.4** (*Reducibility*) Theorem 2.1 and 2.3 are proved under the assumptions of Theorem 1.1, in particular when $\kappa$ is irreducible with respect to $\mu$. We see however that this condition does not play any role in the minimization of Proposition 6.1 and of Lemma 7.1. It is indeed Theorem 2.6 that excludes the admissibility of a minimizer of the form $(\lambda_{c^*}, \delta_{\mu-c^*})$ when $\Sigma(\kappa, \mu) > 1$ for $\widetilde{I}$, as $\alpha^* = \delta_{\mu-c^*}$ may not be connectable. It is straightforward to see that the optimal macroscopic mass in this case takes the form $\widetilde{\alpha}^* = \sum_n \delta_{y^{(n)}}$, with $y^{(n)}(\cdot) = (\mu - c^*)(\cdot \cap S^{(n)})$, for each irreducible class $S^{(n)}$.

### 7.4 The Flory equation

As we explained in Sect. 2.4 the graph model studied in this paper has an important connection with a certain inhomogeneous coagulation process. In this section we prove that the statistics of the limiting microscopic cluster distribution, i.e., the minimizer of the rate function $I$, satisfy the Flory equation (2.14), the related deterministic PDE, over the entire time interval $[0, \infty)$, before and after the gelation time $t_c = 1/\Sigma(\kappa, \mu)$. We prove it in the case of a finite type set $\mathcal{S}$.

We fix an irreducible symmetric matrix $\kappa$ on the finite type space $\mathcal{S}$. Recall from Proposition 6.2 the explicit formula

$$\lambda_k(c; \kappa) = \tau(k; \kappa) \prod_{r \in \mathcal{S}} \frac{(c_r e^{-(\kappa c)_r})^{k_r}}{k_r!}, \qquad k \in \mathbb{N}_0^{\mathcal{S}}, c \in (0, \infty)^{\mathcal{S}}, \qquad (7.16)$$

for the minimizer of the microscopic rate function $I_{\text{Mi}}$ (see also (6.3)). With this notation we stressed the dependence on $\kappa$, since we consider now $t\kappa$ instead of $\kappa$, where $t \in [0, \infty)$ is a time instant, writing $\lambda(\mu; t\kappa)$. Note that $\lambda(\mu; t\kappa)$ is the minimizer both for $t < t_c$ and for $t \geq t_c$, as the characteristic equation $c_r^*(t) e^{-t(\kappa c^*(t))_r} = \mu_r e^{-t(\kappa \mu)_r}$ (where $c^*(t)$ is the $c^*$ of Proposition 6.2 for $t\kappa$ instead of $\kappa$) ensures that $\lambda(c^*(t); t\kappa) = \lambda(\mu; t\kappa)$. Note that $c(\lambda(c^*(t); t\kappa)) = c^*(t)$. We now show that the function $t \mapsto \lambda(\mu; t\kappa)$ solves the Flory equation that corresponds to our model.

**Lemma 7.5** *The function $t \mapsto \lambda(\mu; t\kappa)$ is a solution of*

$$\frac{d}{dt} l_k(t) = \frac{1}{2} \sum_{m + \widetilde{m} = k} l_m(t) l_{\widetilde{m}}(t) \langle m, \kappa \widetilde{m} \rangle - l_k(t) \sum_m l_m(0) \langle k, \kappa m \rangle, \quad \text{for } k \in \mathbb{N}_0^{\mathcal{S}}, \qquad (7.17)$$

*with initial condition $\lambda(\mu; 0) = \sum_{r \in \mathcal{S}} \mu_r \mathbb{1}_{\{e_r\}}$.*

*Proof* The initial condition is easily checked.

Since any tree contributing to the term $\tau(k; t\kappa)$ has exactly $|k| - 1$ edges, we can rewrite $\lambda(\mu; t\kappa)$ for any $t \geq 0$ as

$$\lambda_k(\mu; t\kappa) = t^{|k|-1} \tau(k; \kappa) e^{-t\langle k, \kappa \mu \rangle} \prod_{r \in \mathcal{S}} \frac{\mu_r^{k_r}}{k_r!}, \quad \text{for } k \in \mathbb{N}_0^{\mathcal{S}}. \qquad (7.18)$$

Abreviating $\lambda(t) := \lambda(\mu; t\kappa)$, we get that

$$\frac{d}{dt} \lambda_k(t) = (|k| - 1) \frac{1}{t} \lambda_k(t) - \langle k, \kappa \mu \rangle \lambda_k(t). \qquad (7.19)$$

Now, we study the first summand of the r.h.s. of (7.17). By first inserting (7.18) and then using the recursive Eq. (6.7) from Lemma 4.3, we have that

$$\frac{1}{2} \sum_{m+\widetilde{m}=k} \lambda_m(t)\lambda_{\widetilde{m}}(t)\langle m, \kappa\widetilde{m}\rangle$$

$$= \frac{1}{2} t^{|k|-2} e^{-t\langle k, \kappa\mu\rangle} \Big(\prod_u \frac{\mu_u^{k_u}}{k_u!}\Big) \sum_{r,s} \kappa(r,s) \sum_{m+\widetilde{m}=k} \Big(\prod_u \frac{k_u!}{m_u!\widetilde{m}_u!}\Big) \tau(m) m_r \tau(\widetilde{m}) \widetilde{m}_s$$

$$= (|k|-1)\frac{1}{t}\lambda_k(t).$$

Furthermore, we have that

$$\lambda_k(t) \sum_m \lambda_m(0)\langle k, \kappa m\rangle = \lambda_k(t)\langle k, \kappa\mu\rangle,$$

which implies the claim. □

**Acknowledgements** This research has been funded by the Deutsche Forschungsgemeinschaft (DFG) through grant CRC 1114 "Scaling Cascades in Complex Systems", Project C08, and Grant SPP2265 "Random Geometric Systems", Project P01.

**Funding** Open Access funding enabled and organized by Projekt DEAL.

**Publisher's Note** Springer Nature remains neutral with regard to jurisdictional claims in published maps and institutional affiliations.